\documentclass[11pt, a4paper, reqno]{amsart}

\usepackage{anysize}
\usepackage[utf8]{inputenc}
\usepackage{amsmath}
\usepackage{pdflscape,adjustbox}
\usepackage{amsthm}
\usepackage{mathtools}
\usepackage{amssymb}
\usepackage{hyperref}
\usepackage{xcolor}
\usepackage{mdframed}
\usepackage{mathrsfs}
\usepackage{enumerate}
\usepackage{complexity}		

\usepackage{tikz}
\usetikzlibrary{patterns.meta,calc}
\usetikzlibrary{decorations.pathreplacing,calligraphy}
\usepackage{bm}
\usepackage{dsfont}						
\usepackage{enumitem}
\usepackage[colorinlistoftodos]{todonotes}

\definecolor{munsell}{rgb}{0.0, 0.5, 0.69}

\marginsize{2.5cm}{2.5cm}{2cm}{2cm}

\newtheorem{theorem}{Theorem}
\newtheorem{conjecture}{Conjecture}
\newtheorem*{question}{Question}
\newtheorem{corollary}[theorem]{Corollary}
\newtheorem{lemma}[theorem]{Lemma}

\theoremstyle{definition}
\newtheorem{definition}{Definition}
\AfterEndEnvironment{definition}{\noindent\ignorespaces}
\AfterEndEnvironment{theorem}{\noindent\ignorespaces}
\AfterEndEnvironment{conjecture}{\noindent\ignorespaces}
\AfterEndEnvironment{lemma}{\noindent\ignorespaces}
\AfterEndEnvironment{proof}{\noindent\ignorespaces}
\AfterEndEnvironment{corollary}{\noindent\ignorespaces}
\AfterEndEnvironment{proposition}{\noindent\ignorespaces}
\newtheorem*{remark}{Remark}

\newcommand{\pred}{\rho_r}
\newcommand{\pblue}{\rho_b}

\begin{document}

\title{Ramsey Numbers of Most Even-Linked Double Stars}

\author{Chaitanya Karamchedu, Mithra Karamchedu, Maria Klawe, Andr\'es R.\ Vindas-Mel\'endez}
\address{Department of Computer Science\\ University of Maryland, College Park}
\email{cdkaram@umd.edu}

\address{Department of Computer Science and Department of Mathematics\\ Harvey Mudd College \vskip 0pt Department of Computer Science\\ Columbia University}
\email{mithra@cs.columbia.edu}

\address{Math for America}
\email{klawe@mathforamerica.org}

\address{Department of Mathematics, Harvey Mudd College, \url{https://math.hmc.edu/arvm/}}
\email{avindasmelendez@g.hmc.edu}

\begin{abstract}
    For $n \geq m$, the \emph{linked double star} $S_c(n, m)$ is the graph formed by taking two disjoint stars $K_{1, n}$ and $K_{1, m}$ and linking their centers with a path on $c$ vertices. In this paper, we investigate the Ramsey numbers of \emph{even-linked double stars}, that is, where $c = 2p$ for some $p \in \mathbb{N}$. Determining the Ramsey numbers of the $p = 1$ case---classical double stars---is a longstanding open question \cite{GrossmanHararyKlawe}; the $p = 2$ case was resolved by Burr and Erd\H{o}s in 1976 \cite{BurrErdos}, and Liu and Chen made progress for general even-linked double stars in 2026 \cite{LiuChen}. We extend  Burr and Erd\H{o}s' formula for $p = 2$, proving that the Ramsey number of $S_{2p}(n, m)$ is given by \[r(S_{2p}(n, m)) = \max(2n + 2p - 1, n + 2m + 3p - 1),\]when $m \geq p -2 \geq 0$. Our upper bound argument hinges upon the existence (or nonexistence) of an object that we call a \emph{polarized pair}. We also establish a matching lower bound that applies for all $n, m, p \geq 2$, and as such, we conjecture that this formula applies for all even-linked double stars where $p \geq 2$. 
\end{abstract}
\maketitle
\section{Introduction}
Given a graph $G$, the \emph{Ramsey number} $r(G)$ is the smallest positive integer $r$ such that every $2$-edge-coloring of the complete graph $K_r$ contains a monochromatic copy of $G$.
Determining Ramsey numbers exactly is well known to be extremely difficult and has attracted the attention of mathematicians since the publication of Erd\H{o}s and Szekeres' seminal 1935 paper \cite{ErdosSzekeres}.
The problem of establishing the Ramsey numbers of various classes of graphs quickly becomes computationally infeasible, even for relatively small instances.
As a result, the Ramsey numbers of most graph families are established through analytical methods rather than exhaustive computation.  
Exact values of Ramsey numbers have been established for various families of graphs, including paths, wheels, books, and stars; for surveys on known results of Ramsey numbers, see \cite{ConlonFoxSudakov, Radziszowski}.

In 1972, Harary determined the Ramsey numbers of stars \cite{Harary}, and in 1979, Grossman (acknowledging the prior work of Rosta and Burr) identified those of any pair of disjoint stars \cite{Grossman}.
In the same year, Grossman, Harary, and Klawe initiated the study of the Ramsey numbers of \emph{double stars}, graphs formed by taking two disjoint stars and adding an edge between their centers \cite{GrossmanHararyKlawe}.
They determined the Ramsey numbers for most such graphs, leaving only a small gap of cases unresolved. 
Progress on the remaining cases has been achieved by Norin, Sun, and Zhao in 2016 \cite{NorinSunZhao}, as well as by Flores Dub\'o and Stein in 2025 \cite{ FloresStein}.

\begin{figure}[ht]
\label{fig:double_star}
\centering
\begin{center}
\begin{tikzpicture}[transform shape,scale = 0.5, every node/.style={scale=0.7}]
    \node[circle,fill,inner sep=0.3cm] (center) at (0,0) {};
\foreach \phi in {1,...,7}{
    \node[circle,fill,inner sep=0.3cm]  (v_\phi) at (360/8  * \phi:3cm) {};
         \draw[black] (v_\phi) -- (center);
      }
    \node[circle,fill,inner sep=0.3cm] (center5) at (15,0) {};

    \begin{scope}[xshift=15cm]
    \foreach \phi in {1,...,5}{
    \node[circle,fill,inner sep=0.3cm]  (w_\phi) at (36 + 360/5  * \phi:3cm) {};
         \draw[black] (w_\phi) -- (center5);
      }
      \end{scope}
       \node[circle,fill,inner sep=0.3cm] (center2) at (3,0) {};
        \node[circle,fill,inner sep=0.3cm] (center3) at (6,0) {};
         \node[circle,fill,inner sep=0.3cm] (center4) at (9,0) {};
    \draw[black] (center) -- (center2) -- (center3) -- (center4) -- (w_2);
   \end{tikzpicture} 
\end{center}
\caption{The even-linked double star $S_6(7,4)$, with Ramsey number $r(S_6(7, 4)) = 23$.}
\end{figure}
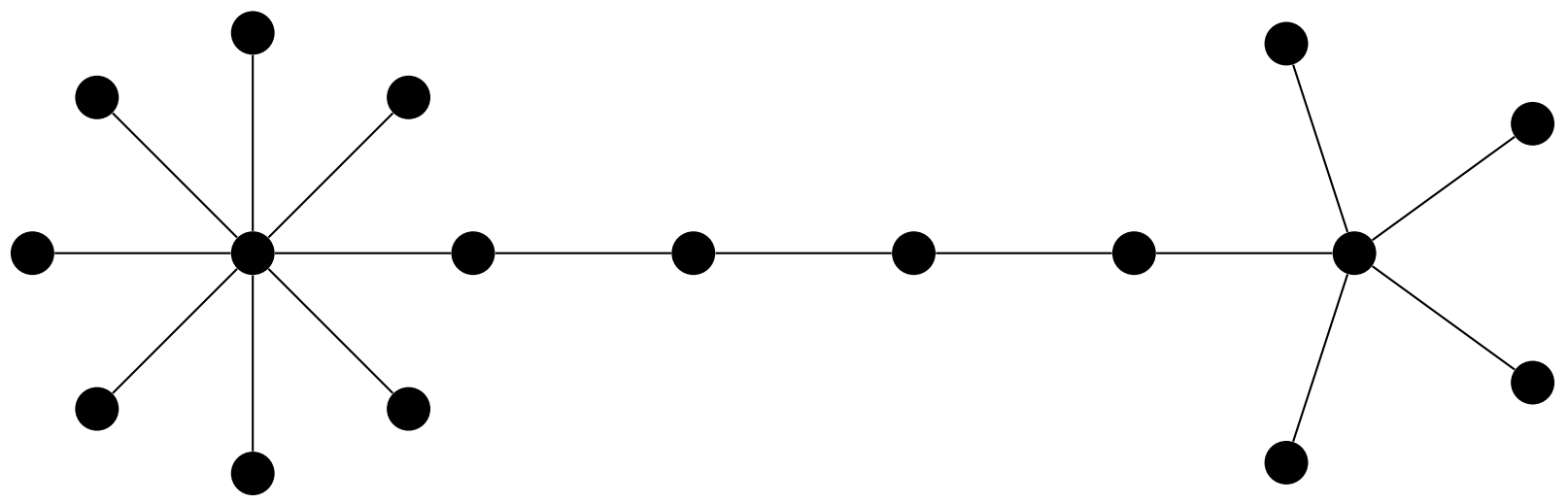

Motivated by the results of Grossman et al.\ and Norin et al., Karamchedu and Klawe in 2022 investigated the Ramsey numbers of a generalization of double stars, known as \emph{linked double stars} \cite{KaramcheduKlawe}.
A linked double star is defined as the union of two disjoint stars with a path joining their centers. 
In particular, for $n \geq m \geq 0$, the linked double star $S_c(n,m)$ is the graph on the vertices $\{ a_1, \ldots, a_c, v_1, \ldots, v_n, w_1, \ldots w_m  \}$ with the edges
\[\big\{ \{a_1, v_i\}, \{a_1, a_2\}, \{a_2, a_3\}, \ldots, \{a_{c-1}, a_c\}, \{a_c, w_j\} : 1 \leq i \leq n, 1 \leq j \leq m\big\}.\]
Note that above, $a_1$ and $a_c$ are the centers of the $n$-star and $m$-star, respectively, and the path from $a_1$ to $a_c$ is the link between the two stars. 
See Figure \ref{fig:double_star} for an example of a linked double star, namely, $S_6(7,4)$.

In \cite{KaramcheduKlawe}, Karamchedu and Klawe studied the Ramsey numbers of \emph{odd-linked double stars}, that is, when $c=2p+1$ for $p\in \mathbb{N}$ in the definition above. 
They determined the Ramsey numbers of two groups of odd-linked double stars, finding that
\begin{enumerate}
    \item $r(S_{2p + 1}(n,2)) = n + 3p + 3$, if $n \leq p-2$. 
    \item $r(S_{2p + 1}(n,m)) =2n + 2m + 2p - 1$, if $n \geq 2p + 1$.
\end{enumerate}
In this paper, we investigate the Ramsey numbers of \emph{even-linked double stars} $S_{2p}(n,m)$, for $p \geq 3$. 
Note that the $p = 2$ case was solved by Burr and Erd\H{o}s in 1976  \cite{BurrErdos}, who proved that 
\[r(S_4(n, m)) = \max(2n + 3, n + 2m + 5).\] 
The case of $p = 1$ is that of ``classical\rq\rq\ double stars, first studied by Grossman et al.\ in \cite{GrossmanHararyKlawe}. 
In 2026, Liu and Chen published the first results for the general even-linked case \cite{LiuChen}, showing that
\begin{enumerate}
    \item $r(S_{2p}(n,2)) = n + 3p + 3$, if $n = 2$, or $n \geq 3$ and $p \geq n + 1$,
    \item $r(S_{2p}(n, m)) \leq 2(n + m + p - 1)$, if $n \geq 2p \geq m$.
\end{enumerate}
In this paper, we prove that
\[r(S_{2p}(n, m)) = \max(2n + 2p - 1, n+2m + 3p -1),\]
when $m \geq p - 2\geq 0$.

\subsection{Outline}

In Section \ref{sec:lower_bounds}, we prove a lower bound on the Ramsey numbers of even-linked double stars, which hold irrespective of the relative sizes of the stars and link.

In Section \ref{sec:upper_bounds}, we establish an upper bound in the regime where $n \geq m \geq p-2$ and $p \geq 4$; this upper bound matches the lower bound described in Section \ref{sec:lower_bounds}, thus solving the corresponding Ramsey numbers exactly. The upper bound argument is organized into two sections. Specifically, we define a structure called a \textit{polarized pair} (Definition \ref{def:polarized_pair}) and then condition on whether our 2-coloring of $K_r$ admits such a pair. In Section \ref{sec:polarized_pair_exists}, we prove that if the 2-coloring has a polarized pair, then it must admit a monochromatic $S_{2p}(n,m)$; Section \ref{sec:no_polarized_pair_exists} shows that we still have a monochromatic $S_{2p}(n,m)$ if the 2-coloring does not have a polarized pair. The arguments in Sections \ref{sec:polarized_pair_exists} and \ref{sec:no_polarized_pair_exists} assume that $p \geq 4$, so we handle the $p=3$ case separately in Appendix \ref{sec:p=3}. In order to help the reader navigate the upper bound arguments, we have included a dependency diagram in Figure \ref{fig:even-linked-dependencies}.

In Section \ref{sec:conclusion}, we conclude by discussing a few open problems suggested by this work.

%

\usetikzlibrary{arrows.meta,calc}

\begin{landscape}
\begin{figure}[p]
\centering
\begingroup
\newcommand{\ELDPanelShift}{1.2}
\rmfamily
\definecolor{eldark}{HTML}{202B3C}
\definecolor{elgray}{HTML}{637183}
\definecolor{elbluepanel}{HTML}{F3F7FF}
\definecolor{elblueborder}{HTML}{D5DEEF}
\definecolor{elbluebox}{HTML}{DFEDFF}
\definecolor{elgreenpanel}{HTML}{F6FCF9}
\definecolor{elgreenborder}{HTML}{D0E6DC}
\definecolor{elgreenbox}{HTML}{EAF8F1}
\definecolor{elamberpanel}{HTML}{FFFBF3}
\definecolor{elamberborder}{HTML}{E4D4B7}
\definecolor{elamberbox}{HTML}{FFF0D2}
\definecolor{elneutral}{HTML}{F3F5F8}
\definecolor{elborder}{HTML}{D4DCE9}

\begin{adjustbox}{max totalsize={\linewidth}{0.94\textheight}}
\begin{tikzpicture}[x=1cm,y=-1cm,
  elbox/.style={rectangle,draw=black,
    line width=1.2pt,minimum height=1.42cm,inner xsep=.18cm,
    inner ysep=.16cm,outer sep=0pt,align=flush center,text=eldark,
    font=\rmfamily\fontsize{12.5}{15}\selectfont},
  elarrow/.style={draw=elgray,line width=1.0pt,
    -{Stealth[length=2.5mm,width=1.7mm]},sharp corners,line join=miter,
    line cap=butt,shorten <=0.6pt,shorten >=0.6pt},
  elarrowupper/.style={elarrow,shorten >=1.0pt},
  eljoinarrow/.style={elarrow,shorten <=0pt},
  elline/.style={draw=elgray,line width=1.0pt,sharp corners,line join=miter,
    line cap=butt,shorten <=0.6pt,shorten >=0pt},
  ellabel/.style={font=\rmfamily\fontsize{9.5}{11}\selectfont,
    text=elgray,inner sep=2pt,align=flush center}]
\useasboundingbox (0,0) rectangle ({27.1+\ELDPanelShift},22.5);

\def\ELDBlueLeftX{3.2}
\def\ELDBlueRightX{9.3}
\def\ELDBlueColumnWidth{4.55}
\pgfmathsetmacro{\ELDBlueCenterX}{(\ELDBlueLeftX+\ELDBlueRightX)/2}
\pgfmathsetmacro{\ELDBlueFullWidth}{\ELDBlueRightX-\ELDBlueLeftX+\ELDBlueColumnWidth}

\def\ELDGreenLeftX{16.65}
\def\ELDGreenRightX{23.55}
\def\ELDGreenColumnWidth{5.1}
\pgfmathsetmacro{\ELDGreenCenterX}{(\ELDGreenLeftX+\ELDGreenRightX)/2}
\pgfmathsetmacro{\ELDGreenFullWidth}{\ELDGreenRightX-\ELDGreenLeftX+\ELDGreenColumnWidth}

\newcommand{\ELDBox}[6]{%
  \node[elbox,text width=#3cm,
    minimum width=\dimexpr#3cm+0.36cm\relax,fill=#4] (#1) at (#2)
    {\textbf{#5}\if\relax\detokenize{#6}\relax\else\\[4pt]{\fontsize{11}{13}\selectfont #6}\fi};%
}

\begin{scope}[yshift=2cm]
\path[fill=elbluepanel,draw=elblueborder,line width=1pt,rounded corners=9pt]
  (0,2.0) rectangle (12.5,20.1);
\node[text=eldark,font=\rmfamily\fontsize{15}{18}\selectfont\bfseries]
  at (6.25,2.65) {Polarized Pair Exists};
\node[text=elgray,font=\rmfamily\fontsize{10.5}{13}\selectfont]
  at (6.25,3.22) {Section \ref{sec:polarized_pair_exists}\quad$\bullet$\quad $p\ge4$};

\ELDBox{l5}{\ELDBlueLeftX,4.65}{\ELDBlueColumnWidth}{elbluebox}{Lemma \ref{lem:helper3}}{}
\ELDBox{l9}{\ELDBlueRightX,4.65}{\ELDBlueColumnWidth}{elbluebox}{Lemma \ref{lem:helper2}}{}
\ELDBox{l6}{\ELDBlueLeftX,9.35}{\ELDBlueColumnWidth}{elbluebox}{Lemma \ref{lem:W_i>=2p-2}}{$|W_R|,|W_B|\ge2p-2$}
\ELDBox{l10}{\ELDBlueRightX,9.35}{\ELDBlueColumnWidth}{elbluebox}{Lemma \ref{lem:helper_red_cycle_case}}{}
\ELDBox{l78}{\ELDBlueLeftX,14.10}{\ELDBlueColumnWidth}{elbluebox}{Lemmas \ref{lem:smallW_B}--\ref{lem:n=m=p-2}}{$|W_B|\le p-3$}
\ELDBox{l12}{\ELDBlueRightX,14.10}{\ELDBlueColumnWidth}{elbluebox}{Lemma \ref{lem:medium_W_B}}{$p-2\le |W_B|\le2p-3$}
\ELDBox{l13}{\ELDBlueCenterX,18.85}{\ELDBlueFullWidth}{elbluebox}{Lemma \ref{lem:polarized_pair}: Polarized-Pair Case}{A monochromatic $S_{2p}(n,m)$ exists when $m\ge p-2$.}

\draw[elarrow] (l5.south) -- (l6.north);
\draw[elarrow] (l9.south) -- (l10.north);
\draw[elarrow] (l5.east) -- (\ELDBlueCenterX,4.65) -- (\ELDBlueCenterX,14.10) -- (l12.west);
\draw[elarrow] (l10.south) -- (l12.north);
\draw[elline] (l6.south) -- (\ELDBlueLeftX,10.45) -- (.35,10.45)
  -- (.35,17.15) -- (\ELDBlueRightX,17.15);
\draw[elline] (l78.south) -- (\ELDBlueLeftX,17.15);
\draw[elline] (l12.south) -- (\ELDBlueRightX,17.15);
\fill[elgray] (\ELDBlueCenterX,17.15) circle (1.35pt);
\draw[eljoinarrow] (\ELDBlueCenterX,17.15) -- (l13.north);

\begin{scope}[xshift=\ELDPanelShift cm]
\path[fill=elgreenpanel,draw=elgreenborder,line width=1pt,rounded corners=9pt]
  (13.1,2.0) rectangle (27.1,20.1);
\node[text=eldark,font=\rmfamily\fontsize{15}{18}\selectfont\bfseries]
  at (20.1,2.65) {No Polarized Pair Exists};
\node[text=elgray,font=\rmfamily\fontsize{10.5}{13}\selectfont]
  at (20.1,3.22) {Section \ref{sec:no_polarized_pair_exists}\quad$\bullet$\quad $p\ge4$};

\ELDBox{tools}{\ELDGreenCenterX,4.65}{\ELDGreenFullWidth}{elgreenbox}{Lemmas \ref{lem:min_red_blue_degree}--\ref{lem:max_blue_degree}; Corollary \ref{cor:max_blue_degree}}{}
\ELDBox{l1819}{\ELDGreenLeftX,7.0}{\ELDGreenColumnWidth}{elgreenbox}{Lemmas \ref{lem:p_leq_k_leq_m_plus_p_minus_1}--\ref{lem:k_geq_m_plus_p}}{$k\ge p$}
\ELDBox{l21}{\ELDGreenRightX,7.0}{\ELDGreenColumnWidth}{elgreenbox}{Lemma \ref{lem:wrsl_helper}}{}
\ELDBox{l22}{\ELDGreenRightX,9.35}{\ELDGreenColumnWidth}{elgreenbox}{Lemma \ref{lem:wrsl}}{}
\ELDBox{l2325}{\ELDGreenRightX,11.70}{\ELDGreenColumnWidth}{elgreenbox}{Lemmas \ref{lem:k_1_pminus1}--\ref{lem:k_1_pminus1_edge_case}}{$0\le k\le p-1$}
\ELDBox{l26}{\ELDGreenRightX,14.10}{\ELDGreenColumnWidth}{elgreenbox}{Lemma \ref{lem:wrapper2}}{}
\ELDBox{l28}{\ELDGreenLeftX,14.10}{\ELDGreenColumnWidth}{elgreenbox}{Lemma \ref{lem:little_helper}}{}
\ELDBox{l27}{\ELDGreenRightX,16.45}{\ELDGreenColumnWidth}{elgreenbox}{Lemma \ref{lem:large_n_wrapper}}{$n\ge2m+p-1$, for all $k$}
\ELDBox{l29}{\ELDGreenCenterX,18.85}{\ELDGreenFullWidth}{elgreenbox}{Lemma \ref{lem:no_polarized_pair}: No-Polarized-Pair Case}{A monochromatic $S_{2p}(n,m)$ exists when $m\ge p-3$.}

\draw[elarrow] (l1819.north |- tools.south) -- (l1819.north);
\draw[elarrow] (l21.north |- tools.south) -- (l21.north);
\draw[elarrow] (l21.east) -- (26.70,7.0) -- (26.70,11.70) -- (l2325.east);
\draw[elarrow] (l22.south) -- (l2325.north);
\draw[elarrow] (l2325.south) -- (l26.north);
\draw[elarrow] (l26.south) -- (l27.north);
\draw[elarrow] (l1819.east) -- (20.15,7.0) -- (20.15,16.45) -- (l27.west);
\draw[elarrow] (l27.south) -- (l27.south |- l29.north);
\draw[elarrow] (l28.south) -- (l28.south |- l29.north);

\end{scope}
\end{scope}

\node[elbox,text width=7.2cm,fill=elamberbox,minimum height=2.10cm,
  inner ysep=.28cm,anchor=east] (p3)
  at ([xshift=-.1pt]{27.1+\ELDPanelShift},20.60)
  {\textbf{Appendix \ref{sec:p=3}: $p=3$}\\[6pt]
   {\fontsize{11}{13}\selectfont Upper bound for $S_6(n,m)$}};
\draw[elarrow] (p3.north) -- (p3.north |- l29.south);
\node[ellabel,anchor=west] at ($(p3.north |- 0,18.40)+(.10,0)$) {induction base};

\node[elbox,fill=white,draw=black,line width=2pt,text=black,text width=16cm,
  minimum height=3.00cm,anchor=west] (upper)
  at ([xshift=.5pt]0,20.60)
  {\textbf{\fontsize{16}{19}\selectfont Upper Bound (Lemma \ref{lem:upper_bound} for $p\ge4$)}\\[6pt]
   {\fontsize{14}{17}\selectfont
   $r(S_{2p}(n,m))\le\max(2n+2p-1,\,n+2m+3p-1)$}\\[4pt]
   {\fontsize{11}{13}\selectfont $n\ge m\ge p-2\ge0$}};
\draw[elarrowupper] (l13.south) -- (l13.south |- upper.north);
\draw[elarrowupper] (l29.south) -- (l29.south |- 0,18.35) -- (upper.north |- 0,18.35) -- (upper.north);
\draw[elarrowupper] (p3.west) -- (upper.east);

\end{tikzpicture}
\end{adjustbox}
\endgroup

\caption{Dependencies for the upper bound on $r(S_{2p}(n,m))$.}
\label{fig:even-linked-dependencies}
\end{figure}
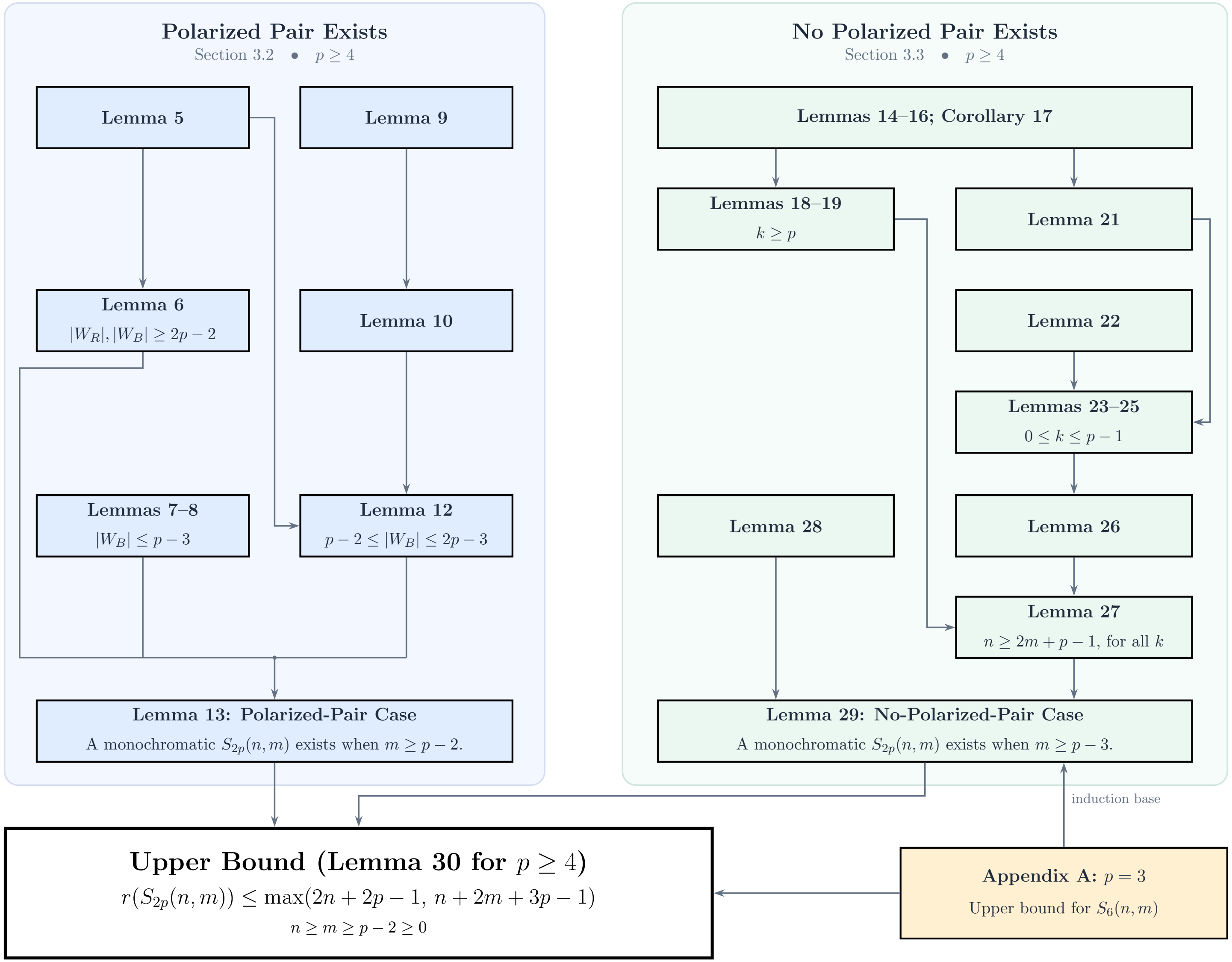
\end{landscape}

\subsection{AI Usage Disclosure} 
We used GPT-5.5 Pro to prove Lemmas \ref{lem:wrsl_helper}--\ref{lem:k_equals_0}, which together represent the case where both $n \geq 2m + p$ and $k \leq p - 1$ in the no-polarized-pair-exists section.  We then took inspiration from the proof of Lemma \ref{lem:k_1_pminus1} to produce Lemma \ref{lem:k_1_pminus1_edge_case}---which relaxes the $n \geq 2m + p$ assumption---without the use of AI. Prior to consulting the model, we suspected that the proof strategy for the $k \leq p - 1$ case would adapt elements of Yu and Li's argument for the Ramsey numbers of brooms (see \cite{YuLi}). We thus prompted GPT-5.5 Pro by asking it to draw inspiration from Yu and Li's techniques, which it did successfully.

As such, there are similarities between Yu and Li's proof and our $k \leq p - 1$ case that are worth commenting on. Both arguments make use of extremal numbers and, specifically, resolved cases of the Erd\H{o}s--S\'os conjecture on trees.\footnote{In the time since we completed our proof, a pre-release version of GPT-6 Astra actually proved the entire Erd\H{o}s--S\'os conjecture, for all trees (as part of the FrontierMath Erd\H{o}s benchmark) \cite{AdamczewskiBloom, astra}. However, our argument here only makes use of the case of spider graphs (which include brooms), proved by Fan, Hong, and Liu in 2019 \cite{Spiders}; see Lemma \ref{Extremal bund on brooms}.} In particular, Yu and Li use the path graph case as a helper lemma (Lemma 8 in their paper, citing Erd\H{o}s and Gallai), whereas we utilize the case of broom graphs (see Lemma \ref{Extremal bund on brooms}). In addition, Lemma \ref{lem:wrsl} here resembles another helper lemma that Yu and Li use to bound the number of edges in a particular bipartite graph (Lemma 9 in their paper, citing Gy\'arf\'as, Rousseau, and Schelp). Yu and Li use these lemmas to show that their $2$-coloring must admit a monochromatic broom graph; if we assume otherwise, then there exists a certain edge quantity whose upper and lower bounds conflict, a contradiction. Lemmas \ref{lem:k_1_pminus1} and \ref{lem:k_1_pminus1_edge_case} in this paper apply a similar contradiction strategy, replacing Yu and Li's helper lemmas with our analogous ones.

Our original argument handled even-linked double stars $S_{2p}(n, m)$ where $m \geq p - 1$. However, GPT-6 Astra highlighted a few simple modifications that extend this result to $m \geq p - 2$. Specifically, the model suggested that we isolate Lemma \ref{lem:helper3}---which was previously part of Lemma \ref{lem:smallW_B}'s argument---as its own statement, so that we could reuse it in the strengthened version of Lemma \ref{lem:medium_W_B}. It also produced the auxiliary Lemma \ref{lem:n=m=p-2} to handle the $n = m = p -2$ case for one particular graph regime.

All of the other lemmas were proved without the use of AI. We also used GPT-5.5 Pro and GPT-6 Astra to audit the full argument for mistakes and simplifications, but we ultimately stand by the correctness of this paper. In addition, we used GPT-5.5 Pro and GPT-6 Astra to help produce Figures \ref{fig:even-linked-dependencies} and \ref{fig:WR-WB}--\ref{fig:appendix_fig} from our illustrations; however, none of the writing in this paper was generated by AI.

\section{Lower Bounds}
\label{sec:lower_bounds}
In this section we establish lower bounds on the Ramsey numbers of the even-linked double stars. 
These lower bounds are derived from \cite[Lemma 1]{Burr_survey}.

\begin{lemma}\label{Burr Lower Bound}
Let $G = X \sqcup Y$ be a connected bipartite graph on $g$ vertices, such that $|X| = k$, $|Y| = \ell$, and $k \geq \ell$. Then,
\[r(G) \geq \max(2k - 1, g + \ell - 1).\]
\end{lemma}

\begin{proof}
We will first prove that $r(G) \geq 2k - 1$ by constructing a $2$-coloring of $K_{2k - 2}$ that avoids a monochromatic copy of $G$. Specifically, let $K_{k-1} \sqcup K_{k-1}$ be the red subgraph of our $2$-coloring, so that the blue subgraph is the complete bipartite graph $K_{k - 1, k - 1}$ between the two disjoint red $K_{k - 1}$ subgraphs. This $2$-coloring does not admit a red copy of $G$, since $G$ is connected and requires at least $g = k + \ell > k - 1$ vertices; we also avoid a blue copy of $G$ since neither partition has the required $k$ vertices needed for $X$. As a result, $r(G) \geq 2k - 1$.\\

\begin{figure}[h]
\centering
\begin{tikzpicture}[scale=0.5]


\foreach \x in {1,...,11}{%
    \pgfmathsetmacro{\ang}{180-(\x-1)*360/11}
    \coordinate (N-\x) at (\ang:5.4cm);
}

\begin{scope}[shift={(16,0)}]
    \foreach \x in {1,...,11}{%
        \pgfmathsetmacro{\ang}{(\x-1)*360/11}
        \coordinate (R-\x) at (\ang:5.4cm);
    }
\end{scope}

\foreach \y in {1,...,11}{%
    \foreach \z in {1,...,11}{%
        \draw[cyan!30, line width=0.25pt]
            (R-\y) -- (N-\z);
    }
}


\foreach \x in {2,...,11}{%
    \foreach \y in {2,...,11}{%
        \ifnum\x<\y
            \draw[red, thick] (N-\x) -- (N-\y);
        \fi
    }
}

\foreach \x in {2,...,11}{%
    \foreach \y in {2,...,11}{%
        \ifnum\x<\y
            \draw[red, thick] (R-\x) -- (R-\y);
        \fi
    }
}


\foreach \y in {2,...,11}{%
    \draw[red, very thick, loosely dotted]
        (N-1) -- (N-\y);
}

\foreach \y in {2,...,11}{%
    \draw[red, very thick, loosely dotted]
        (R-1) -- (R-\y);
}


\foreach \x in {2,...,11}{%
    \fill[black] (N-\x) circle (0.2cm);
}

\foreach \x in {2,...,11}{%
    \fill[black] (R-\x) circle (0.2cm);
}

\draw[black, fill=white, line width=0.6pt]
    (N-1) circle (0.2cm);

\draw[black, fill=white, line width=0.6pt]
    (R-1) circle (0.2cm);


\draw [
    pen colour=red,
    red,
    very thick,
    decorate,
    decoration={
        calligraphic brace,
        amplitude=7.5pt,
        mirror
    }
]
    (-4,-6) -- (4,-6)
    node[pos=0.5, anchor=north, yshift=-0.35cm]
    {Red $K_{k-1}$};

\draw [
    pen colour=blue,
    blue,
    very thick,
    decorate,
    decoration={
        calligraphic brace,
        amplitude=7.5pt,
        mirror
    }
]
    (5,-6) -- (11,-6)
    node[pos=0.5, anchor=north, yshift=-0.35cm]
    {Blue $K_{k-1,k-1}$};

\draw [
    pen colour=red,
    red,
    very thick,
    decorate,
    decoration={
        calligraphic brace,
        amplitude=7.5pt,
        mirror
    }
]
    (12,-6) -- (20,-6)
    node[pos=0.5, anchor=north, yshift=-0.35cm]
    {Red $K_{k-1}$};

\end{tikzpicture}

\caption{The first lower bound construction of Lemma
\ref{Burr Lower Bound}. The second construction replaces the two red $K_{k - 1}$ subgraphs with a red $K_{g - 1}$ and a red $K_{\ell - 1}$.}
\label{fig:lowerbound}
\end{figure}
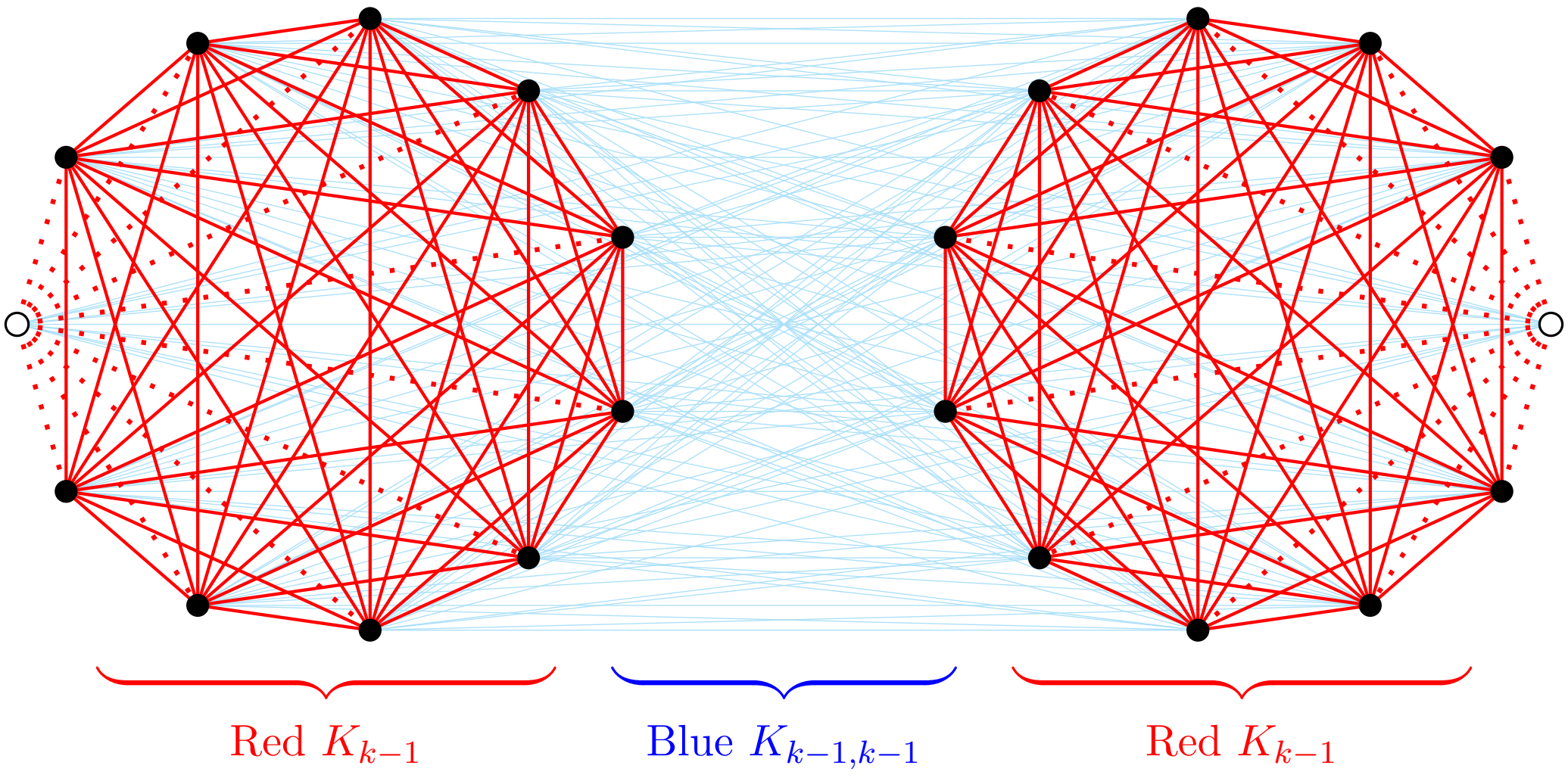
Next, we will prove that $r(G) \geq g + \ell - 1$, by constructing a $2$-coloring of $K_{g + \ell - 2}$ that avoids a monochromatic copy of $G$. 
Let $K_{g-1} \sqcup K_{\ell-1}$ be the red subgraph of our $2$-coloring, so that the blue subgraph is a complete bipartite graph $K_{g - 1,\ell - 1}$ between the two red subgraphs.
This $2$-coloring does not admit a red copy of $G$, since $G$ is connected and requires at least $g$ vertices.
Similarly, there cannot be a blue copy of $G$ since the smaller partition of the blue subgraph does not have the required $\ell$ vertices. Hence, $r(G) \geq g + \ell - 1$, so we conclude that
\[r(G) \geq \max(2k - 1, g + \ell - 1),\]
as desired.
\end{proof}



We thus obtain the following lower bound for the Ramsey numbers of even-linked double stars. 

\begin{lemma}\label{lem:lower_bounds}
The Ramsey numbers of the even-linked double stars satisfy
\[r(S_{2p}(n,m)) \geq \max(2n + 2p - 1,n + 2m + 3p - 1) =
\begin{cases}
2n + 2p - 1, & \text{for } n \geq 2m + p, \\
n+2m+3p-1,  & \text{for } n \leq 2m + p.
\end{cases}
\]
\end{lemma}
\begin{proof}
    This follows immediately from Lemma \ref{Burr Lower Bound}.
\end{proof}


\section{Upper Bounds}
\label{sec:upper_bounds}
For any even-linked double star $S_{2p}(n, m)$, we define \[r = \max(2n + 2p - 1, n + 2m + 3p - 1).\]

The structure of this section is as follows. In Section \ref{sec:polarized_pair_exists}, we show that if a $2$-coloring of $K_r$ contains an $n$-polarized pair, then it admits a monochromatic $S_{2p}(n, m)$. In Section \ref{sec:no_polarized_pair_exists}, we show that any $2$-coloring of $K_r$ still admits a monochromatic $S_{2p}(n, m)$ even if it does not contain an $n$-polarized pair. These two cases imply Lemma \ref{lem:upper_bound}, which gives us our desired upper bound.

\subsection{Preliminaries}

\begin{definition}
    Given a red/blue edge-colored graph $G$ and a vertex $v \in V(G)$, let $R(v)$ and $B(v)$ denote the red and blue neighborhoods of $v$, respectively. Define $$\text{red-}d(v) = |R(v)| \quad \text{and} \quad \text{blue-}d(v) = |B(v)|.$$ For a particular subset of vertices $W \subset V(G)$, we also define $R_W(v) = R(v) \cap W$ and $B_W(v) = B(v) \cap W$; we then have that  $\text{red-}d_W(v) = |R_W(v)|$ and $\text{blue-}d_W(v) = |B_W(v)|$.
\end{definition}

\begin{definition}\label{def:polarized_pair}
    Given a red/blue edge-colored graph $G$ and an integer $n$, an $n$-\emph{polarized pair} is a pair of vertices $\pred, \pblue \in V(G)$ with $|R(\pred) \cap B(\pblue)| \geq n + 1$. Since $n$ is clear in context (i.e.\ the size of the $n$-star), we simply refer to $\pred$ and $\pblue$ as a \emph{polarized pair} for the remainder of the paper. Note that the letter ``$r$\rq\rq\ in the subscript of ``$\pred$\rq\rq\ does not refer to the variable $r$ defined above, but rather the color red.
\end{definition}

\begin{definition}
    The \emph{star graph} $K_{1, n}$ is the order-$(n + 1)$ graph formed by taking a central vertex and connecting it to $n$ ``spoke\rq\rq\ vertices. We also refer to the star graph $K_{1, n}$ as an \emph{$n$-star}.
\end{definition}

\begin{definition}
    The \emph{broom graph} $S_t(n, 0)$ is the graph formed by taking an $n$-star and connecting its center to one of the endpoints of the path graph $P_{t - 1}$. Informally, a broom graph is a double star with $m = 0$ (or equivalently, $m = 1$, with path length one smaller), hence the notation. Note that the broom $S_1(n, 0)$ is simply the $n$-star $K_{1, n}$.
\end{definition}

\begin{remark}
    The study of the Ramsey numbers of brooms was initiated by Burr and Erd\H{o}s, and fully resolved by Yu and Li \cite{BurrErdos, YuLi}. Since broom graphs represent the $m = 0$ and $m = 1$ cases, we assume that $m \geq 2$ throughout the paper.
\end{remark}

\begin{remark}
    Appendix \ref{sec:p=3} handles the $p = 3$ case for all values of $n$ and $m$. Therefore, in all of the following lemmas, we assume that $p \geq 4$ without explicitly stating so.
\end{remark}
\subsection{Polarized Pair Exists}
\label{sec:polarized_pair_exists}
In this section, we will show that any $2$-coloring of $K_r$ admits a monochromatic $S_{2p}(n, m)$ if it contains an $n$-polarized pair. As such, consider any fixed $2$-coloring of $K_r$ with polarized pair $\{\pred, \pblue\}$.
\subsubsection{Notation}
\begin{figure}[h]
\centering
\begin{tikzpicture}[
    line cap=round,
    outline/.style={black,line width=0.85pt},
    rededge/.style={red,line width=0.53pt},
    blueedge/.style={blue,line width=0.53pt},
    redfan/.style={red,line width=0.64pt},
    bluefan/.style={blue,line width=0.64pt},
    dotteddivider/.style={
        black!60,
        dash pattern=on 3pt off 4pt,
        line width=0.8pt
    },
    every node/.style={font=\large}
]

\def\Wleft{0.88}
\def\Wright{8.32}
\def\Wtop{2.9}

\def\WcutL{3.16}
\def\WcutR{6.04}

\def\Abot{4.45}
\def\Atop{6.55}
\def\Aleft{2.2}
\def\Aright{7.0}


\def\spacing{0.48}
\def\shift{1.20}


\foreach \i in {0,...,10} {

    \pgfmathsetmacro{\ax}{\Aleft+\spacing*\i}
    \pgfmathsetmacro{\redwx}{\ax-\shift}
    \pgfmathsetmacro{\bluewx}{\ax+\shift}

    \draw[rededge]
        ({\redwx},\Wtop)
        --
        ({\ax},\Abot);

    \draw[blueedge]
        ({\bluewx},\Wtop)
        --
        ({\ax},\Abot);
}

\coordinate (pR) at (2.85,8.55);
\coordinate (pB) at (6.35,8.55);

\foreach \i in {0,...,7} {
    \draw[redfan]
        (pR) -- ({\Aleft+0.25+0.62*\i},\Atop);
}

\foreach \i in {0,...,7} {
    \draw[bluefan]
        (pB) -- ({\Aleft+0.25+0.62*\i},\Atop);
}

\draw[outline,ultra thick,fill=white]
    (\Aleft,\Abot) rectangle (\Aright,\Atop);

\node at (4.6,5.5) {$|A|=n+1$};

\draw[outline,ultra thick,fill=white]
    (\Wleft,0) rectangle (\Wright,\Wtop);

\draw[dotteddivider]
    (\WcutL,0.045) -- (\WcutL,\Wtop);

\draw[dotteddivider]
    (\WcutR,0.045) -- (\WcutR,\Wtop);

\node at (2.0,1.45) {$W'$};
\node[right] at ({\Wright+0.25},1.45) {$W$};

\draw[outline,fill=black] (pR) circle (2.0pt);
\draw[outline,fill=black] (pB) circle (2.0pt);

\node[above=6pt] at (pR) {$\pred$};
\node[above=6pt] at (pB) {$\pblue$};


\draw[pen colour = {red},
    red,
    very thick,
    decorate,
    decoration={
        calligraphic brace,
        mirror,
        amplitude=13pt,
        raise=4pt,
        aspect=0.35
    },
]
    (\Wleft,-0.05) -- (\WcutR,-0.05)
    node[pos=0.35,below=16pt,red] {$W_R$};

\draw[pen colour = {blue},
    blue,
    very thick,
    decorate,
    decoration={
        calligraphic brace,
        mirror,
        amplitude=13pt,
        raise=4pt,
        aspect=0.65
    },
]
    (\WcutL,-0.5) -- (\Wright,-0.5)
    node[pos=0.65,below=16pt,blue] {$W_B$};

\end{tikzpicture}

\caption{The sets $A$, $W$, $W_R$, and $W_B$,
and the connections between them.}
\label{fig:WR-WB}
\end{figure}
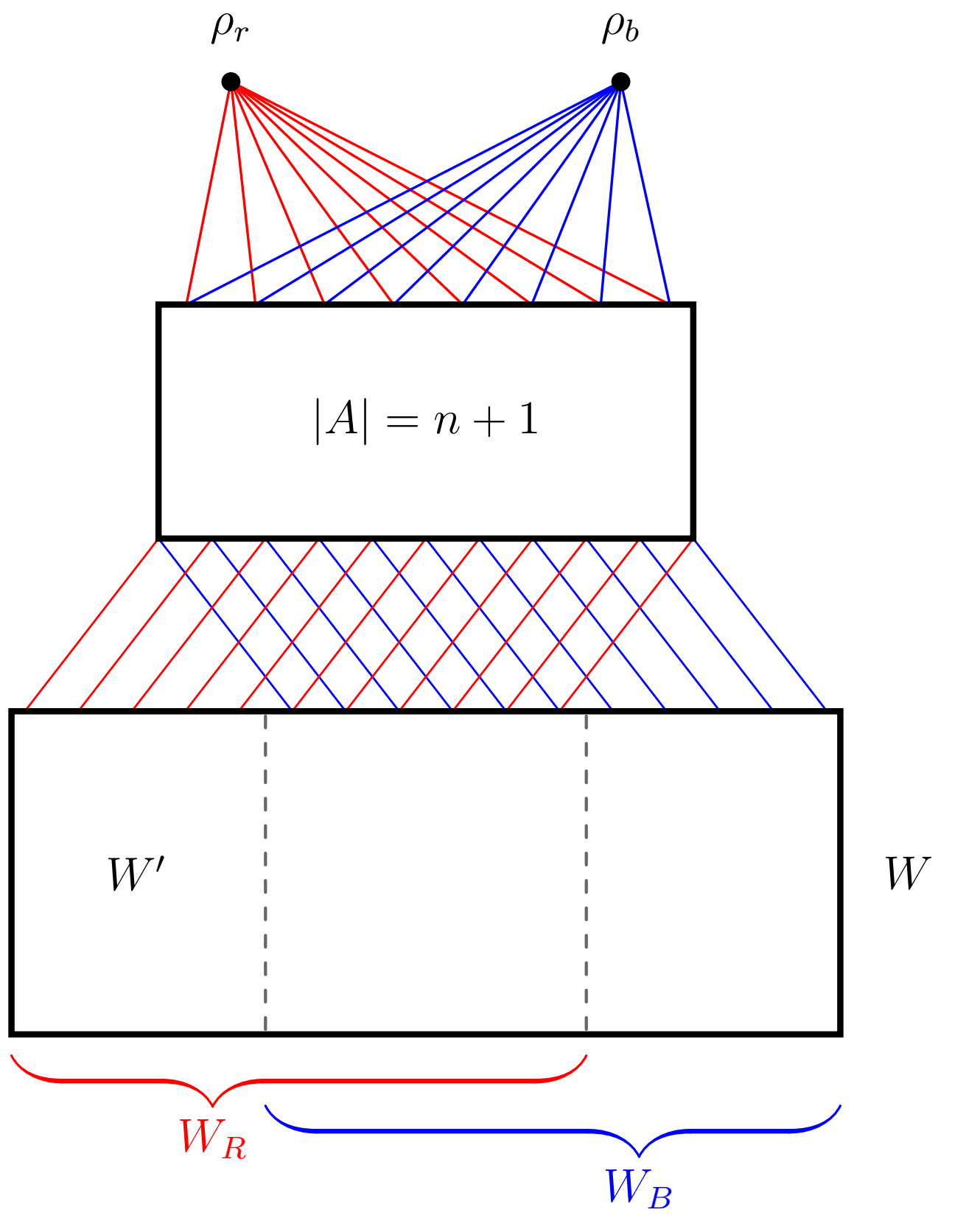

We define the following sets with respect to our chosen $2$-coloring; we will reuse them throughout the section. See Figure \ref{fig:WR-WB} for these sets and the relationships between them.
\begin{itemize}
    \item Let $A$ be any subset of $n + 1$ vertices in $R(\pred) \cap B(\pblue)$. Note that we fix $|A| = n + 1$ for convenience, even though the full red/blue neighborhood shared by $\pred$ and $\pblue$ may be larger.
    \item Let $W = V(K_r)\setminus (A \cup \{\pred, \pblue\})$, the set of vertices not in $A \cup \{\pred, \pblue\}$.
    \item  Let $W_R = \big\{w \in W : \{a, w\}\text{ is red for some }a \in A\big\}$. Similarly, define $W_B = \big\{w \in W : \{a, w\}\text{ is blue for some }a \in A\big\}$, and observe that $W = W_R \cup W_B$.
    \item  Finally, let $W^\prime = W \setminus W_B$, the set of vertices connected to $A$ only by red edges.
\end{itemize}

We proceed by first noting two useful results.
\begin{lemma}[Due to Yan and Peng, \cite{Yan_and_Peng}]\label{bipartite-ramsey-cycles}
Let $\ell$ be any integer such that $\ell \geq 2p$. Then, any $2$-coloring of the edges of a $K_{\ell, \ell}$ admits a monochromatic cycle on $2p$ vertices.
\end{lemma}

\begin{lemma}[Due to Jackson, \cite{Jackson_bipartite_cycles}]\label{Jackson_cycles}
    Let $k \geq 1$, and let $G$ be a simple bipartite graph with bipartition $A \sqcup B$, where $|A| = a \geq 2$, $|B| = b \geq k$, and each vertex in $A$ has degree at least $k$. 
    If $ a \leq k$ and $b \leq 2k -2$, then $G$ contains a cycle on $2a$ vertices.
\end{lemma}

\begin{lemma}
    \label{lem:helper3}
     Suppose there exist disjoint sets $R,S, T \subset W$ with $|R| = |S| = |T| = p-1$, such that
    \begin{enumerate}
        \item $R$ and $S$ are both red-complete to $A$,
        \item $T$ is blue-complete to $A$,
        \item There is a red cycle $C$ on $2p-2$ vertices between $S$ and $T$,
        \item There is a blue cycle $C'$ on $2p-2$ vertices between $R$ and $T$.
    \end{enumerate} 
    If our $2$-coloring of $K_r$ has a polarized pair $\{\pred, \pblue\}$, then it admits a monochromatic $S_{2p}(n, m)$. (By symmetry, the roles of red and blue can be interchanged.)
\end{lemma}
\begin{proof}
     Figure \ref{fig:fig5} shows the sets $A$, $R$, $S$, and $T$ with the conditions described above.
    
    To prove this result, we will show that our 2-coloring must have a series of properties in order to avoid a monochromatic $S_{2p}(n, m)$; but we can then show that these properties together also imply the existence of a monochromatic $S_{2p}(n, m)$.
    
    Let $Q = W \setminus (C \cup C^\prime)$, and note that $\vert Q\vert = |W| - 3(p - 1) \geq 2m - 1$. If $\text{red-}d_{Q}(t) \geq m$ for some $t \in T$, then we complete a red $S_{2p}(n,m)$ as follows: use $t$ as the center of our $m$-star (with spokes in $Q$), and follow $C$ for $2p-2$ vertices from $t$ (including $t$), ending at some vertex $s \in S$; connect $s \rightarrow a \rightarrow \pred$, with $\pred$ as the center of our $n$-star (with spokes in $A \setminus \{a\}$).
    Therefore, we may assume that every vertex $t \in T$ satisfies $\text{blue-}d_Q(t) \geq m$.
    
    Suppose that some vertex $r \in R$ has $\text{blue-}d_Q(r) \geq m$. 
    We then complete a blue $S_{2p}(n, m)$ as follows: use $r$ as the center of our $m$-star (with spokes in $Q$), and follow the cycle $C^\prime$ from $r$ for $2p - 2$ vertices (including $r$), ending at some vertex $t \in T$; connect $t \rightarrow a \rightarrow \pblue$, using $\pblue$ as the center of our $n$-star (with spokes in $A \setminus \{a\}$).
    Hence, we may assume that every vertex $r \in R$ has $\text{blue-}d_Q(r) \leq m - 1$ and thus $\text{red-}d_Q(r) \geq m$.

    Next, suppose that there is a red edge between some vertex $r \in R$ and $s \in S$.
    Then, we can form a red $S_{2p}(n, m)$ as follows: use $r$ as the center of our $m$-star (with spokes in $Q$), and connect $r$ to $s$; follow $C$ from $s$ through $2p-3$ vertices (including $s$), ending at some $s' \in S$; finally, connect $s'\rightarrow a \rightarrow \pred$, with $\pred$ as the center of our $n$-star (with spokes in $A \setminus \{a\}$). 
    As a result, we may assume that there is a blue complete bipartite graph between $R$ and $S$.

    Similarly, suppose that there is a blue edge between some vertex $s \in S$ and $t \in T$.
    We can then complete a blue $S_{2p}(n, m)$ as follows: use some other vertex $t' \in T$ as the center of our $m$-star (with spokes in $Q$), and connect $t'$ to $r$ for some $r \in R$; next, use the blue complete bipartite graph between $R$ and $S$ to complete a blue path on $2p-4$ vertices between $R$ and $S$ that ends at $s$; finally, connect $ s \rightarrow t \rightarrow a \rightarrow \pblue$, using $\pblue$ as the center of our $n$-star (with spokes in $A \setminus \{a\}$). 
    Such a vertex $a$ must again exist by our prior reasoning. 
    Hence, we may assume that there is a red complete bipartite graph between $S$ and $T$.

    Now, suppose that there exists some blue edge $\{t_1, t_2\}$ inside $T$'s induced subgraph.
    Then, we can complete a blue $S_{2p}(n,m)$ as follows: use $t_1$ as our blue $m$-star (with spokes in $Q$), and connect $t_1 \rightarrow t_2 \rightarrow r$ for some $r \in R$; using the blue complete bipartite graph between $R$ and $S$, create a blue path on $2p-5$ vertices starting at $r$, ending at some other vertex $r' \in R$ such that $\{r', t'\}$ is blue for some $t' \neq t_1, t_2$; connect $r' \rightarrow t' \rightarrow a \rightarrow \pblue$ (we know that this is possible, since $C'$ is a blue cycle between $R$ and $T$); finally, use $\pblue$ as the center of our $n$-star (with spokes in $A \setminus \{a\}$). 
    Thus, we may assume that there are no blue edges anywhere inside $T$, so $T$ defines a red complete graph.

    Likewise, suppose that there exists some blue edge $\{s_1, s_2\}$ inside the subgraph induced by $S$. We can then complete a blue $S_{2p}(n,m)$ as follows: use some vertex $t \in T$ as our $m$-star (with spokes in $Q$); connect $t$ to $r$, for some $r \in R$; using the blue complete bipartite graph between $R$ and $S$, complete a blue path on $2p-6$ vertices that starts at $r$ and ends at $s_1$, but does not use $s_2$; then append $s_1 \rightarrow s_2 \rightarrow r' \rightarrow t' \rightarrow a \rightarrow \pblue$ to the end of this path, choosing $r' \neq r$ and $t' \neq t$; use $\pblue$ as the center of our $n$-star (with spokes in $A \setminus \{a\}$).
    Hence, we may assume that there are no blue edges inside $S$, and so $S$ also defines a red complete graph. Because $T$ is also a red complete graph and there is a red complete bipartite graph between $S$ and $T$, the set $S \cup T$ induces a red complete graph on $2p-2$ vertices.

    Suppose now that there is a red edge between $\pred$ and some vertex $s \in S$.
    We then form a red $S_{2p}(n, m)$ as follows: use $\pred$ as the center of our $n$-star (with spokes in $A$), connect $\pred$ to $s$, and follow the cycle $C$ from $s$ for $2p - 3$ vertices (including $s$), ending at some vertex $s^\prime \in S$; connect $s^\prime$ to some $a \in A$, and connect $a$ to some $r \in R$, using $r$ as the center of our $m$-star (with spokes in $Q$).
    Hence, we may assume that every edge from $\pred$ to $S$ is blue.

     Likewise, suppose that there is a blue edge between $\pred$ and some vertex $r \in R$. Since $r \in C^\prime$, we have that $r$ is blue-adjacent to some vertex $t \in T$. We then complete a blue $S_{2p}(n, m)$ as follows: use $t$ as the center of our $m$-star (with spokes in $Q$), and connect $t \rightarrow r \rightarrow \pred$; next, connect $\pred$ to some $s \in S$, and cycle through the blue complete bipartite graph between $S$ and $R$ for $2p - 6$ vertices (including $s$), ending at some vertex $r^\prime \in R$; connect $r^\prime$ to its neighbor $t^\prime \in T$ via $C^\prime$, and connect $t^\prime \rightarrow a \rightarrow \pblue$, using $\pblue$ as the center of our $n$-star.
    As a result, we may assume that every edge from $\pred$ to $R$ is red.

    But we can now complete a red $S_{2p}(n, m)$. 
    Let $\{a, s\}$ be any red edge between $A$ and $S$. 
    Then, choose any vertex $r \in R$ to be the center of our $m$-star (with spokes in $Q$), and connect $r \rightarrow a \rightarrow s$; using the red complete graph on $S \cup T$, form a red path on $2p - 4$ vertices from $s$ to some other vertex $s^\prime \in S$; connect $s^\prime$ to some $a^\prime \in A$ (with $a' \neq a$), and connect $a'$ to $\pred$. 
    Since $\pred$ has $n + 1$ red neighbors in $A$ and $p - 1$ in $R$---of which we have only exhausted $a, a^\prime \in A$ and $r \in R$---we can still use $\pred$ as the center of our $n$-star.
    This finishes a red $S_{2p}(n, m)$, as desired.
\end{proof}

\begin{lemma}
    \label{lem:W_i>=2p-2}
    Suppose $|W_R|, |W_B|\geq2p-2$ and $m \geq \lceil p/2\rceil$. If our $2$-coloring of $K_r$ has a polarized pair $\{\pred, \pblue\}$, then it admits a monochromatic $S_{2p}(n,m)$.
\end{lemma}

\begin{proof}
    We will first show that we can find disjoint sets $R, S, T \subset W$ with $|R| = |S| = |T| = p - 1$ satisfying the four conditions of Lemma \ref{lem:helper3}.
    
    If $m \geq \lceil p/2\rceil$, then $|W| \geq 2m + 3p - 4 \geq 4p - 4$. 
    Since $|W_R| \geq 2p-2$, $|W_B| \geq 2p-2$, and $|W|\geq 4p-4$, we can identify two disjoint sets $Z_R \subseteq W_R$ and $Z_B \subseteq W_B$, both of size $2p-2$: create $Z_R$ by taking as many vertices as needed from $W_R \setminus W_B$; if more vertices are necessary to reach size $2p - 2$, choose the remainder from $W_R \cap W_B$.
    If we can form $Z_R$ entirely from $W_R \setminus W_B$, then we can take $Z_B$ to be any size-$(2p - 2)$ subset of $W_B$. 
    If not, then we must exhaust the vertices in $W_R \setminus W_B$ to form $Z_R$, so $W \setminus Z_R \subseteq W_B$ (i.e.\ the remaining vertices all belong to $W_B$).
    Since $|W \setminus Z_R| = |W| - (2p - 2) \geq 2p - 2$, we can choose $Z_B$ as any size-$(2p - 2)$ subset of $W \setminus Z_R$. 

      By Lemma \ref{bipartite-ramsey-cycles}, there exists a monochromatic cycle on $2p-2$ vertices in the bipartite graph between $Z_R$ and $Z_B$; let $Q$ be such a cycle.
      If $Q$ is red, we use the vertices of $Q \cap Z_R$ to connect back to $A$.
      Similarly, if $Q$ is blue, we use $Q \cap Z_B$.

     Without loss of generality, let us assume that $Q$ is red. 
     Let $X = Q \cap Z_R$ and $Y = Q \cap Z_B$ be the halves of $Q$'s vertices that belong to $Z_R$ and $Z_B$, respectively; note that $\vert X \vert = \vert Y\vert = p - 1$. 
     If any vertex $y \in Y$ has red degree $\text{red-}d_{W \setminus Q}(y) \geq m$, then we can complete a red $S_{2p}(n,m)$ as follows: use $y$ as the center of our $m$-star (with spokes in $W \setminus Q$), and follow $Q$ back from $y$ for $2p-2$ vertices (including $y$); since $y \in Y$, we will end at some vertex $x \in X$; finally, connect $x \rightarrow a \rightarrow \pred$ and complete the $n$-star at $\pred$ (with spokes in $A \setminus \{a\}$). 
     Note that such a vertex $a \in A$ with $\{x,a\}$ being red must exist, since $X \subset Z_R \subseteq  W_R$. 
     Hence, we may assume that $\text{red-}d_{W\setminus Q}(y) \leq m -1$ for all $y \in Y$. 
     Since $|W \setminus Q| = \max(n - 2, 2m + p - 2)$, it follows that
    \[\text{blue-}d_{W \setminus Q}(y) \geq |W\setminus Q| - (m-1) = \max(n - m - 1, m + p - 1) \geq |Y|\] for all $y \in Y$. 
    
    Further note that \[|W \setminus Q| \leq 2\cdot \text{blue-}d_{W \setminus Q}(y) -2\] for all $y \in Y$. 
    As a result, Lemma \ref{Jackson_cycles} with $a=p-1, b= |W\setminus Q|$, and $k= |W \setminus Q| -m + 1$ gives us a blue cycle on $2p-2$ vertices between $Y$ and $W \setminus Q$; call this cycle $Q'$, and let $S = (W \setminus Q) \cap Q'$.
    Because every $y \in Y$ has $\text{blue-}d_{W \setminus Q}(y) \geq m + p -1$, every vertex $q \in Q' \cap Y$ defines a blue $S_{2p-2}(m,0)$: we use $q$ as the center of our $m$-star (with spokes in $W \setminus (Q \cup Q^\prime)$), follow the cycle $Q^\prime$ for $2p - 2$ vertices (including $q$), ending with the tail vertex in $S$. 
    As such, if there is a blue edge between some vertex $s \in S$ and $a \in A$, then we can complete a blue $S_{2p}(n,m)$ as follows: use $\pblue$ as the $n$-star (with spokes in $A \setminus \{a\}$), connect $\pblue \rightarrow a \rightarrow s$, and then follow the cycle $Q^\prime$ from $s$ for $2p - 2$  vertices (including $s$), ending at some vertex $q \in Q^\prime \cap Y$; use $q$ as the center of our $m$-star (with spokes in $W \setminus (Q \cup Q^\prime)$).
    Consequently, we may assume that there is a red complete bipartite graph between $S$ and $A$. Figure \ref{fig:fig2} below shows the relationships between these sets and subgraphs.
    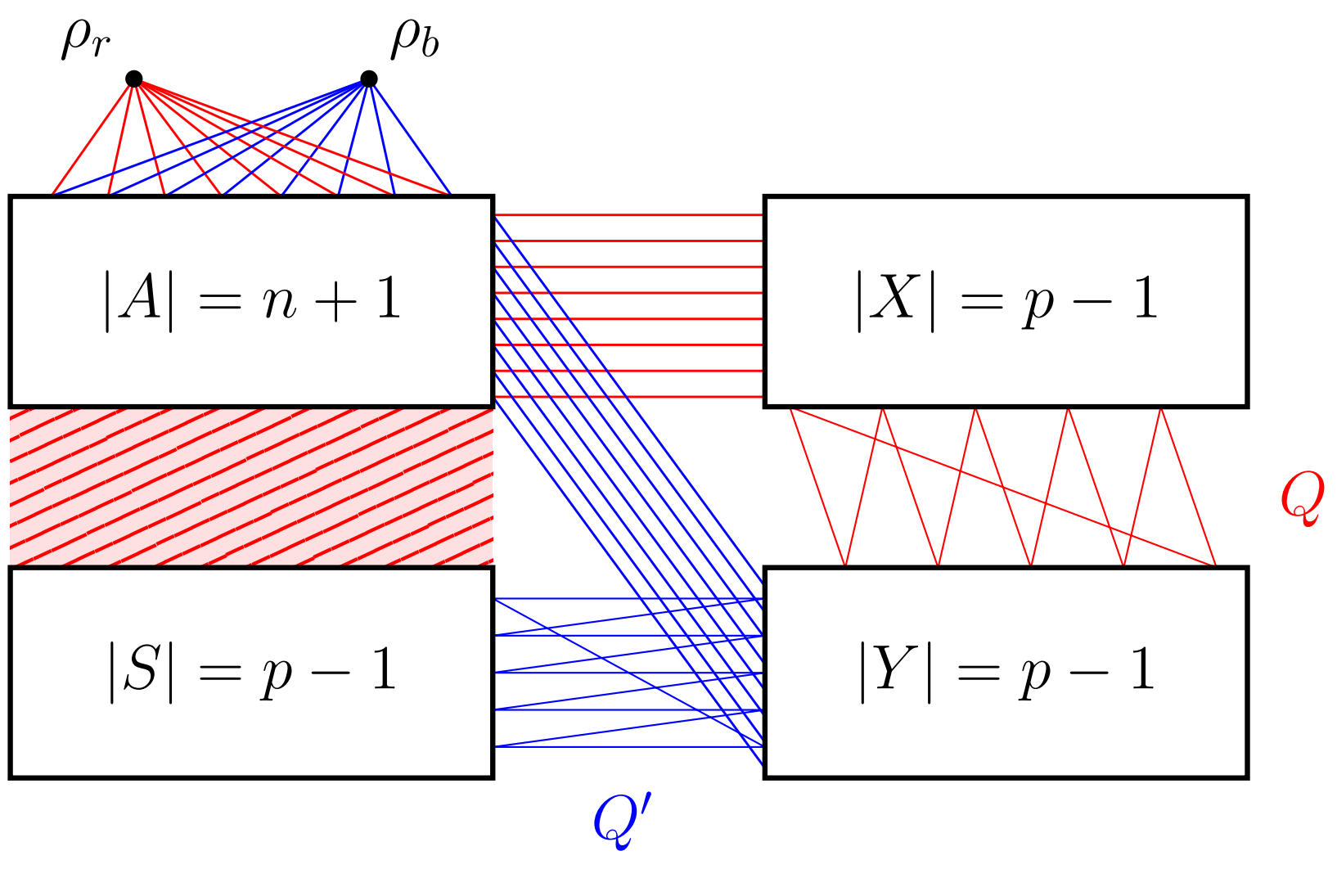
\begin{figure}[h]
    \centering
    \begin{tikzpicture}[
    box/.style={draw=black, very thick, fill=white},
    redfan/.style={red, line width=0.55pt},
    bluefan/.style={blue, line width=0.55pt},
    redcycle/.style={red},
    bluecycle/.style={blue},
    redhatch/.style={
        pattern={Lines[angle=25,distance=4.5pt,line width=.8pt]},
        pattern color=red
    },
    bluehatch/.style={
        pattern={Lines[angle=-25,distance=4.5pt,line width=.8pt]},
        pattern color=blue
    },
    every node/.style={font=\Large}
]

\def\W{3.9}
\def\H{1.7}
\def\Xgap{2.2}
\def\Ygap{1.3}

\pgfmathsetmacro{\Axc}{\W/2}
\pgfmathsetmacro{\Atopy}{2*\H+\Ygap}

\coordinate (Rsw) at (0,0);
\coordinate (Rne) at ({\W},{\H});

\coordinate (Asw) at (0,{\H+\Ygap});
\coordinate (Ane) at ({\W},{2*\H+\Ygap});

\coordinate (Tsw) at ({\W+\Xgap},0);
\coordinate (Tne) at ({2*\W+\Xgap},{\H});

\coordinate (Ssw) at ({\W+\Xgap},{\H+\Ygap});
\coordinate (Sne) at ({2*\W+\Xgap},{2*\H+\Ygap});


\path[fill=red!12] (0,{\H}) rectangle ({\W},{\H+\Ygap});
\path[redhatch]    (0,{\H}) rectangle ({\W},{\H+\Ygap});

\coordinate (ax1) at ({\W},{\H+\Ygap+0.08});
\coordinate (ax2) at ({\W},{\H+\Ygap+0.29});
\coordinate (ax3) at ({\W},{\H+\Ygap+0.50});
\coordinate (ax4) at ({\W},{\H+\Ygap+0.71});
\coordinate (ax5) at ({\W},{\H+\Ygap+0.92});
\coordinate (ax6) at ({\W},{\H+\Ygap+1.13});
\coordinate (ax7) at ({\W},{\H+\Ygap+1.34});
\coordinate (ax8) at ({\W},{\H+\Ygap+1.55});

\coordinate (xa1) at ({\W+\Xgap},{\H+\Ygap+0.08});
\coordinate (xa2) at ({\W+\Xgap},{\H+\Ygap+0.29});
\coordinate (xa3) at ({\W+\Xgap},{\H+\Ygap+0.50});
\coordinate (xa4) at ({\W+\Xgap},{\H+\Ygap+0.71});
\coordinate (xa5) at ({\W+\Xgap},{\H+\Ygap+0.92});
\coordinate (xa6) at ({\W+\Xgap},{\H+\Ygap+1.13});
\coordinate (xa7) at ({\W+\Xgap},{\H+\Ygap+1.34});
\coordinate (xa8) at ({\W+\Xgap},{\H+\Ygap+1.55});

\draw[redfan] (ax1) -- (xa1);
\draw[redfan] (ax2) -- (xa2);
\draw[redfan] (ax3) -- (xa3);
\draw[redfan] (ax4) -- (xa4);
\draw[redfan] (ax5) -- (xa5);
\draw[redfan] (ax6) -- (xa6);
\draw[redfan] (ax7) -- (xa7);
\draw[redfan] (ax8) -- (xa8);

\coordinate (ay1) at ({\W},{\H+\Ygap+0.08});
\coordinate (ay2) at ({\W},{\H+\Ygap+0.29});
\coordinate (ay3) at ({\W},{\H+\Ygap+0.50});
\coordinate (ay4) at ({\W},{\H+\Ygap+0.71});
\coordinate (ay5) at ({\W},{\H+\Ygap+0.92});
\coordinate (ay6) at ({\W},{\H+\Ygap+1.13});
\coordinate (ay7) at ({\W},{\H+\Ygap+1.34});
\coordinate (ay8) at ({\W},{\H+\Ygap+1.55});

\coordinate (ya1) at ({\W+\Xgap},{0.08});
\coordinate (ya2) at ({\W+\Xgap},{0.29});
\coordinate (ya3) at ({\W+\Xgap},{0.50});
\coordinate (ya4) at ({\W+\Xgap},{0.71});
\coordinate (ya5) at ({\W+\Xgap},{0.92});
\coordinate (ya6) at ({\W+\Xgap},{1.13});
\coordinate (ya7) at ({\W+\Xgap},{1.34});
\coordinate (ya8) at ({\W+\Xgap},{1.55});

\draw[bluefan] (ay1) -- (ya1);
\draw[bluefan] (ay2) -- (ya2);
\draw[bluefan] (ay3) -- (ya3);
\draw[bluefan] (ay4) -- (ya4);
\draw[bluefan] (ay5) -- (ya5);
\draw[bluefan] (ay6) -- (ya6);
\draw[bluefan] (ay7) -- (ya7);
\draw[bluefan] (ay8) -- (ya8);

\def\fansep{0.95}
\def\fanygap{0.95}

\coordinate (u)  at ({\Axc-\fansep},{\Atopy+\fanygap});
\coordinate (up) at ({\Axc+\fansep},{\Atopy+\fanygap});

\foreach \dx in {-1.62,-1.16,-0.70,-0.24,0.24,0.70,1.16,1.62}{
    \pgfmathsetmacro{\rx}{\Axc+\dx}
    \pgfmathsetmacro{\bx}{\Axc-\dx}
    \draw[redfan]  (u)  -- ({\rx},{\Atopy});
    \draw[bluefan] (up) -- ({\bx},{\Atopy});
}

\coordinate (s1) at ({\W+\Xgap+0.20},{\H+\Ygap});
\coordinate (s2) at ({\W+\Xgap+0.95},{\H+\Ygap});
\coordinate (s3) at ({\W+\Xgap+1.70},{\H+\Ygap});
\coordinate (s4) at ({\W+\Xgap+2.45},{\H+\Ygap});
\coordinate (s5) at ({\W+\Xgap+3.20},{\H+\Ygap});

\coordinate (t1) at ({\W+\Xgap+0.65},{\H});
\coordinate (t2) at ({\W+\Xgap+1.40},{\H});
\coordinate (t3) at ({\W+\Xgap+2.15},{\H});
\coordinate (t4) at ({\W+\Xgap+2.90},{\H});
\coordinate (t5) at ({\W+\Xgap+3.65},{\H});

\draw[redcycle]
    (s1)--(t1)--(s2)--(t2)--(s3)--(t3)--(s4)--(t4)--(s5)--(t5)--cycle;

\coordinate (r1) at ({\W},1.45);
\coordinate (r2) at ({\W},1.15);
\coordinate (r3) at ({\W},0.85);
\coordinate (r4) at ({\W},0.55);
\coordinate (r5) at ({\W},0.25);

\coordinate (q1) at ({\W+\Xgap},1.45);
\coordinate (q2) at ({\W+\Xgap},1.15);
\coordinate (q3) at ({\W+\Xgap},0.85);
\coordinate (q4) at ({\W+\Xgap},0.55);
\coordinate (q5) at ({\W+\Xgap},0.25);

\draw[bluecycle]
    (r1)--(q1)--(r2)--(q2)--(r3)--(q3)--(r4)--(q4)--(r5)--(q5)--cycle;

\draw[box] (Asw) rectangle (Ane);
\draw[box] (Rsw) rectangle (Rne);
\draw[box] (Ssw) rectangle (Sne);
\draw[box] (Tsw) rectangle (Tne);

\node at ($(Asw)!0.5!(Ane)$) {$|A| = n+1$};
\node at ($(Rsw)!0.5!(Rne)$) {$|S| = p-1$};
\node at ($(Ssw)!0.5!(Sne)$) {$|X| = p-1$};
\node at ($(Tsw)!0.5!(Tne)$) {$|Y| = p-1$};

\fill (u) circle (2.0pt);
\fill (up) circle (2.0pt);

\node[above left=1pt] at (u) {$\pred$};
\node[above right=1pt] at (up) {$\pblue$};

\node[red] at ({2*\W+\Xgap+0.45},{\H+0.55}) {$Q$};
\node[blue] at ({\W+1.05},-0.35) {$Q^\prime$};

\end{tikzpicture}
        \caption{The sets $A$, $S$, $X$, and $Y$, and the subgraphs between them.}
        \label{fig:fig2}
\end{figure}

    Observe that the set $Y$ and cycle $Q^\prime$ have the same relationship to $A$ as the set $X$ and cycle $Q$, only with reversed colors.
    As a result, we may disregard the set $X$ and repeat the preceding argument with red and blue interchanged. Suppose there exists a vertex $s \in S$ with $\text{blue-}d_{W\setminus Q'}(s) \geq m$.
    We then complete a blue $S_{2p}(n,m)$ as follows: use $s$ as the center of our $m$-star (with spokes in $W \setminus Q^\prime$), follow the cycle $Q'$ for $2p-2$ vertices from $s$ (including $s$), ending at some vertex $y \in Y$; connect $y \rightarrow a \rightarrow \pblue$, using $\pblue$ as our $n$-star (with spokes in $A \setminus \{a\}$).
    As before, note that such a vertex $a \in A$ with $\{y, a\}$ being blue must exist, since $Y \subset Z_B \subseteq W_B$.
    Therefore, we may assume that $\text{blue-}d_{W \setminus Q'}(s) \leq m-1$.
    By mirroring our previous argument (which applies Lemma \ref{Jackson_cycles}), we find that there exists a \textit{red} cycle $C$ between $S$ and some set $T \subset W \setminus Q'$, with a blue complete bipartite graph between $T$ and $A$.
    The resulting configuration is illustrated in Figure \ref{fig:fig3}.

    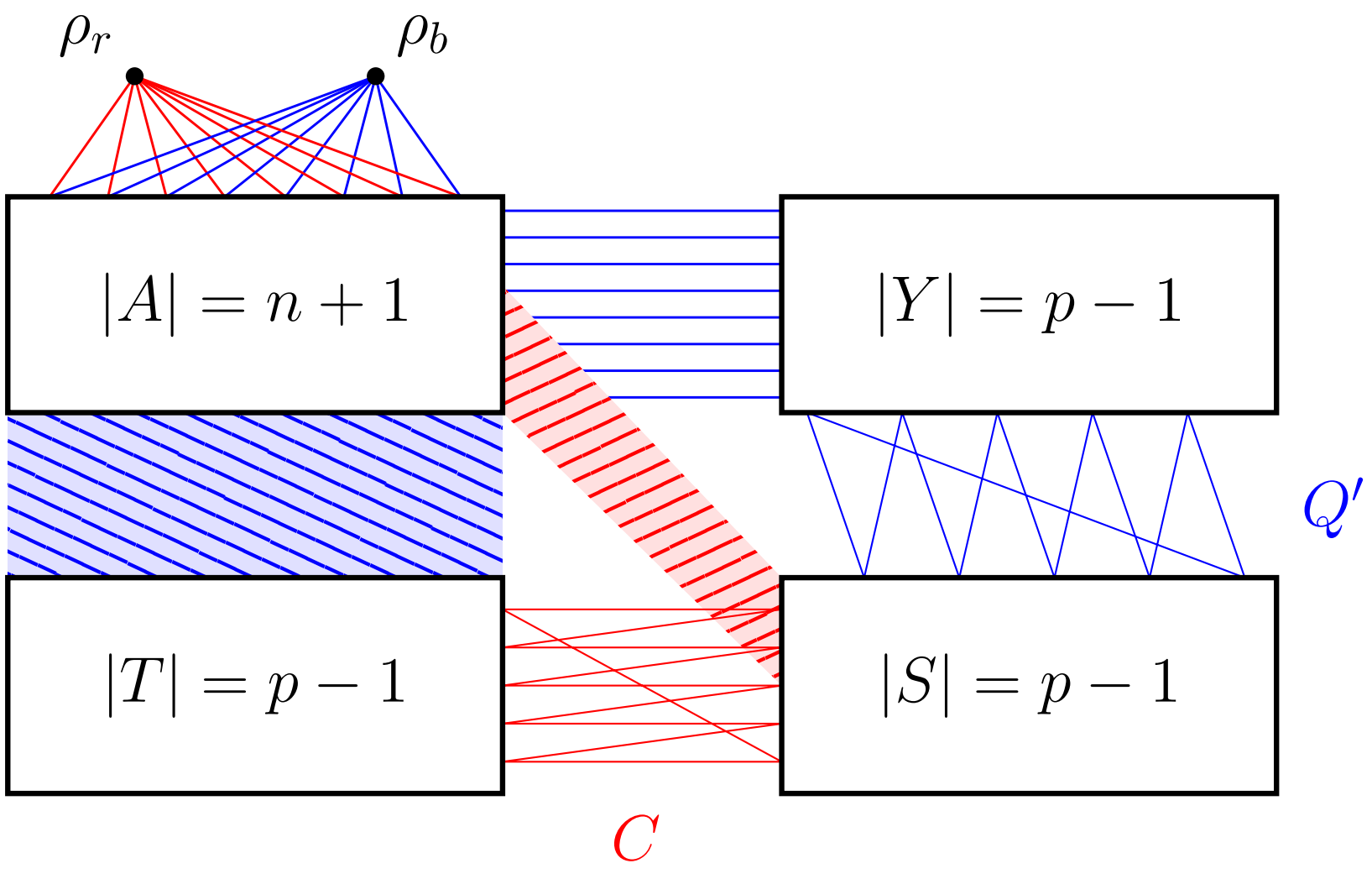
\begin{figure}[h]
\centering
\begin{tikzpicture}[
    box/.style={draw=black, very thick, fill=white},
    redfan/.style={red, line width=0.55pt},
    bluefan/.style={blue, line width=0.55pt},
    redcycle/.style={red},
    bluecycle/.style={blue},
    redhatch/.style={
        pattern={Lines[angle=25,distance=4.5pt,line width=.8pt]},
        pattern color=red
    },
    bluehatch/.style={
        pattern={Lines[angle=-25,distance=4.5pt,line width=.8pt]},
        pattern color=blue
    },
    every node/.style={font=\Large}
]

\def\W{3.9}
\def\H{1.7}
\def\Xgap{2.2}
\def\Ygap{1.3}

\pgfmathsetmacro{\Axc}{\W/2}
\pgfmathsetmacro{\Atopy}{2*\H+\Ygap}

\coordinate (Rsw) at (0,0);
\coordinate (Rne) at ({\W},{\H});

\coordinate (Asw) at (0,{\H+\Ygap});
\coordinate (Ane) at ({\W},{2*\H+\Ygap});

\coordinate (Tsw) at ({\W+\Xgap},0);
\coordinate (Tne) at ({2*\W+\Xgap},{\H});

\coordinate (Ssw) at ({\W+\Xgap},{\H+\Ygap});
\coordinate (Sne) at ({2*\W+\Xgap},{2*\H+\Ygap});


\path[fill=blue!12] (0,{\H}) rectangle ({\W},{\H+\Ygap});
\path[bluehatch]    (0,{\H}) rectangle ({\W},{\H+\Ygap});

\coordinate (ay1) at ({\W},{\H+\Ygap+0.12});
\coordinate (ay2) at ({\W},{\H+\Ygap+0.33});
\coordinate (ay3) at ({\W},{\H+\Ygap+0.54});
\coordinate (ay4) at ({\W},{\H+\Ygap+0.75});
\coordinate (ay5) at ({\W},{\H+\Ygap+0.96});
\coordinate (ay6) at ({\W},{\H+\Ygap+1.17});
\coordinate (ay7) at ({\W},{\H+\Ygap+1.38});
\coordinate (ay8) at ({\W},{\H+\Ygap+1.59});

\coordinate (ya1) at ({\W+\Xgap},{\H+\Ygap+0.12});
\coordinate (ya2) at ({\W+\Xgap},{\H+\Ygap+0.33});
\coordinate (ya3) at ({\W+\Xgap},{\H+\Ygap+0.54});
\coordinate (ya4) at ({\W+\Xgap},{\H+\Ygap+0.75});
\coordinate (ya5) at ({\W+\Xgap},{\H+\Ygap+0.96});
\coordinate (ya6) at ({\W+\Xgap},{\H+\Ygap+1.17});
\coordinate (ya7) at ({\W+\Xgap},{\H+\Ygap+1.38});
\coordinate (ya8) at ({\W+\Xgap},{\H+\Ygap+1.59});

\draw[bluefan] (ay1) -- (ya1);
\draw[bluefan] (ay2) -- (ya2);
\draw[bluefan] (ay3) -- (ya3);
\draw[bluefan] (ay4) -- (ya4);
\draw[bluefan] (ay5) -- (ya5);
\draw[bluefan] (ay6) -- (ya6);
\draw[bluefan] (ay7) -- (ya7);
\draw[bluefan] (ay8) -- (ya8);

\path[fill=red!12]
    ({\W},{\H+\Ygap+1.00}) --
    ({\W+\Xgap},{\H}) --
    ({\W+\Xgap},{\H-0.85}) --
    ({\W},{\H+\Ygap}) --
    cycle;
\path[redhatch]
    ({\W},{\H+\Ygap+1.00}) --
    ({\W+\Xgap},{\H}) --
    ({\W+\Xgap},{\H-0.85}) --
    ({\W},{\H+\Ygap}) --
    cycle;

\def\fansep{0.95}
\def\fanygap{0.95}

\coordinate (u)  at ({\Axc-\fansep},{\Atopy+\fanygap});
\coordinate (up) at ({\Axc+\fansep},{\Atopy+\fanygap});

\foreach \dx in {-1.62,-1.16,-0.70,-0.24,0.24,0.70,1.16,1.62}{
    \pgfmathsetmacro{\rx}{\Axc+\dx}
    \pgfmathsetmacro{\bx}{\Axc-\dx}
    \draw[redfan]  (u)  -- ({\rx},{\Atopy});
    \draw[bluefan] (up) -- ({\bx},{\Atopy});
}

\coordinate (s1) at ({\W+\Xgap+0.20},{\H+\Ygap});
\coordinate (s2) at ({\W+\Xgap+0.95},{\H+\Ygap});
\coordinate (s3) at ({\W+\Xgap+1.70},{\H+\Ygap});
\coordinate (s4) at ({\W+\Xgap+2.45},{\H+\Ygap});
\coordinate (s5) at ({\W+\Xgap+3.20},{\H+\Ygap});

\coordinate (t1) at ({\W+\Xgap+0.65},{\H});
\coordinate (t2) at ({\W+\Xgap+1.40},{\H});
\coordinate (t3) at ({\W+\Xgap+2.15},{\H});
\coordinate (t4) at ({\W+\Xgap+2.90},{\H});
\coordinate (t5) at ({\W+\Xgap+3.65},{\H});

\draw[bluecycle]
    (s1)--(t1)--(s2)--(t2)--(s3)--(t3)--(s4)--(t4)--(s5)--(t5)--cycle;

\coordinate (r1) at ({\W},1.45);
\coordinate (r2) at ({\W},1.15);
\coordinate (r3) at ({\W},0.85);
\coordinate (r4) at ({\W},0.55);
\coordinate (r5) at ({\W},0.25);

\coordinate (q1) at ({\W+\Xgap},1.45);
\coordinate (q2) at ({\W+\Xgap},1.15);
\coordinate (q3) at ({\W+\Xgap},0.85);
\coordinate (q4) at ({\W+\Xgap},0.55);
\coordinate (q5) at ({\W+\Xgap},0.25);

\draw[redcycle]
    (r1)--(q1)--(r2)--(q2)--(r3)--(q3)--(r4)--(q4)--(r5)--(q5)--cycle;

\draw[box] (Asw) rectangle (Ane);
\draw[box] (Rsw) rectangle (Rne);
\draw[box] (Ssw) rectangle (Sne);
\draw[box] (Tsw) rectangle (Tne);

\node at ($(Asw)!0.5!(Ane)$) {$|A| = n+1$};
\node at ($(Rsw)!0.5!(Rne)$) {$|T| = p-1$};
\node at ($(Ssw)!0.5!(Sne)$) {$|Y| = p-1$};
\node at ($(Tsw)!0.5!(Tne)$) {$|S| = p-1$};

\fill (u) circle (2.0pt);
\fill (up) circle (2.0pt);

\node[above left=1pt] at (u) {$\pred$};
\node[above right=1pt] at (up) {$\pblue$};

\node[blue] at ({2*\W+\Xgap+0.45},{\H+0.55}) {$Q^\prime$};
\node[red] at ({\W+1.05},-0.35) {$C$};

\end{tikzpicture}
\caption{The sets $A$, $S$, $T$, and $Y$, and the subgraphs between them.}
\label{fig:fig3}
\end{figure}

    We now disregard the set $Y$. When we repeat the above reasoning once more---using the set $S$ and the cycle $C$---we introduce a final set $R \subset W \setminus{C}$, such that there is a blue cycle $C'$ between $T$ and $R$ and a red complete bipartite graph between $A$ and $R$. The resulting configuration is illustrated in Figure \ref{fig:fig5}.
    
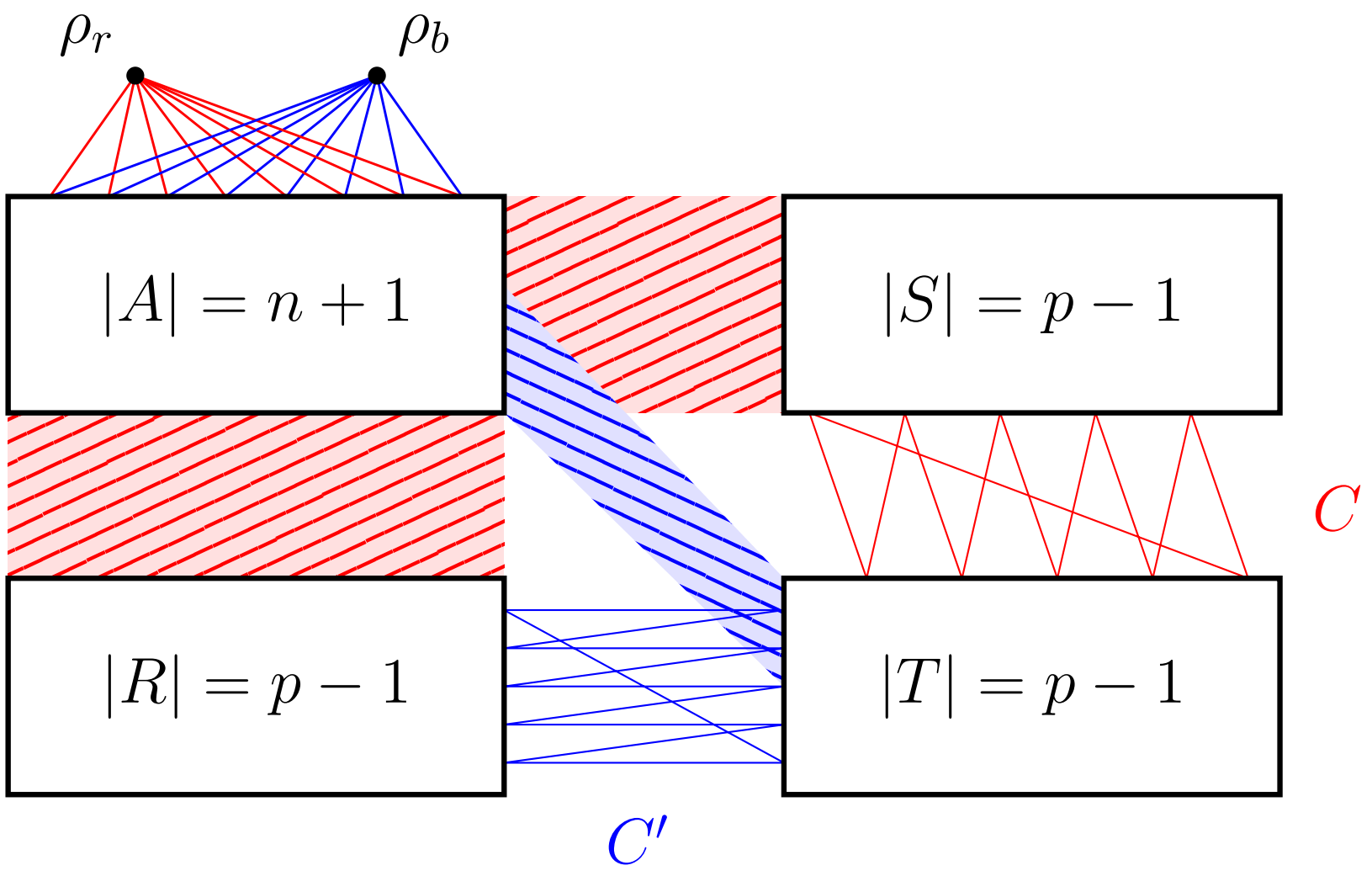
\begin{figure}[h]
    \centering
    \begin{tikzpicture}[
    box/.style={draw=black, very thick, fill=white},
    redfan/.style={red, line width=0.55pt},
    bluefan/.style={blue, line width=0.55pt},
    redcycle/.style={red},
    bluecycle/.style={blue},
    redhatch/.style={
        pattern={Lines[angle=25,distance=4.5pt,line width=.8pt]},
        pattern color=red
    },
    bluehatch/.style={
        pattern={Lines[angle=-25,distance=4.5pt,line width=.8pt]},
        pattern color=blue
    },
    every node/.style={font=\Large}
]

\def\W{3.9}
\def\H{1.7}
\def\Xgap{2.2}
\def\Ygap{1.3}

\pgfmathsetmacro{\Axc}{\W/2}
\pgfmathsetmacro{\Atopy}{2*\H+\Ygap}

\coordinate (Rsw) at (0,0);
\coordinate (Rne) at ({\W},{\H});

\coordinate (Asw) at (0,{\H+\Ygap});
\coordinate (Ane) at ({\W},{2*\H+\Ygap});

\coordinate (Tsw) at ({\W+\Xgap},0);
\coordinate (Tne) at ({2*\W+\Xgap},{\H});

\coordinate (Ssw) at ({\W+\Xgap},{\H+\Ygap});
\coordinate (Sne) at ({2*\W+\Xgap},{2*\H+\Ygap});


\path[fill=red!12] (0,{\H}) rectangle ({\W},{\H+\Ygap});
\path[redhatch]    (0,{\H}) rectangle ({\W},{\H+\Ygap});

\path[fill=red!12]
    ({\W},{\H+\Ygap})
    rectangle
    ({\W+\Xgap},{2*\H+\Ygap});
\path[redhatch]
    ({\W},{\H+\Ygap})
    rectangle
    ({\W+\Xgap},{2*\H+\Ygap});

\path[fill=blue!12]
    ({\W},{\H+\Ygap+1.00}) --
    ({\W+\Xgap},{\H}) --
    ({\W+\Xgap},{\H-0.85}) --
    ({\W},{\H+\Ygap}) --
    cycle;
\path[bluehatch]
    ({\W},{\H+\Ygap+1.00}) --
    ({\W+\Xgap},{\H}) --
    ({\W+\Xgap},{\H-0.85}) --
    ({\W},{\H+\Ygap}) --
    cycle;

\def\fansep{0.95}
\def\fanygap{0.95}

\coordinate (u)  at ({\Axc-\fansep},{\Atopy+\fanygap});
\coordinate (up) at ({\Axc+\fansep},{\Atopy+\fanygap});

\foreach \dx in {-1.62,-1.16,-0.70,-0.24,0.24,0.70,1.16,1.62}{
    \pgfmathsetmacro{\rx}{\Axc+\dx}
    \pgfmathsetmacro{\bx}{\Axc-\dx}
    \draw[redfan]  (u)  -- ({\rx},{\Atopy});
    \draw[bluefan] (up) -- ({\bx},{\Atopy});
}

\coordinate (s1) at ({\W+\Xgap+0.20},{\H+\Ygap});
\coordinate (s2) at ({\W+\Xgap+0.95},{\H+\Ygap});
\coordinate (s3) at ({\W+\Xgap+1.70},{\H+\Ygap});
\coordinate (s4) at ({\W+\Xgap+2.45},{\H+\Ygap});
\coordinate (s5) at ({\W+\Xgap+3.20},{\H+\Ygap});

\coordinate (t1) at ({\W+\Xgap+0.65},{\H});
\coordinate (t2) at ({\W+\Xgap+1.40},{\H});
\coordinate (t3) at ({\W+\Xgap+2.15},{\H});
\coordinate (t4) at ({\W+\Xgap+2.90},{\H});
\coordinate (t5) at ({\W+\Xgap+3.65},{\H});

\draw[redcycle]
    (s1)--(t1)--(s2)--(t2)--(s3)--(t3)--(s4)--(t4)--(s5)--(t5)--cycle;

\coordinate (r1) at ({\W},1.45);
\coordinate (r2) at ({\W},1.15);
\coordinate (r3) at ({\W},0.85);
\coordinate (r4) at ({\W},0.55);
\coordinate (r5) at ({\W},0.25);

\coordinate (q1) at ({\W+\Xgap},1.45);
\coordinate (q2) at ({\W+\Xgap},1.15);
\coordinate (q3) at ({\W+\Xgap},0.85);
\coordinate (q4) at ({\W+\Xgap},0.55);
\coordinate (q5) at ({\W+\Xgap},0.25);

\draw[bluecycle]
    (r1)--(q1)--(r2)--(q2)--(r3)--(q3)--(r4)--(q4)--(r5)--(q5)--cycle;

\draw[box] (Asw) rectangle (Ane);
\draw[box] (Rsw) rectangle (Rne);
\draw[box] (Ssw) rectangle (Sne);
\draw[box] (Tsw) rectangle (Tne);

\node at ($(Asw)!0.5!(Ane)$) {$|A| = n+1$};
\node at ($(Rsw)!0.5!(Rne)$) {$|R| = p-1$};
\node at ($(Ssw)!0.5!(Sne)$) {$|S| = p-1$};
\node at ($(Tsw)!0.5!(Tne)$) {$|T| = p-1$};

\fill (u) circle (2.0pt);
\fill (up) circle (2.0pt);

\node[above left=1pt] at (u) {$\pred$};
\node[above right=1pt] at (up) {$\pblue$};

\node[red] at ({2*\W+\Xgap+0.45},{\H+0.55}) {$C$};
\node[blue] at ({\W+1.05},-0.35) {$C'$};

\end{tikzpicture}
        \caption{The sets $A$, $R$, $S$, and $T$, and the subgraphs between them.}
                \label{fig:fig5}
\end{figure}
    
    The sets $R,S$, and $T$ satisfy the desired four conditions:
    \begin{enumerate}
        \item $R$ and $S$ are both red-complete to $A$,
        \item $T$ is blue-complete to $A$,
        \item There is a red cycle $C$ on $2p-2$ vertices between $S$ and $T$,
        \item There is a blue cycle $C'$ on $2p-2$ vertices between $R$ and $T$.
    \end{enumerate}
By Lemma \ref{lem:helper3}, our $2$-coloring of $K_r$ admits a monochromatic $S_{2p}(n, m)$, a contradiction.
    We conclude that if $\vert W_R\vert, \vert W_B\vert \geq 2p - 2$ and $m \geq \lceil p/2\rceil,$ then our 2-coloring admits a monochromatic $S_{2p}(n, m)$, as desired.    
\end{proof}
\begin{remark}
    Lemma \ref{lem:W_i>=2p-2} shows us that if both $|W_R|$ and $|W_B|$ are large (i.e.\ at least $2p - 2$), then our 2-coloring of $K_r$ admits a monochromatic $S_{2p}(n, m)$. 
Without loss of generality, we may thus assume that $|W_B| \leq 2p - 3$. 
The following lemmas and helpers below prove that our 2-coloring still admits a monochromatic $S_{2p}(n, m)$. 
Lemmas \ref{lem:smallW_B} and \ref{lem:n=m=p-2} below handle the case where $|W_B| \leq p - 3$ and Lemma \ref{lem:medium_W_B} handles $p - 2 \leq |W_B| \leq 2p - 3$.
\end{remark}

\begin{lemma}
    \label{lem:smallW_B}
    Suppose $n \geq p-1$ and $|W_B| \leq p-3$. If our $2$-coloring of $K_r$ has a polarized pair $\{\pred, \pblue\}$, then it admits a monochromatic $S_{2p}(n, m)$.
\end{lemma}
\begin{proof}
    We write that $|W_B| = p-3-k$, for some $k \geq 0$. Because there are $|W_B|$ vertices 
    that are connected to some vertex in $A$ by a blue edge, the remaining $|W^\prime| = |W| - |W_B| = \max(n + p - 1 + k, 2m + 2p -1 + k)$ vertices in $W$ are red-complete to $A$.
    
    Since $A \subseteq R(\pred)$, there exists a red complete bipartite graph between $A$ and $W^\prime \cup \{\pred\}$, with sizes $n+1$ and $\max(n + p + k, 2m + 2p + k)$, respectively. If $n \geq m + p - 1$ and thus $|A| \geq m + p$, then we can form a red $S_{2p}(n, m)$ as follows. We first use the complete bipartite graph to create a path on $2p$ vertices starting at $\pred$ and ending at some vertex $a \in A$. This path exhausts $p$ vertices from each partition; therefore, $\pred$ has at least $|A| - p \geq m$ neighbors in $A$ not used by the path, so we can use $\pred$ as the center of the $m$-star. Similarly, $a$ has at least $|W^\prime \cup \{\pred\}| - p \geq \max(n + k, 2m + p + k) \geq n$ neighbors in $W^\prime \cup \{\pred\}$ not used by the path, so we can use $a$ as the center of the $n$-star.
    
    Hence, we may assume that $n \leq m + p - 2$. It follows that \[|W' \cup \{ \pred \}| = 2m + 2p + k \geq n + m + p + k + 2.\]
    Suppose there is a red edge $\{w_1, w_2\}$ inside $W^\prime \cup \{\pred\}$. Then, choose $p$ distinct vertices $a_1, a_2, \ldots, a_p \in A$ and $p$ distinct vertices $w_1, w_2, \ldots, w_p \in W^\prime \cup \{\pred\}$, and form the red path
    \[a_1 - w_1 - w_2 - a_2 - w_3 - \cdots - w_p - a_p\]
    on $2p$ vertices (using the red complete bipartite graph between $A$ and $W^\prime \cup \{\pred\}$). This path leaves behind at least $|W^\prime \cup \{\pred\}| - p \geq 2m + p + k \geq m + n$ unused vertices in $W^\prime \cup \{\pred\}$. Hence, we can complete a red $S_{2p}(n, m)$, with $a_1$ as the center of our $n$-star and $a_p$ as the center of our $m$-star (with spokes in $W^\prime \cup \{\pred\}$), or vice versa.
    Therefore, we may assume that $W^\prime \cup \{ \pred \}$ induces a blue complete graph. Suppose now that $\pblue$ is blue-adjacent to some vertex $w \in W^\prime \cup \{\pred\}$. By our polarized pair construction, we know that $\pblue$ is blue-adjacent to every vertex in $A$. As a result, we can form a blue $S_{2p}(n, m)$ as follows: use $\pblue$ as the center of the $n$-star (with spokes in $A$); connect $\pblue$ to $w$, and then use $w$ as the tail of a blue $S_{2p - 1}(m, 0)$ in $W^\prime \cup \{\pred\}$. Hence, we may assume that $\pblue$ is red-adjacent to every vertex in $W^\prime \cup \{\pred\}$.
    
    Finally, consider $W_B$. By definition, $W_B$ is the set of vertices that are connected to $A$ by a blue edge. Thus, if there is a blue edge between some vertex $w_1 \in W_B$ and $w_2 \in W' \cup \{ \pred \}$, then we can complete a blue $S_{2p}(n, m)$ as follows: form a blue broom $S_{2p - 3}(m, 0)$ in $W^\prime \cup \{\pred\}$ with tail vertex $w_2$, and connect $w_2 \rightarrow w_1 \rightarrow a$, for some blue neighbor $a \in B_A(w_1)$; connect $a$ to $\pblue$, and use $\pblue$ as the center of the $n$-star (with spokes in $A \setminus \{a\}$). As a result, we may assume that there is a red complete bipartite graph between $W_B$ and $W' \cup \{ \pred \}$. 
    
    If $n \leq m + k$, then $W' \cup \{ \pred \}$ admits a blue $S_{2p}(n,m)$. On the other hand, if $n \geq m +k +1$, then $|A \cup W_B \cup \{ \pblue \}| \geq m + p$ and $|W' \cup \{ \pred \}| \geq n + m + p + k + 2 > n + p$, so the complete bipartite graph admits a red $S_{2p}(n,m)$. In either case, our $2$-coloring admits a monochromatic $S_{2p}(n,m)$, as desired.
\end{proof}

\begin{lemma}
    \label{lem:n=m=p-2}
    Suppose $n = m = p - 2$ and $|W_B| \leq p - 3$. If our $2$-coloring of $K_r$ has a polarized pair $\{\pred, \pblue\}$, then it admits a monochromatic $S_{2p}(n, m)$.
\end{lemma}
\begin{proof}
    We write that $|W_B| = p-3-k$, for some $k \geq 0$. 

    Suppose that $W^\prime \cup \{\pred\}$ either admits a red $P_4$ or a red $P_2 \sqcup P_3$; call this red subgraph $U$. Now, choose $(p + 1) - |U|$ singleton vertices in $W^\prime \cup \{\pred\}$ that do not belong to $U$ (note that this is possible, since $(p + 1) - |U| \geq (p + 1) - 5 \geq p - 4 \geq 0$). These singleton vertices along with $U$ create a red forest on $p - 2$ connected components, say, $F_1, F_2, \ldots, F_{p - 2}$; the subgraph $U$ either occupies the component $F_1$ or both $F_1$ and $F_2$, depending on whether $U = P_4$ or $U = P_2 \sqcup P_3$, respectively. Since $A$ is red-complete to $W^\prime \cup \{\pred\}$, we can form the path
    \[P = a_1 - F_1 - a_2 - F_2 - \cdots - a_{p - 2} - F_{p - 2} - a_{p - 1}\]
    on $2p$ vertices, where $a_1,a_2, \ldots, a_{p - 1} \in A$. This path occupies $p + 1$ vertices in $W^\prime \cup \{\pred\}$, leaving behind $|W^\prime \cup \{\pred\}| - (p + 1) = (4p - 4 + k) - (p + 1) \geq 3p -5 \geq 2p - 4 = n + m$ unused vertices in $W^\prime \cup \{\pred\}$. Hence, we can use $a_1$ as the center of our $n$-star and $a_{p - 1}$ as the center of our $m$-star (or vice versa, with spokes in $(W^\prime \cup \{\pred\}) \setminus V(P)$), completing a red $S_{2p}(n, m)$.

    Therefore, we may assume that $W^\prime \cup \{\pred\}$ neither contains a red $P_4$ nor a red $P_2 \sqcup P_3$. Since every red connected component of $W^\prime \cup \{\pred\}$ avoids $P_4$ as a subgraph,\footnote{Not to be confused with a \emph{cograph} or ``$P_4$-free graph,\rq\rq\ which only avoids $P_4$ as an \emph{induced} subgraph.} each such component must either be a singleton, a triangle, or a star (possibly on two vertices, i.e.\ $P_2$). However, if any one component contains a $P_3$, then no other component can contain a $P_2$ and they are thus  singletons. As a result, the red subgraph induced by $W^\prime \cup \{\pred\}$ is either a matching (possibly empty) with added singletons, a triangle with singletons, or a star with singletons. We will consider each case below.

    Let us begin with the matching case. First, choose some subset $S \subseteq W^\prime \cup \{\pred\}$ on $4p - 4$ vertices. We will partition $S$ into two disjoint sets $X$ and $Y$. If $S$ contains vertices that are already matched in $W^\prime \cup \{\pred\}$, then leave them paired; pair the remaining unmatched vertices (i.e.\ without a partner in $S$) arbitrarily. This creates a total of $(4p - 4)/2 = 2p - 2$ pairs. Add half of these pairs to $X$ and half to $Y$, so that $|X| = |Y| = 2p - 2$. Notice that $X$ is blue-complete to $Y$, since any red-adjacent pair in $S$ was added to the same partition. We can use this blue complete graph to form a blue path $P$ on $2p$ vertices with endpoints $x \in X$ and $y \in Y$, leaving behind $(2p - 2) - p = p - 2$ unused vertices on either side. But $n = m = p - 2$, so we can use $x$ as the center of our $n$-star (with spokes in $Y \setminus P$) and $y$ as the center of our $m$-star (with spokes in $X \setminus P$), or vice versa, completing a blue $S_{2p}(n, m)$.

    Now, consider the triangle case. Every edge in our $2$-coloring of $W^\prime \cup \{\pred\}$ is blue except for the red triangle; in other words, $W^\prime \cup \{\pred\}$ consists of a red $K_3$ that is blue-complete to a blue complete graph on $|W^\prime \cup \{\pred\} | - 3 = (4p - 4 + k) -3 = 4p - 7 + k$ vertices. Since $4p - 7 + k \geq (n - 1) + (m - 2) + 2p$, we can use this blue $K_{4p - 7 + k}$ to form a blue $S_{2p}(n - 1, m - 2)$. We can supply the remaining spokes of the $n$- and $m$-stars here using the vertices of the red $K_3$, finishing a blue $S_{2p}(n, m)$.

    Finally, suppose that $W^\prime \cup \{\pred\}$ induces a red star. Every edge in our $2$-coloring of $W^\prime \cup \{\pred\}$ is blue except for this red star, so deleting its center induces a blue complete graph on $|W^\prime \cup \{\pred\}| - 1 = 4p - 5 + k$ vertices. Choose some subset $S^\prime$ of exactly $4p - 5 = n + m + 2p - 1$ vertices from this blue complete graph. If there is a blue edge from any vertex $s \in S^\prime$ to $x \in V(K_r) \setminus S^\prime$, then we can complete a blue $S_{2p}(n, m)$: we use 
    $S^\prime$ to form a blue $S_{2p}(n - 1, m)$ with $(n - 1)$-star centered at $s$, and then use $x$ to finish the $n$-star. Hence, we may assume that $S^\prime$ is red-complete to $V(K_r) \setminus S^\prime$. But $|S^\prime| \geq n + p$ and $|V(K_r) \setminus S^\prime| = r - (4p - 5) = r - (n + m + 2p - 1) = m + p$, so we have a red $K_{n + p, m + p}$, which yields a red $S_{2p}(n, m)$.

    In all three cases, our $2$-coloring admits a monochromatic $S_{2p}(n, m)$, as desired.
\end{proof}

\begin{lemma}\label{lem:helper2}
    Suppose $|W_B| \leq 2p-3$. If there is a red path on $2p-1$ vertices in $W$ with both endpoints in $W'$ and also at most $m + p - 1$ vertices in $W'$, then our 2-coloring of $K_r$ admits a red $S_{2p}(n,m)$.
\end{lemma}

\begin{proof}
    The proof here is straightforward. By the definition of $W^\prime = W \setminus W_B$, there exists a red complete bipartite graph between $W^\prime$ and $A$.
    
    Now, assume that our $2$-coloring of $K_r$ admits such a red path $P$ with endpoints $w_1,  w_2 \in W^\prime$. Since the path from $w_1$ to $w_2$ only uses at most $m + p - 1$ vertices in $W^\prime$, there are at least $m$ vertices left unused in $W^\prime$. We can thus complete a red $S_{2p}(n,m)$ as follows: choose any vertex $a \in A$, and use $a$ as the center of our $m$-star (with spokes in $W^\prime \setminus V(P)$); connect $a$ to $w_1$, and follow the path on $2p - 1$ vertices from $w_1$ to $w_2$; finally, use $w_2$ as the center of our $n$-star, with spokes in $A \setminus \{a\}$ (note that $|A| = n + 1$). Hence, our $2$-coloring of $K_r$ admits a red $S_{2p}(n, m)$, as desired.
\end{proof}

\begin{lemma}\label{lem:helper_red_cycle_case}
     Suppose $p - 2 \leq |W_B| \leq 2p-3$. 
     If $W^\prime$ contains a red cycle on $2p - 2$ vertices and our $2$-coloring of $K_r$ has a polarized pair $\{\pred, \pblue\}$, then it admits a monochromatic $S_{2p}(n, m)$.
\end{lemma}

\begin{proof}
    Let $C$ be our red cycle on $2p - 2$ vertices in $W'$.
    If any vertex $c\in C$ satisfies $\text{red-}d_{W \setminus C}(c) \geq m$, then we can complete a red $S_{2p}(n, m)$ as follows: use $c$ as the center of the $m$-star (with spokes in $W \setminus C$), and follow $C$ back for $2p-2$ vertices from $c$ (including $c$), ending at some other vertex $c' \in C$; connect $c' \rightarrow a \rightarrow \pred$, using $\pred$ as the center of the $n$-star (with spokes in $A \setminus \{a\}$). 
    Thus, we may assume that $\text{red-}d_{W\setminus C} (c) \leq m -1$ for all $c \in C$, so
    \[\text{blue-}d_{W \setminus C} (c) \geq |W| - |C| - (m-1) \geq 2m + 3p -4 - (2p-2) - (m-1) = m + p - 1.\]
    Because $C$ is an even cycle, it has two classes of $p - 1$ alternating vertices on the cycle; let $S$ be either of these classes.  Since $p \geq 4$, we have that $|W \setminus C| \geq 2m + 2$ and thus $|W \setminus C| \leq 2(|W \setminus C| - m + 1) - 2$. As a result, Lemma \ref{Jackson_cycles} with $a = p -1$, $b = |W \setminus C|$, and $k = |W \setminus C| - m + 1$ gives us a blue cycle on $2p-2$ vertices between $S$ and a subset $T \subset W \setminus C$, where $|T| = p - 1$. Call this blue cycle $C^\prime$, so that $S = C \cap C^\prime$.

    Since $S \subset C$, we have that $\text{blue-}d_{W \setminus C}(s) \geq m + p - 1 \geq m + 3$ for all $s \in S$. Thus, every vertex $s \in S$ forms an $m$-star in $W \setminus (C \cup C^\prime)$ (setting aside all $p-1$ vertices in $T$). If there exists a blue edge between some vertex $t \in T$ and $a \in A$, then we can complete a blue $S_{2p}(n, m)$ as follows: use $\pblue$ as the center of our $n$-star (with spokes in $A \setminus \{a\}$), and connect $\pblue \rightarrow a \rightarrow t$; follow $C'$ for $2p -2$ vertices from $t$ (including $t$) to some vertex $s \in S$, and use $s$ as the center of our $m$-star (with spokes in $W \setminus (C \cup C^\prime)$). Hence, we may assume that $T \cap W_B = \emptyset$, meaning that $T$ is red-complete to $A$, and thus $T \subset W^\prime$.
    
    Similarly, suppose there is a blue edge between some vertex $s \in S$ and $w \in W_B$. We can then complete a blue $S_{2p}(n,m)$ in the same manner: use $\pblue$ as the center of our $n$-star (with spokes in $A \setminus \{a\}$), and connect $\pblue \rightarrow a \rightarrow w \rightarrow s$; follow $C'$ for $2p-3$ vertices (including $s$), ending at some other vertex $s' \in S$, which we use as the $m$-star. Note that $s'$ has at least $m+p-1$ blue neighbors in $W\setminus C$. Our path occupies $p-1$ vertices in $W \setminus C$, so excluding these vertices leaves at least $(m+p-1) - (p- 1) = m$ unused blue neighbors of $s^\prime$. Thus, we can use $s^\prime$ as the center of our $m$-star (with spokes in $W \setminus C$), completing a blue $S_{2p}(n, m)$. As a result, we may assume that $S$ is red-complete to $W_B$.
    
    Using this complete bipartite graph, we can form a red path on $2p-3$ vertices, with both endpoints in $S$. 
    But recall that $S$ is only one of $C$'s two alternating classes of vertices. Both of the endpoints in $S$ have two red neighbors in the other alternating class, so we can extend our path by two vertices, appending one to each side. This creates a red path on $2p - 1$ vertices with both endpoints in $W^\prime$. Because this path also has $p + 1 \leq m + p - 1$ vertices in $W^\prime$, Lemma \ref{lem:helper2} gives us a red $S_{2p}(n,m)$, as desired.
\end{proof}
\begin{lemma}[Due to Rosta, \cite{Rosta1}, and Faudree and Schelp, \cite{FAUDREE1974313}]\label{lem:ramsey_numbers_of_cycles}
    Let $C_k$ denote a cycle graph on $k$ vertices. The Ramsey numbers of cycles are given by
    
    $$r(C_k, C_{\ell}) = \begin{cases}
        2\ell -1 , & \text{for } 3\leq k \leq \ell \text{, }k \text{ odd, } (k, \ell) \neq (3,3), \\
        \ell + k/2 - 1, & \text{for } 4\leq k \leq \ell \text{, }k \text{ and } \ell \text{ even, } (k, \ell) \neq (4,4),\\
        \max(\ell + k/2 - 1, 2k -1), & \text{for } 4\leq k < \ell \text{, } k \text{ even and }\ell \text{ odd}. 
        \end{cases}$$

\end{lemma}
\newpage
\begin{lemma}
    \label{lem:medium_W_B}
    Suppose $m \geq p - 2$ and $p-2 \leq |W_B| \leq 2p-3$. If our $2$-coloring of $K_r$ has a polarized pair $\{\pred, \pblue\}$, then it admits a monochromatic $S_{2p}(n,m)$.
\end{lemma}

\begin{proof}
    If $m \geq p-2$, then $|W' \cup \{\pred\}| \geq 2m + p \geq 3p -4$. By Lemma \ref{lem:ramsey_numbers_of_cycles}, $r(C_{2p - 2}) = 3p - 4 \leq |W^\prime \cup \{\pred\}|$, so $W^\prime \cup \{\pred\}$ admits a monochromatic cycle on $2p - 2$ vertices; call this cycle $C$. If $C$ uses the vertex $\pred$, then for the rest of the lemma take $\pred$ to be some other vertex $\pred^\text{new} \in W^\prime \setminus C$, which is necessarily also red-complete to $A$; in this case, exclude $\pred = \pred^\text{new}$ but include the original $\pred^\text{old}$ in $W$ and $W^\prime$. If $C$ is red, then Lemma \ref{lem:helper_red_cycle_case} gives us a monochromatic $S_{2p}(n, m)$. Thus, we may assume that $C$ is blue.

    Now, suppose that there exist $p-1$ vertices on $C$ adjacent to $W_B$ by at least one blue edge; call this subset of vertices $T$. Consider, across all choices of the vertex $t \in T$, the set $S$ of vertices reachable by following $C$ for $2p - 3$ vertices from $t$ (including $t$) along some fixed orientation of $C$. Note that $|S| = p - 1$, since our construction of $S$ simply performs a cyclic shift of $T$.  If any vertex $s \in S$ has $\text{blue-}d_{W \setminus C}(s) \geq m$, then we can form a blue $S_{2p}(n, m)$ as follows: follow the cycle $C$ for $2p - 3$ vertices from $s$ (including $s$), ending at some vertex $t \in T$ and only omitting the vertex $z \in C$ (a blue neighbor of $s$); connect $t \rightarrow w \rightarrow a \rightarrow \pblue$, where $w \in B_{W_B}(t)$ and $a \in A$; use $\pblue$ as the center of our $n$-star (with spokes in $A \setminus \{a\}$) and $s$ as the center of our $m$-star (with $m - 1$ spokes in $W \setminus (C \cup \{w\})$ and $z$ as the last spoke).
 
    By similar reasoning to Lemma \ref{lem:helper_red_cycle_case} (which invokes Lemma \ref{Jackson_cycles}), there exists a red cycle $C'$ on $2p-2$ vertices between $S$ and $R \subset W \setminus C$, such that $|R| = p - 1$ and $\text{red-}d_{W \setminus (C \cup C^\prime)}(s) \geq m$ for all $s \in S$ (i.e.\ $s$ forms a red $m$-star that does not intersect $C^\prime$). If there is a red edge between any vertex $r \in R$ and $a \in A$, then we can complete a red $S_{2p}(n, m)$ as follows: use $\pred$ as the center of our $n$-star (with spokes in $A \setminus \{a\}$), and connect $\pred \rightarrow a \rightarrow r$; follow the cycle $C^\prime$ from $r$ for $2p - 2$ vertices from $r$ (including $r$), ending at some vertex $s \in S$; use $s$ as the center of our $m$-star (with spokes in $W \setminus (C \cup C^\prime)$). Hence, we may assume that $R$ is blue-complete to $A$, and thus $R \subset W_B$.

    If any vertex $r \in R$ has $\text{red-}d_{W \setminus C^\prime}(r) \geq m$, then we can complete a red $S_{2p}(n, m)$ as follows: use $r$ as the center of our $m$-star (with spokes in $W \setminus C^\prime$), and follow the cycle $C^\prime$ for $2p - 2$ vertices from $r$ (including $r$), ending at some vertex $s \in S$; connect $s \rightarrow a \rightarrow \pred$ for some  $a \in A$, and use $\pred$ as the center of our $n$-star (with spokes in $A \setminus \{a\}$). Consequently, we may assume that $\text{red-}d_{W \setminus C^\prime}(r) \leq m - 1$ and thus $\text{blue-}d_{W \setminus C^\prime}(r) \geq |W \setminus C^\prime| - (m - 1)$ for all $r \in R$. By similar reasoning to before, there exists a blue cycle $C^{\prime\prime}$ on $2p - 2$ vertices between $R$ and $D \subset W \setminus C^\prime$, where $|D| = p - 1$; we also have that $\text{blue-}d_{W \setminus (C^\prime \cup C^{\prime\prime})}(r) \geq m$ for all $r \in R$. By applying an analogous $S_{2p}(n, m)$ construction---with $D$ replacing $R$, and blue instead of red---we may assume that $D$ is red-complete to $A$. But we now have three pairwise disjoint sets $D, R, S \subset W$ with $|D| = |R| = |S| = p - 1$, such that 
    \begin{enumerate}
        \item $D$ and $S$ are both red-complete to $A$,
        \item $R$ is blue-complete to $A$,
        \item There is a red cycle $C^\prime$ on $2p - 2$ vertices between $R$ and $S$,
        \item There is a blue cycle $C^{\prime\prime}$ on $2p - 2$ vertices between $R$ and $D$.
    \end{enumerate}
    By Lemma 
    \ref{lem:helper3} (with $D$, $S$, and $R$ in place of the lemma's $R$, $S$, and $T$), our $2$-coloring of $K_r$ admits a monochromatic $S_{2p}(n, m)$.

    Therefore, our original assumption---that there exist $p - 1$ vertices on $C$ blue-adjacent to some vertex in $W_B$---is false. In other words, there exists a red complete bipartite graph between $W_B$ and some subset $S^\prime$ of at least $p$ vertices on $C$. If $|W_B| \geq p - 1$, then this red complete bipartite graph produces a red path on $2p-1$ vertices, with both endpoints in $S^\prime \subset W^\prime$ as well as $p \leq m + p - 1$ vertices in $W^\prime$; thus, Lemma \ref{lem:helper2} gives us a red $S_{2p}(n, m)$. As a result, we may assume that $|W_B| = p - 2$.

     Finally, suppose there exists a red edge $\{x, y\}$ inside $W^\prime$, and choose any two vertices $a_1, a_2 \in A$. If $|S^\prime \setminus \{x, y\}| \geq p - 1$, then let $S^{\prime\prime}$ be any subset of $p - 1$ vertices from $S' \setminus \{x, y\}$. Because $|W_B| = p - 2$, there exists a red path on $2p - 3$ vertices between $W_B$ and $S^{\prime\prime}$ with both endpoints $w_1^{\prime\prime}, w_2^{\prime\prime} \in S^{\prime\prime} \subset W^\prime$. Since there is a complete red bipartite graph between $W^\prime$ and $A$, we also have that $w_1^{\prime\prime}$ and $w_2^{\prime\prime}$ are red-adjacent to $a_1$ and $a_2$. Thus, we can complete a red $S_{2p}(n, m)$ as follows: use $x$ as the center of our $n$-star (with one spoke at $y$, and $n-1$ spokes in $A \setminus \{a_1, a_2\}$); connect $x \rightarrow a_1 \rightarrow w_1^{\prime\prime}$, and follow the path on $2p - 3$ vertices from $w_1^{\prime\prime}$ to $w_2^{\prime\prime}$; connect $w_2^{\prime\prime}$ to $a_2$, and use $a_2$ as the center of our $m$-star (with spokes in $W^\prime \setminus (S^{\prime\prime} \cup \{x, y\})$). Note that $|W^\prime \setminus (S^{\prime\prime} \cup \{x, y\})| \geq |W^\prime| - (p - 1) - 2 \geq 2m - 2 \geq m$.
     
     On the other hand, suppose that $|S^\prime \setminus \{x, y\}| \leq p -2$; since $|S^\prime| \geq p$, this means that $|S^\prime| = p$ and $x, y \in S^\prime$. As such, we write that $S^\prime \setminus \{x, y\} = \{s_1, s_2, \ldots, s_{p - 2}\}$ and $W_B = \{b_1, b_2, \ldots, b_{p - 2}\}$. We can now form the red path
      \[a_1 - z - a_2 - s_1 - b_1 - s_2 - b_2 - \cdots - s_{p - 2} - b_{p - 2} - x\]
      on $2p$ vertices, where $z \in W^\prime \setminus S^\prime$. This path uses $p$ vertices in $S^\prime \cup \{z\}$, leaving behind $|W^\prime| - p - 1 \geq  2m - 2 \geq m$ unused vertices in $W^\prime$. Therefore, we can complete a red $S_{2p}(n, m)$, with $a_1$ as the center of our $m$-star (with spokes in $W^\prime \setminus (S^\prime \cup \{z\})$) and $x$ as the center of our $n$-star (with $n - 1$ spokes in $A \setminus \{a_1, a_2\}$ and $y$ as the last spoke).
      
      Hence, we may assume that there are no red edges inside $W^\prime$, and so $W^\prime$ induces a blue complete graph. By similar reasoning, there cannot exist a red edge between $\pblue$ and some vertex $x \in W^\prime$, since we can use $\pblue$ to supplement the red $n$-star at $x$ (in the same fashion as $y$ above). As a result, $\pblue$ is blue-complete to $W^\prime$. But we can now complete a blue $S_{2p}(n, m)$ by choosing any blue broom $S_{2p - 1}(m, 0)$ in $W^\prime$, and then connecting its tail vertex to $\pblue$, which we use as the center of our $n$-star (with spokes in $A$). We conclude that our $2$-coloring of $K_r$ admits a monochromatic $S_{2p}(n, m)$, as desired.
\end{proof}

\begin{lemma}
    \label{lem:polarized_pair}
    Suppose $m \geq p - 2$. If our 2-coloring of $K_r$ has a polarized pair, then it admits a monochromatic $S_{2p}(n, m)$.
\end{lemma}
\begin{proof}
    If $|W_R| \geq 2p - 2$ and $|W_B| \geq 2p - 2$, then Lemma \ref{lem:W_i>=2p-2} gives us a monochromatic $S_{2p}(n, m)$. Without loss of generality, then, suppose that $|W_B| \leq 2p - 3$. If $|W_B| \leq p - 3$ and $n \geq p - 1$, then Lemma \ref{lem:smallW_B} also yields a monochromatic $S_{2p}(n, m)$; if $|W_B| \leq p - 3$ but $n= m = p -2$, then Lemma 
    \ref{lem:n=m=p-2} does the same.
    Finally, if $p - 2 \leq |W_B| \leq 2p - 3$, then Lemma \ref{lem:medium_W_B} provides a monochromatic $S_{2p}(n, m)$. In all cases, we find that our 2-coloring of $K_r$ admits a monochromatic $S_{2p}(n, m)$, as desired.
\end{proof}

\subsection{No Polarized Pair Exists}
\label{sec:no_polarized_pair_exists}
In this section, we will show that if our $2$-coloring of $K_r$ does not admit a polarized pair, then it must admit a monochromatic $S_{2p}(n, m)$. Equivalently, we prove that any $2$-coloring of $K_r$ admits a polarized pair or a monochromatic $S_{2p}(n, m)$.
\subsubsection{Notation}
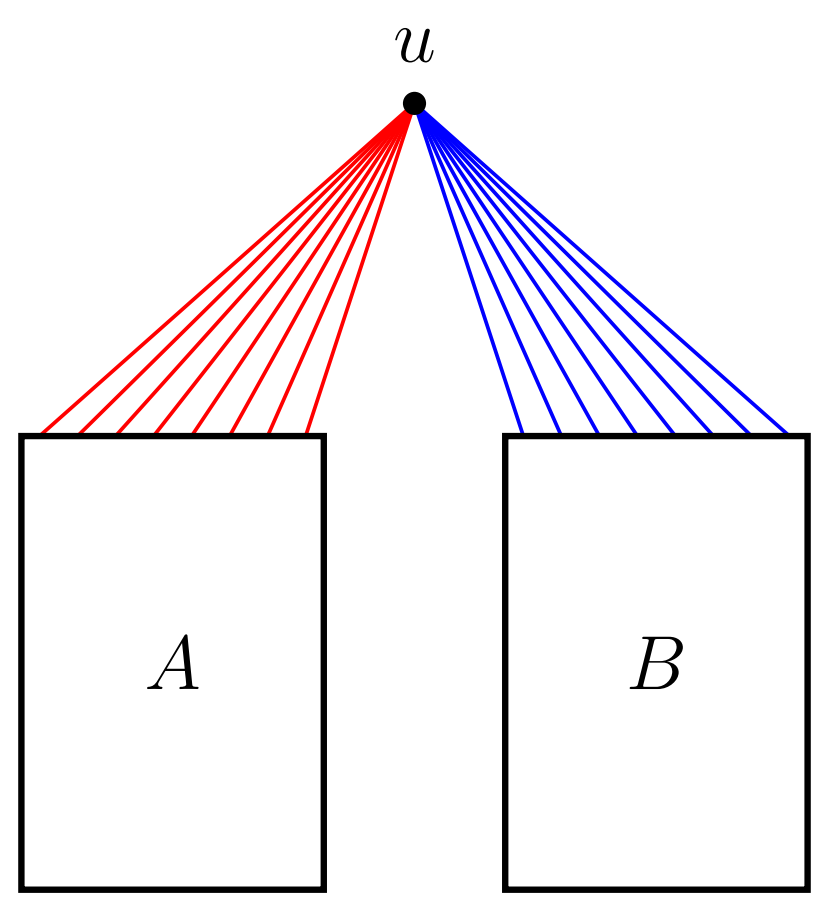
\begin{figure}[h]
\centering
\begin{tikzpicture}[
    box/.style={draw=black, very thick, fill=white},
    vertex/.style={circle, draw=black, fill=white, inner sep=1.4pt},
    rededge/.style={red, line width=0.7pt},
    blueedge/.style={blue, line width=0.7pt},
    every node/.style={font=\Large}
]

\def\W{2}
\def\H{3}
\def\Gap{1.2}
\def\FanHeight{1.7}


\coordinate (Asw) at (0,0);
\coordinate (Ane) at ({\W},{\H});

\coordinate (Bsw) at ({\W+\Gap},0);
\coordinate (Bne) at ({2*\W+\Gap},{\H});

\coordinate (u) at ({\W+\Gap/2},{\H+\FanHeight+0.5});

\foreach \i in {0,...,7} {
    \draw[rededge]
        (u) -- ({0.12+(\W-0.24)*\i/7},{\H});
}

\foreach \i in {0,...,7} {
    \draw[blueedge]
        (u) -- ({\W+\Gap+0.12+(\W-0.24)*\i/7},{\H});
}

\draw[box] (Asw) rectangle (Ane);
\draw[box] (Bsw) rectangle (Bne);

\node at ({\W/2},{\H/2}) {$A$};
\node at ({1.5*\W+\Gap},{\H/2}) {$B$};

\node[vertex, fill=black] at (u) {};
\node[above=4pt] at (u) {$u$};

\end{tikzpicture}

\caption{The sets $A$ and $B$, and their connections to the vertex $u$.}
\label{fig:AB-u}
\end{figure}

Once again, we define the following sets and variables with respect to our chosen $2$-coloring; see Figure \ref{fig:AB-u}.
\begin{itemize}
    \item Let $u$ be any vertex in our $2$-coloring with maximum monochromatic degree; without loss of generality, assume that $u$ is maximal in the color red.
    \item Let $A = R(u)$, the red neighborhood of $u$. Since $r \geq 2n + 2p - 1$ and $|A| \geq (r - 1)/2$, we write that $|A| = n + p - 1 + k$, where $k \geq 0$.
    \item Let $B = V(K_r) \setminus (A \cup \{u\})$, the remaining vertices outside of $u$ and $A$.
\end{itemize}

We proceed with Lemma \ref{lem:min_red_blue_degree} below.
\begin{lemma}\label{lem:min_red_blue_degree}
    For all vertices $v \in V$,  $\text{red-}d(v) \geq |B|$ and $\text{blue-}d(v) \geq |B|$.
\end{lemma}

\begin{proof}
    Suppose, without loss of generality, that there exists a vertex $v \in V$ with $\text{red-}d(v) \leq |B|-1$. Then, $\text{blue-}d(v) \geq |A|+1$, which contradicts our choice of $u$ as a vertex with maximum monochromatic degree. Therefore, every vertex $v \in V$ has $\text{red-}d(v) \geq |B|$ and similarly $\text{blue-}d(v) \geq |B|$, as desired.
\end{proof}
\newpage
\begin{lemma}\label{lem:min_path_degree}
    Let $G$ be a graph with minimum degree at least $s$, and let $z$ be any vertex in $G$. Then, $G$ contains a path on $s+1$ vertices with $z$ as an endpoint.
\end{lemma}

\begin{proof}
    Let $P$ be a path of maximal length starting at $z$, and call the other endpoint $z^\prime$. Since $P$ is maximal, each of the neighbors of $z^\prime$ must already be on the path; otherwise, we could append one of these neighbors to $P$ at $z^\prime$, making a longer path. But $z^\prime$ has at least $\delta(G) \geq s$ neighbors, so $P$ includes at least $s + 1$ vertices (including $z^\prime$), as desired.
\end{proof}

\begin{lemma}\label{lem:max_blue_degree}
    If our 2-coloring of $K_r$ does not admit a polarized pair, then $\text{blue-}d_A(v) \leq n$ and thus $\text{red-}d_A(v) \geq p + k - 2$ for all $v \in V$.
\end{lemma}

\begin{proof}
    If there exists a vertex $v \in V$ with $\text{blue-}d_A(v) \geq n + 1$, then $|R(u) \cap B(v)| \geq n + 1$, so $\{u, v\}$ is a polarized pair. The lemma statement follows immediately.
\end{proof}

\begin{corollary}\label{cor:max_blue_degree}
    If our $2$-coloring of $K_r$ does not admit a polarized pair, then every vertex $a \in A$ is the endpoint of a red path on $p+k-1$ vertices in $A$.
\end{corollary}
\begin{proof}
    By Lemma \ref{lem:max_blue_degree}, $\text{blue-}d_A(a) \leq n$ for all $a \in A$, and so \[\text{red-}d_A(a) \geq |A \setminus \{ a \}| - n  = p +k - 2.\] By Lemma \ref{lem:min_path_degree}, every vertex $a \in A$ is therefore the endpoint of a red path on $p+k-1$  vertices contained in $A$.
\end{proof}

\begin{lemma}
    \label{lem:p_leq_k_leq_m_plus_p_minus_1}
    If $p \leq k \leq m + p -1$, then our 2-coloring of $K_r$ either admits a polarized pair or a monochromatic $S_{2p}(n,m)$.
\end{lemma}

\begin{proof}
    Suppose that our $2$-coloring of $K_r$ does not admit a polarized pair. By Corollary \ref{cor:max_blue_degree}, every vertex $a \in A$ is the endpoint of a red path on $2p - 1$ vertices. If $p \leq k \leq m + p - 1$, then
    \begin{equation*}
        \begin{split}
            |B| = (r - 1) - |A| &= \max(2n + 2p -2, n + 2m + 3p -2) - |A| \\
            &= \max(n + p -1 - k, 2m + 2p -1 - k) \\
            & \geq m + p.
        \end{split}
    \end{equation*}
    It follows from Lemma \ref{lem:min_red_blue_degree} that $\text{red-}d(v) \geq m + p$ for all $v \in V$. Now, suppose that $\text{red-}d_B(x) \geq m - (k - p)$ for some vertex $x \in A$. Then, let $J$ be any subset of $R_A(x)$ with $|J| = k - p$. Because every vertex $a \in A$ has $\text{red-}d_A(a) \geq p- 2 + k$, we have that $\text{red-}d_{A \setminus J}(a) \geq (p - 2 + k) - (k - p) = 2p - 2$. By Lemma \ref{lem:min_path_degree}, the red subgraph induced by $A \setminus J$ admits a red path on $2p - 1$ vertices with endpoint $x$. We can then complete a red $S_{2p}(n, m)$ as follows. First, connect the other endpoint of our red path in $A \setminus J$ to $u$, extending it to $2p$ vertices. Then, use $x$ as the center of our $m$-star, with at most $k - p$ spokes in $J$ and at least $m - (k - p)$ spokes in $B$ (i.e.\ exhaust $J$ first before using $B$). So far, we have used at most $|J| + (2p - 1) = (k - p) + (2p - 1) = p - 1 + k$ vertices in $A$, leaving behind at least $|A| - (p - 1 + k) = n$ unused vertices. As a result, we can use $u$ as the center of our $n$-star, with these remaining vertices in $A$ as its spokes.

    Therefore, we may assume that $\text{red-}d_B(a) \leq m - (k-p) - 1$ for all $a \in A$. Since $|B| = \max(n + p - 1 - k, 2m + 2p - 1 - k)$, it follows that 
    \[\text{blue-}d_B(a) \geq  \max(n - m, m + p)\] for all $a \in A$. Now, suppose there exists a vertex $b \in B$ such that $\text{red-}d_B(b) \geq m - (k - p + 1)$, and let $a^\prime$ be any red neighbor in $R_A(b)$. We can then form a red $S_{2p}(n, m)$ in similar fashion to above. We first choose a subset $J^\prime \subset R_A(b) \setminus \{a^\prime\}$ with $|J^\prime| = k - p + 1$. Since $\text{red-}d_A(a) \geq p - 2 + k$ for all $a \in A$, we have that $\text{red-}d_{A \setminus J^\prime}(a) \geq (p - 2 + k) - (k - p + 1) = 2p - 3$. By Lemma \ref{lem:min_path_degree}, the red subgraph induced by $A \setminus J^\prime$ contains a red path on $2p - 2$ vertices with endpoint $a^\prime$. 
    We then extend this red path to $2p$ vertices by connecting $b$ to $a^\prime$, following the path from $a^\prime$ to its other endpoint, and connecting this endpoint to $u$. As before, we use $b$ as the center of our $m$-star, with at most $k - p + 1$ spokes in $J^\prime$ (which we exhaust first) and $m - (k - p + 1)$ spokes in $B$. Once again, we have used at most $|J^\prime| + (2p - 2) = (k - p + 1) + (2p - 2) = p - 1 + k$ vertices in $A$, leaving behind precisely $|A| - (p - 1 + k)$ vertices in $A$ to complete the $n$-star at $u$.

    Thus, we may assume that $\text{red-}d_B(b) \leq m - (k - p + 1) -1 $ for all $b \in B$, and so \[\text{blue-}d_B(b) \geq |B \setminus \{ b \}| - (m - (k - p + 1)  -1 ) = \max(n - m, m + p).\]
    But we can now find a vertex in $B$ with large blue degree, contradicting the maximality of $u$. To see this, observe that $\text{blue-}d_B(a) \geq \max(n-m, m+ p)$ for all $a \in A$, so there are at least $|A|\cdot \max(n-m, m+p)$ blue edges between $A$ and $B$. As a result, there exists a vertex $b^\prime \in B$ with $$\text{blue-}d_A(b^\prime) \geq \frac{|A|\cdot \max(n-m, m+p)}{|B|} = \frac{(n +p - 1 + k) \cdot \max(n-m, m+p)}{\max(n + p -1 - k, 2m + 2p - 1- k)}.$$ 
    Since $\text{blue-}d_B(b^\prime) \geq \max(n - m - 1, m + p -1)$, we have that
    \begin{equation}
        \label{eq:blue-d}
        \begin{split}
            \text{blue-}d(b^\prime) &= \text{blue-}d_A(b^\prime) + \text{blue-}d_B(b^\prime) + |\{ u \}| \\
            & \geq \frac{(n +p - 1 + k) \cdot \max(n-m, m+p)}{\max(n + p -1 - k, 2m + 2p - 1- k)} + \max(n - m - 1, m + p -1) + 1.
        \end{split}
    \end{equation} Let us assume, for the purpose of contradiction, that $n \geq 2m + p.$ Then, (\ref{eq:blue-d}) above simplifies to \begin{equation}
        \label{eq:blue-d2}
        \text{blue-}d(b^\prime) \geq \frac{(n + p - 1 + k)(n-m)}{n + p -1 - k} + n - m.
    \end{equation}
    Because $n \geq 2m + p$, $(n + p - 1)(2m + p - n - 1) \leq 0$. Since $k^2 \geq 1$, it follows that
    \[(n + p - 1)(2m + p - n - 1) < k^2.\]
    It is a matter of formality to show, upon rearranging, that
    \[\frac{(n + p - 1 + k)(n-m)}{n + p -1 - k} + n - m > n+ p - 1 + k.\]
   Combining this inequality with (\ref{eq:blue-d2}) above gives
     $\text{blue-}d(b^\prime) > n + p - 1 + k = |A|$. Hence, $b^\prime$ contradicts the maximality of $u$, so our initial assumption that $n \geq 2m + p$ is false.

    Therefore, $n < 2m + p$. In this regime, (\ref{eq:blue-d}) simplifies to 
    \begin{equation}
    \label{eq:blue-d3}
        \text{blue-}d(b^\prime) \geq \frac{(n + p - 1 + k)(m +  p)}{2m + 2p -1 - k} + m + p. 
    \end{equation}
    Since $k \leq m + p - 1$, we have that $k(2m + p - n) \leq (m + p - 1)(2m + p - n)$ and thus $k(2m + p - n) + 1 \leq (m + p - 1)(2m + p - n) + 1$. It is not hard to see, then, that
    \[k(2m + p - n) + 1 < k^2 +  (m + p)(2m + p - n) + n + p.\]
    But rearranging this inequality yields
    \[\frac{(n + p - 1 + k)(m +  p)}{2m + 2p -1 - k} + m + p > n + p - 1 + k,\]
   so $\text{blue-}d(b^\prime) > |A|$ by (\ref{eq:blue-d3}), a contradiction. Therefore, our assumption that $\text{red-}d_B(b) \leq m - (k - p + 1) - 1$ for all $b \in B$ is false. We conclude that if our $2$-coloring of $K_r$ does not have a polarized pair, then it admits a  monochromatic $S_{2p}(n, m)$. The lemma statement follows immediately.
\end{proof}

\begin{lemma}
    \label{lem:k_geq_m_plus_p}
    If $k \geq m +p$, then our 2-coloring of $K_r$ either admits a polarized pair or a monochromatic $S_{2p}(n,m)$.
\end{lemma}

\begin{proof}
    If $k \geq m + p$, then $|A| \geq n + m + 2p - 1$. By Corollary \ref{cor:max_blue_degree}, we know that every vertex in $A$ is the endpoint of a red path on $2p-1$ vertices. Choose any such path, and let $a \in A$ be one of its endpoints. By Lemma \ref{lem:max_blue_degree}, we have that $\text{blue-}d_A(a) \leq n$, and so $\text{red-}d_A(a) \geq n + m + 2p - 2 - n = m + 2p -2$. Thus, $a$ is connected to at least $m$ other vertices in $A$ not on the path, forming a red broom $S_{2p-1}(m,0)$ in $A$. Connecting the tail of this broom to $u$ completes a red $S_{2p}(n,m)$.
\end{proof}

\begin{lemma}[Due to Fan, Hong, and Liu, \cite{Spiders}]\label{Extremal bund on brooms}
    If a graph $G$ on $N$ vertices does not contain the broom $S_t(m,0)$, then $$e(G) \leq \frac{N\cdot (m +t-2)}{2}.$$
    \end{lemma}
\begin{remark}
    Lemma \ref{Extremal bund on brooms} asserts that the Erd\H{o}s--S\'os conjecture is true for broom graphs. In their paper, Fan, Hong, and Liu actually prove a more general statement: that the Erd\H{o}s--S\'os conjecture is true for all \emph{spiders}, trees with at most one vertex of degree at least 3 and all others with degree at most 2 \cite{Spiders}. It is not hard to see that broom graphs are spiders, and hence  Lemma \ref{Extremal bund on brooms} follows.
\end{remark}

\begin{lemma}\label{lem:wrsl_helper}
    Consider some fixed 2-coloring of $K_r$ with $r = \max(2n + 2p -1, n + 2m + 3p -1)$ and define $u$, $A$, and $B$ as usual, with $|A| = n + p - 1 + k$. Suppose that $1 \leq k \leq p-1$, and let $t = p-k$. In addition, suppose there exists a vertex $a' \in A$ with $\text{red-}d_B(a') \geq m+ t$, and define $S = B \cap R(a')$. Finally, suppose that there is no polarized pair in this 2-coloring. Then, if there is a red path $P$ in $A \cup B$ satisfying the following properties:
    \begin{itemize}
        \item $P$ starts at a vertex of $S$,
        \item $P$ contains $t$ vertices in $B$,
        \item $P$ uses no internal edges in $A$,
        \item $P$ does not use $a'$,
    \end{itemize} then our $2$-coloring admits a red $S_{2p}(n,m)$.
\end{lemma}

\begin{proof}
    We will modify $P$ to achieve a more useful path. First, observe that since $S \subseteq R(a')$ and $P$ does not use $a'$, we can append $a'$ to the endpoint of $P$ belonging to $S$. Now consider the opposite endpoint of $P$; if this endpoint is a vertex in $A$, delete it. This leaves us with a modified path $P'$ such that
    \begin{itemize}
        \item $P'$ has one endpoint in $B$, and one endpoint equal to $a'$,
        \item $P'$ contains $t$ vertices in $B$,
        \item $P'$ uses no red edges inside $A$.
    \end{itemize}
    Since $P'$ starts in $B$, ends at $a' \in A$, and contains no edges inside $A$, it must be the case that $P'$ uses at most $t$ vertices of $A$ (in the case where $P'$ simply alternates between $A$ and $B$). Let $r_A = |P' \cap A|$, and let $z$ be the endpoint of $P'$ in $B$. By assumption, our $2$-coloring of $K_r$ does not admit a polarized pair, so $\text{blue-}d_A(z) \leq n$ and thus $\text{red-}d_A(z) \geq p + k - 1$. Since $P'$ uses $r_A$ vertices in $A$, the vertex $z$ has at least $p + k - 1- r_A$ red neighbors in $A \setminus P'$. Furthermore, after deleting the $r_A$ vertices of $A$ used by $P'$, it follows from Lemma \ref{lem:max_blue_degree} that  every remaining vertex of $A$ still has red degree at least $p + k - 2 - r_A$ into $A \setminus P'$.
    
    Now, consider the red subgraph induced by $A \setminus P^\prime$. By Lemma \ref{lem:min_path_degree}, there is a red path $Q \subset A \setminus P^\prime$ on $p + k - 1 - r_A$ vertices that begins at some vertex in the red neighborhood $R_{A \setminus P^\prime}(z)$. Therefore, we can form a path by connecting $P'$ to $Q$ via $z$, and connecting the terminal endpoint of $Q$ to $u$. This forms a red path on \[|P^\prime|+ |Q| + |\{ u \}| = (t+r_A) + (p+k-1 -r_A) + 1 = p+k+t = 2p\] vertices, with $u$ and $a'$ as endpoints. This path uses $|Q| + r_A = p+k-1$ vertices in $A$, leaving exactly $n$ vertices to form the $n$-star at $u$. The other endpoint $a'$ has $\text{red-}d_B(a') \geq m + t$; but since our path only uses $t$ vertices in $B$, there are $m$ unused vertices in $B$ to complete the $m$-star at $a^\prime$. As a result, this path defines a red $S_{2p}(n,m)$, as desired.
\end{proof}

\begin{lemma}\label{lem:wrsl}
     Let $H$ be a graph with vertex partition \(A\sqcup B\), and assume that there are no edges inside \(A\). Let \(S\subseteq B\). Suppose $t \geq 2$ and that there is no path in \(H\) which starts at a vertex of \(S\) and contains at least \(t\) vertices of \(B\). Let $\mathcal{C}$ be the union of all connected components of $H$ that share vertices with $S$, and write
\[
A_\mathcal{C}=A\cap V(\mathcal{C}),
\qquad
B_\mathcal{C}=B\cap V(\mathcal{C}).
\]
Then,
\begin{equation}
\label{eq:12}
e(A_\mathcal{C},B_\mathcal{C})+2e(B_\mathcal{C})
\le
(t-1)|A_\mathcal{C}|+(3t-6)|B_\mathcal{C}|-(t-2)|S|.
\end{equation}
\end{lemma}
\begin{proof}
    It suffices to prove the claim for each connected component in $\mathcal{C}$, and then add their respective inequalities together to produce the overall result. Let $C \in \mathcal{C}$, meaning that $C$ is a connected component of $H$ and the subset $S_C = S \cap V(C) \neq \emptyset$. We define the \emph{capacity function} $c : V \rightarrow \mathbb{N}$ as follows:
    \[c(v) = \begin{cases}
         t-1 & v \in A,\\
         3t - 6 & v \in B\setminus S_C,\\
         2t - 4 & v \in S_C.
    \end{cases}\]
    By applying induction over $|V(C)|$, we will prove that 
    \[e(A_C, B_C) + 2e(B_C) \leq \sum_{v \in V(C)} c(v).\] It is not hard to see that this is equivalent to (\ref{eq:12}). For each vertex $v \in V(C)$, let $w(v)$ be the \emph{contribution} in weight that $v$ makes to $e(A_C, B_C) + 2e(B_C)$. In particular,
    \[w(v) = \begin{cases}
         d(v) & v \in A,\\
         d_{A_C}(v) + 2d_{B_C}(v) & v \in B,
    \end{cases}\]
    where we implicitly assign edges from $A_C$ to $B_C$ weight 1 and edges internal to $B_C$ weight 2. The induction is as follows.
    \begin{enumerate}[label = (\alph*)]
        \item \underline{Base Cases}:\\[6pt]
        If $|V(C)| = 0$, then (\ref{eq:12}) is trivial, since $|A_C| = |B_C| = |S_C| = 0$, so both the left- and right-hand sides of the inequality equal 0.  If $|V(C)| = 1$, then the left-hand side of (\ref{eq:12}) equals 0. For the right-hand side, we know that $|A_C| = 1$ and $|B_C| = 0$, or vice versa; thus, $(t - 1)|A_C| + (3t - 6)|B_C| - (t - 2)|S_C|$ equals $t - 1 > 0$ (if $|A_C| = 1$) or $(3t - 6) - (t - 2)|S_C| \geq 2t - 4 \geq 0$ (if $|S_C| \leq |B_C| = 1$). The right-hand side is greater than or equal to 0 in either case, so  \[e(A_C, B_C) + 2e(B_C) = 0 \leq (t - 1)|A_C| + (3t - 6)|B_C| - (t - 2)|S_C|,\] as desired.\\
        \item \underline{Induction Step}:\\[6pt]
         For the induction step, we will show that we can always identify a vertex whose weighted degree is bounded above by its capacity. We begin by choosing $P =v_0 -v_1 - \cdots - v_\ell$ to be a certain maximal path which starts in $S_C$ (i.e.\ $v_0 \in S_C$): among the paths that start in $S_C$, $P$ first maximizes the number $\beta$ of its $B$-vertices, and then maximizes its total number of vertices. By our original assumption, $P$ has $\beta \leq t - 1$ vertices that belong to $B$. Now, let $x = v_\ell$ be the terminal endpoint of $P$. By the maximality of $P$, every neighbor of $x$ lies on $P$. If $x$ had a neighbor $y$ not on $P$, then $y$ would either (1) increase the number of $B$-vertices on $P$, or (2) preserve the number of $B$ vertices but make $P$ longer. By our construction of $P$, either occurrence would contradict our choice of $P$ as the maximal path. We consider three cases, where $x \in A$, $x \in B\setminus S_C$, and $x \in S_C$, and show that in each case we can identify a vertex whose weighted degree is bounded above by its capacity.

          First, suppose that $x \in A$. Since $A$ has no internal edges, all of $x$'s neighbors belong to $B$ and lie on the path $P$. Therefore,
         \begin{equation}
             w(x) = d(x) \leq |B \cap V(P)| = \beta \leq t - 1 = c(x),
         \end{equation}
         so $w(x) \leq c(x)$ when $x \in A$. In other words, $c(x)$ is bounded below by the contribution that $x$ makes to $e(A_C, B_C) + 2e(B_C)$---namely, $d(x)$. 
         
         When $x \in B\setminus S_C$, then both endpoints $v_0$ and $x = v_\ell$ of $P$ belong to $B$. As a result, $x$ has at most $\beta - 1$ neighbors in $B$ and at most $\beta -1$ neighbors in $A$. Hence,
         \begin{equation}
             w(x) = d_{A_C}(x) + 2d_{B_C}(x) \leq (\beta - 1) + 2(\beta- 1) = 3\beta - 3 \leq 3t - 6 = c(x),
         \end{equation}
         so $w(x) \leq c(x)$ when $x \in B\setminus S_C$. 
         
         The third case, when $x \in S_C$, is the most complex. For each neighbor $v_i \in N_C(x)$\footnote{By definition, the neighbors of $x$ belong to the connected component $C$, so $N_C(x)$ is the same as $N(x)$.} (all of which belong to $P$), form the path $P_i$ by rotating $P$ around the edge $\{x, v_i\}$; that is,
         \[P_i = v_0 -v_1 - \cdots - v_i -x - v_{\ell - 1} - v_{\ell - 2} - \cdots - v_{i + 1}.\]
         When $i = \ell - 1$, we take $P_i$ to simply equal $P$. Note that the path above is valid, since $x$ is connected to $v_i$ and $v_{\ell - 1}$, and all other edges remain unchanged from $P$. We observe that $P_i$ starts in $S_C$, has the same number $\beta$ of $B$-vertices as $P$, and has the same length as $P$. The path $P_i$ has terminal endpoint $v_{i + 1}$, the successor of $v_i$. If $v_{i + 1}\notin S_C$ for some neighbor $v_i \in N_C(x)$, then the previous two cases show us that $w(v_{i + 1}) \leq c(v_{i + 1})$, so $v_{i + 1}$ satisfies our desired constraint that its weighted degree is at most its capacity.  Let us now assume that $x \in S_C$ and that $v_{i + 1} \in S_C$ for all neighbors $v_i \in N_C(x)$. The function
         \[f(v_i) = v_{i + 1}\]
         is injective from $N_C(x)$ to $V(P) \cap (B \setminus \{v_0\})$, so $d_C(x) \leq \beta - 1$, and $2d_C(x) \leq 2(\beta - 1)$. But $w(x) = d_{A_C}(x) + 2d_{B_C}(x) \leq 2d_C(x) $ and $2(\beta - 1) \leq 2(t - 2) = 2t - 4$. Therefore,
         \[w(x) = d_{A_C}(x) + 2d_{B_C}(x)  \leq 2t - 4 = c(x)\]
         when $x \in S_C$, as desired. In all three cases, we see that we can identify a vertex (either $x$ or some successor $v_{i + 1}$) whose weighted degree is at most its capacity.
         Without loss of generality, delete $x$ (the alternative being $v_{i + 1}$) and all edges incident to $x$ from $C$. Note that $C \setminus \{x\}$ is still connected, since $x$ was chosen to be an endpoint of the path $P$ (resp.\ $P_i$ when we choose $v_{i + 1}$) and all of $x$'s neighbors are still connected via $P \setminus \{x\}$ (resp.\ $P_i \setminus \{v_{i + 1}\}$). As a result,
         \[e\big(A_{C \setminus \{x\}}, B_{C \setminus \{x\}}\big) + 2e\big(B_{C \setminus \{x\}}\big) \leq \sum_{v \in V(C \setminus \{x\})} c(v).\]
         Adding back $x$ gives us that
         \[w(x) + e\big(A_{C \setminus \{x\}}, B_{C \setminus \{x\}}\big) + 2e\big(B_{C \setminus \{x\}}\big) \leq c(x) +  \sum_{v \in V(C \setminus \{x\})} c(v),\]
         so \[e(A_{C}, B_{C}) + 2e(B_{C}) \leq \sum_{v \in V(C)} c(v),\]
         as desired. In the last step, we are able to consolidate 
        $w(x) + e\big(A_{C \setminus \{x\}}, B_{C \setminus \{x\}}\big) + 2e\big(B_{C \setminus \{x\}}\big)$ into $e(A_{C}, B_{C}) + 2e(B_{C})$ precisely because $w(x)$ is the contribution of $x$ to  $e(A_{C}, B_{C}) + 2e(B_{C})$.
    \end{enumerate}
    Therefore, $e(A_C, B_C) + 2e(B_C) \leq \sum_{v \in V(C)} c(v)$ for all connected components $C$ of $H$ that intersect $S$. Expanding the right-hand side and summing over all such components $C \in \mathcal{C}$ gives us (\ref{eq:12}), which is the desired result.
\end{proof}

\begin{lemma}
\label{lem:k_1_pminus1}
    Suppose $1 \leq k \leq p-1$, $n \geq 2m + p$, and $m \geq p - 3$. Then, our 2-coloring of $K_r$ either admits a polarized pair or a monochromatic $S_{2p}(n,m)$.
\end{lemma}

\begin{proof}
    Fix $t = p-k$. Since $1 \leq k \leq p-1$, it follows that $1 \leq t \leq p-1$. We also have that $|A| = n + 2p -t -1$ and $|B| = n + t -1$.

    First, we will show that if $B$ contains a red broom $S_t(m, 0)$, then our 2-coloring of $K_r$ admits a red $S_{2p}(n,m)$. Since $|A| = n + 2p -t -1$, Corollary \ref{cor:max_blue_degree} gives us that every vertex $a \in A$ is the endpoint of a red path on \[p+k-1 = 2p - t -1\]
    vertices. Let $x$ be the tail of our red broom $S_t(m, 0)$. We know that $x$ must be adjacent to some vertex $y \in A$ by a red edge; otherwise, $u$ and $x$ would form a polarized pair. As a result, we can follow our broom to its tail $x$, connect $x$ to $y$, and then follow the red path starting at $y$ for $2p - t -1$ vertices in $A$; this gives a red path on $2p-1$ vertices. Connecting the endpoint of this path to $u$ completes a red $S_{2p}(n,m)$, as there are exactly $n$ vertices left in $A$ to complete the $n$-star.

    We may thus assume that $B$ does not admit a red $S_t(m,0)$. By Lemma \ref{Extremal bund on brooms}, the number $e_R(B)$ of red edges in $B$ therefore satisfies 
    \[e_R(B) \leq \frac{|B| \cdot (m+t -2)}{2}.\]
    Now, consider the set of ordered pairs $(b, v_i)$ where $b \in B$ and the edge $\{b, v_i\}$ is red. We observe that \[\big|\big\{(b, v_i) : b \in B,\ v_i \in V,\text{ and }\{b, v_i\} \text{ is red}\big\}\big| = e_R(A, B) + 2e_R(B),\]
    since the ordered pairs double-count the red edges with both endpoints in $B$. Simultaneously,
    \[\big|\big\{(b, v_i) : b \in B,\ v_i \in V,\text{ and }\{b, v_i\} \text{ is red}\big\}\big| \geq |B|^2,\]
    since we have $|B|$ choices for the left endpoint $b \in B$, and every vertex in $B$ has red degree at least $|B|$ (giving us at least $|B|$ choices for $v_i$). Therefore, 
    \begin{equation}
        \label{eq:countingtrick}
        e_R(A,B) + 2e_R(B) \geq |B|^2.
    \end{equation} Substituting the extremal bound on $e_R(B)$ above, we find that
    \begin{equation*}
        \begin{split}
            e_R(A,B) \geq |B|^2 - 2e_R(B)  &\geq |B|^2 - (m + t -2)|B|\\
            &= |B|(|B| - m - t + 2) \\
            &= (n + t - 1)(n - m + 1),
        \end{split}
    \end{equation*} using the fact that $|B| = n + t -1$. 
    
    We will now use this lower bound to show that there exists some vertex $a^\prime \in A$ such that $\text{red-}d_B(a^\prime) \geq m + t$. Suppose to the contrary that for all $a \in A$, $\text{red-}d_B(a) \leq m + t - 1$ and so $e_R(A,B) \leq |A|(m + t - 1) = (n + 2p - t -1)(m + t - 1)$. Combining this with the lower bound on $e_R(A, B)$ above,
    we have that \[(n + t - 1)(n - m + 1) \leq e_R(A,B) \leq (n + 2p - t -1)(m + t - 1).\] Subtracting the right-hand side from the left gives us the expression \[(n + t - 1)(n - m + 1) - (n + 2p - t -1)(m + t - 1).\] It is straightforward to see that this difference is a strictly increasing function in $n$, say, by evaluating its partial derivative: $$\frac{\partial}{\partial n} \Big[ (n + t - 1)(n - m + 1) - (n + 2p - t -1)(m + t - 1) \Big] = 2n -2m + 1 > 0,$$ using the assumption that $n \geq 2m + p$. At the specific value $n = 2m + p$, the difference above simplifies to $$4m + (p-t)^2 + 3p + t -2 > 0.$$ Since the difference function is positive at $n = 2m + p$ and is increasing, the lower bound on $e_R(A,B)$ exceeds the upper bound when $n \geq 2m + p$, a contradiction. Thus, our initial assumption must have been false, so there exists a vertex $a' \in A$ with $\text{red-}d_B(a') \geq m + t$.

    Now, choose some $a' \in A$ with this property, and let $S = R(a') \cap B$. By our choice of $a^\prime$, we have that $|S| \geq m + t$. We will now introduce an auxiliary graph $H$, the graph with vertex set $A \sqcup B$ whose edges are the red edges inside $B$ and the red edges between $A$ and $B$ (i.e.\ $H$ has no red edges inside $A$). 
    
    By Lemma \ref{lem:wrsl_helper}, we may assume that there is no red path in $H\setminus \{a'\}$ that starts at a vertex of $S$ and contains at least $t$ vertices of $B$. Let $\mathcal{C}$ be the union of the connected components of $H \setminus \{ a' \}$ that share vertices with $S$. 
    
    Suppose that $1 \leq k \leq p-2$ and thus $2 \leq t \leq p-1$; we will handle the case where $k = p - 1$ separately. Lemma \ref{lem:wrsl} gives us that \[e(A_\mathcal{C}, B_\mathcal{C}) + 2e(B_\mathcal{C}) \leq (t-1)|A_\mathcal{C}| + (3t-6)|B_\mathcal{C}| - (t-2)|S|,\] where  $A_\mathcal{C} = A \cap \mathcal{C}$ and $B_\mathcal{C} = B \cap \mathcal{C}$. We observe that the red degree of any vertex in $B$ is unchanged in $H$, since $H$ only removes edges that are internal in $A$, and there are only blue edges between $u$ and $B$. Therefore, $\text{red-}d_H(b) \geq |B|$ for all $b \in B_\mathcal{C}$. Deleting $a'$ removes at most one incident red edge for each $b$ (i.e.\ $\{a^\prime, b\})$, so $\text{red-}d_{H\setminus\{a'\}}(b) \geq |B|-1$ for all $b \in B_\mathcal{C}$. A counting strategy similar to (\ref{eq:countingtrick}) gives us that \[e(A_\mathcal{C}, B_\mathcal{C}) + 2 e(B_\mathcal{C}) \geq (|B|-1)|B_\mathcal{C}|.\] We now proceed as before. Combining our upper and lower bounds for $e(A_\mathcal{C}, B_\mathcal{C}) + 2e(B_\mathcal{C})$ yields the inequality \[(|B|-1)|B_\mathcal{C}| \leq e(A_\mathcal{C}, B_\mathcal{C}) + 2e(B_\mathcal{C}) \leq (t-1)|A_\mathcal{C}| + (3t-6)|B_\mathcal{C}| - (t-2)|S|.\] Rearranging, \[(|B| - 3t + 5)|B_\mathcal{C}| \leq (t-1)|A_\mathcal{C}| - (t-2)|S|.\] 
    We now observe that $|B_\mathcal{C}| \geq |S| \geq m+ t$, and $|A_\mathcal{C}| \leq |A|-1$. Substituting above, \[(|B|-3t + 5)(m+ t) \leq (t-1)(|A|-1) - (t-2)(m+t),\] so \[(|B|-2t + 3)(m + t) \leq (t-1)(|A|-1).\] Since  $|A| = n + 2p - t -1 $ and $|B| = n + t -1$, \[(n -t+2)(m+t) \leq (t-1)(n+2p-t-2).\] But now consider the difference between the left- and the right-hand sides of the inequality above: $$(n -t+2)(m+t)- (t-1)(n+2p-t-2).$$ Evaluating the partial derivative once again, \[\frac{\partial}{\partial n} \Big[ (n -t+2)(m+t)- (t-1)(n+2p-t-2) \Big] = m+1 > 0.\] At the specific value $n = 2m +p$, the difference function simplifies to \[2m^2 + m(p-t+4) -p(2t - 3) + 3t - 2,\] which is positive when $m \geq p - 3$. This inequality is impossible, so we have a contradiction. This closes the case where $2 \leq t \leq p-1$ (equivalently $1 \leq k \leq p-2$).

    Now, suppose $k= p-1$ (i.e.\ $t=1$), so that $|A| = n + 2p-2$. Since $t=1$, there exists a vertex $a' \in A$ such that $\text{red-}d_B(a')\geq m + t = m+1$. Choose any vertex $b \in S$, and recall that $S = B \cap R(a')$. Since $\text{red-}d_A(b) \geq p + k -1 $, we know that $b$ is connected to some vertex $a \in A\setminus \{ a' \}$ by a red edge. By Lemma \ref{lem:max_blue_degree}, we have that $\text{blue-}d_A(a) \leq n$ and thus $\text{red-}d_{A \setminus \{a^\prime\}}(a) \geq |A \setminus \{a^\prime\}| - n - 1 = (n + 2p - 3) - n - 1 = 2p - 4$ for all $a \in A$. It follows from Lemma \ref{lem:min_path_degree} that every vertex in $A \setminus \{ a'\}$ is the endpoint of a red path on $(2p - 4) + 1 = 2p -3$ vertices within $A \setminus \{ a'\}$. Let $P$ be such a red path in $A\setminus \{a'\}$ with endpoint $a$. We form a red $S_{2p}(n, m)$ as follows: connect $a^\prime \rightarrow b \rightarrow a$, and follow the path $P$ to its other endpoint; connecting this endpoint to $u$ completes a red path on $2p$ vertices. This path uses exactly $2p - 2$ vertices in $A$, leaving behind $n$ unused vertices for our $n$-star centered at $u$. Similarly, the path uses exactly one vertex in $S$, leaving behind at least $m$ unused vertices in $S$. We then use $a^\prime$ as the center of our $m$-star, forming a red $S_{2p}(n,m)$, as desired.

    Therefore, if $1 \leq k \leq p-1$ and our 2-coloring of $K_r$ has no polarized pair, it must admit a monochromatic $S_{2p}(n,m)$.
\end{proof}

\begin{lemma}
    \label{lem:k_equals_0}
    If $k = 0$ and $n \geq 2m + p$, then our 2-coloring of $K_r$ either admits a polarized pair or a monochromatic $S_{2p}(n,m)$.
\end{lemma}

\begin{proof}
    If $k = 0$, then $|A| = |B| = n + p - 1$.

    As before, we know that every vertex $b \in B$ is connected to some $a \in A$ by a red edge; if there exists a vertex $b$ without such a red edge, then $b$ and $u$ would form a polarized pair, as desired. If $B$ contains a red $S_{p}(m,0)$, then we can form a red $S_{2p}(n, m)$ as follows: connect the tail vertex of the broom to some vertex $a \in A$; by Corollary \ref{cor:max_blue_degree}, we can follow $a$ for $p - 1$ vertices in $A$ (including $a$); connect the terminal vertex of this path to $u$, completing a red $S_{2p}(n,m)$, with $u$ as the center of our $n$-star.

    Thus, we may assume that $B$ contains no red $S_{p}(m,0)$. By Lemma \ref{Extremal bund on brooms}, we know that \[e_R(B) \leq \frac{|B| \cdot (m + p -2)}{2}.\] Since $\text{red-}d(b) \geq |B|$ for all $b \in B$, we have that \[e_R(A,B) \geq |B|^2 - 2e_R(B) \geq |B|^2 - (m + p - 2)|B|,\]
    by similar reasoning to Lemma \ref{lem:k_1_pminus1}. Substituting $|B| = n + p -1$ gives us that \[e_R(A,B) \geq |B|(n - m + 1) = (n +p-1)(n-m+1).\]

    When $k=0$, the red and blue graphs are ``symmetric.\rq\rq\  In particular, by a very similar argument to the one above, $A$ does not admit a \emph{blue} $S_{p}(m,0)$; otherwise, we could complete a blue $S_{2p}(n,m)$. Repeating the argument above, we have that \[e_B(A,B) \geq |A|(n-m+1)= (n +p-1)(n-m+1),\] so \[e_R(A,B)+ e_B(A,B) \geq 2(n + p -1)(n - m + 1).\] But since any edge between $A$ and $B$ is either red or blue, we know that $e_R(A,B)+ e_B(A,B) = |A|\cdot|B| = (n + p -1)^2$. Hence, it must be the case that $(n + p -1)^2 \geq 2(n + p -1)(n-m+1)$. Solving for $n$, we find that \[n \leq 2m + p -3,\] contradicting the assumption that $n\geq 2m +p$.
    
\end{proof}

\begin{remark}
    Lemma \ref{lem:k_1_pminus1_edge_case} below handles the case where $n = 2m + p - 1$. Note that there is no $k = 0$ case here, by the following reasoning. When $n = 2m + p - 1$, we have that $r = \max(2n + 2p - 1, n + 2m + 3p - 1) = \max(2n + 2p - 1, 2n + 2p) = 2n + 2p$. Removing the maximal vertex $u$ leaves us with $2n + 2p - 1$ vertices, which is odd. Therefore, $|A|$ is strictly larger than $|B|$. But when $k = 0$, it follows that $|A| = n + p - 1$ and $|B| = 2n + 2p - 1 - |A| = n + p$, so $|B| > |A|$, a contradiction. For this reason, the lemma below only considers $1 \leq k \leq p - 1$.
\end{remark}
\begin{lemma}
\label{lem:k_1_pminus1_edge_case}
    Suppose $1 \leq k \leq p-1$, $n = 2m + p-1$, and $m \geq p -3$. Then, our 2-coloring of $K_r$ either admits a polarized pair or a monochromatic $S_{2p}(n,m)$.
\end{lemma}

\begin{proof}

    The proof of this lemma is very similar to that of Lemma \ref{lem:k_1_pminus1}. Since $n = 2m + p- 1$, \[r = \max(2n + 2p - 1, n + 2m + 3p - 1) = 2n + 2p.\] Note that $|A| = n + p -1 + k$ and $|B| = n + p -k$. Fix $t = p-k$. Since $1 \leq k \leq p-1$, it follows that $1 \leq t \leq p-1$. We have that $|A| = n + 2p -t -1$ and $|B| = n + t$.

    First, we will show that if $B$ contains a red broom $S_t(m, 0)$, then our 2-coloring of $K_r$ admits a red $S_{2p}(n,m)$. Since $|A| = n + 2p -t -1$, Corollary \ref{cor:max_blue_degree} gives us that every vertex $a \in A$ is the endpoint of a red path on \[p+k-1 = 2p - t -1\] vertices. Let $x$ be the tail of our red broom $S_t(m, 0)$. We know that $x$ must be red-adjacent to some vertex $y \in A$; otherwise, $u$ and $x$ would form a polarized pair. As a result, we can follow our broom to its tail $x$, connect $x$ to $y$, and then follow the red path starting at $y$ for $2p - t -1$ vertices in $A$; this gives a red path on $2p-1$ vertices. Connecting the endpoint of this path to $u$ completes a red $S_{2p}(n,m)$, as there are exactly $n$ vertices left in $A$ to complete the $n$-star.

    We may thus assume that $B$ does not admit a red $S_t(m,0)$. By Lemma \ref{Extremal bund on brooms}, the number $e_R(B)$ of red edges in $B$ satisfies 
    \[e_R(B) \leq \frac{|B| \cdot (m+t -2)}{2}.\]
    Now, consider the set of ordered pairs $(b, v_i)$ where $b \in B$ and the edge $\{b, v_i\}$ is red. We observe that \[\big|\big\{(b, v_i) : b \in B,\ v_i \in V,\text{ and }\{b, v_i\} \text{ is red}\big\}\big|  = e_R(A, B) + 2e_R(B),\]
    since the ordered pairs double-count the red edges with both endpoints in $B$. Simultaneously,
    \[\big|\{(b, v_i) : b \in B,\ v_i \in V,\text{ and }\{b, v_i\} \text{ is red}\big\}\big|  \geq |B|^2,\]
    since we have $|B|$ choices for the left endpoint $b \in B$, and every vertex in $B$ has red degree at least $|B|$ (giving us at least $|B|$ choices for $v_i$). Therefore, 
    \begin{equation*}
        e_R(A,B) + 2e_R(B) \geq |B|^2.
    \end{equation*} Substituting the extremal bound on $e_R(B)$ above, we find that
    \begin{equation*}
        \begin{split}
            e_R(A,B) \geq |B|^2 - 2e_R(B)  &\geq |B|^2 - (m + t -2)|B|\\
            &= |B|(|B| - m - t + 2) \\
            &= (n + t)(n - m +2),
        \end{split}
    \end{equation*} using the fact that $|B| = n + t$. 
    
    We will now use this lower bound to show that there exists some vertex $a^\prime \in A$ such that $\text{red-}d_B(a^\prime) \geq m + t$. Suppose to the contrary that for all $a \in A$, $\text{red-}d_B(a) \leq m + t - 1$ and so $e_R(A,B) \leq |A|(m + t - 1) = (n + 2p - t -1)(m + t - 1)$. Combining this with the lower bound on $e_R(A, B)$ above,
    we have that \[(n + t)(n - m + 2) \leq e_R(A,B) \leq (n + 2p - t -1)(m + t - 1).\] Subtracting the right-hand side from the left gives us the expression
    \begin{equation*}
        \begin{split}
            (n + t)(n - m + 2) - (n + 2p - t -1)(m + t - 1) &= 5 m + p^2 + p (3 - 2 t) + 2 t + t^2 - 3 \\
            &= 5m + (p -t)^2 + 3p + 2t - 3\\ &> 0,
        \end{split}
    \end{equation*}
   where above we have substituted $n = 2m + p -1 $. But our lower bound on $e_R(A,B)$ exceeds our upper bound, a contradiction. Thus, our initial assumption must have been false, so there exists a vertex $a' \in A$ with $\text{red-}d_B(a') \geq m + t$.

    Choose some $a' \in A$ with this property, and let $S = R(a') \cap B$. By our choice of $a^\prime$, we know that $|S| \geq m + t$. Let $H$ be the graph with vertex set $A \sqcup B$ whose edges are the red edges inside $B$ and the red edges between $A$ and $B$ (i.e.\ $H$ has no red edges inside $A$). 
    
    By Lemma \ref{lem:wrsl_helper}, we may assume that there is no red path in $H\setminus \{a'\}$ that starts at a vertex of $S$ and contains at least $t$ vertices of $B$. Let $\mathcal{C}$ be the union of the connected components of $H \setminus \{ a' \}$ that share vertices with $S$. 
    
    Suppose that $1 \leq k \leq p-2$ and thus $2 \leq t \leq p-1$; we will handle the case where $k = p - 1$ (i.e.\ $t=1$) separately. Lemma \ref{lem:wrsl} gives us that \[e(A_\mathcal{C}, B_\mathcal{C}) + 2e(B_\mathcal{C}) \leq (t-1)|A_\mathcal{C}| + (3t-6)|B_\mathcal{C}| - (t-2)|S|,\] where  $A_\mathcal{C} = A \cap \mathcal{C}$ and $B_\mathcal{C} = B \cap \mathcal{C}$. We observe that the red degree of any vertex in $B$ is unchanged in $H$, since $H$ only removes edges that are internal in $A$, and there are only blue edges between $u$ and $B$. Therefore, $d_H(b) \geq |B|$ for all $b \in B_\mathcal{C}$. Deleting $a'$ removes at most one incident red edge for each $b$ (i.e.\ $\{a^\prime, b\})$, so $d_{H\setminus\{a'\}}(b) \geq |B|-1$ for all $b \in B_\mathcal{C}$. A counting strategy similar to (\ref{eq:countingtrick}) gives us that \[e(A_\mathcal{C}, B_\mathcal{C}) + 2 e(B_\mathcal{C}) \geq (|B|-1)|B_\mathcal{C}|.\] We now proceed as before. Combining our upper and lower bounds for $e(A_\mathcal{C}, B_\mathcal{C}) + 2e(B_\mathcal{C})$ yields the inequality \[(|B|-1)|B_\mathcal{C}| \leq e(A_\mathcal{C}, B_\mathcal{C}) + 2e(B_\mathcal{C}) \leq (t-1)|A_\mathcal{C}| + (3t-6)|B_\mathcal{C}| - (t-2)|S|.\] Rearranging, \[(|B| - 3t + 5)|B_\mathcal{C}| \leq (t-1)|A_\mathcal{C}| - (t-2)|S|.\] 
    We observe that $|B_\mathcal{C}| \geq |S| \geq m+ t$, and $|A_\mathcal{C}| \leq |A|-1$. Substituting above, \[(|B|-3t + 5)(m+ t) \leq (t-1)(|A|-1) - (t-2)(m+t),\] so \[(|B|-2t + 3)(m + t) \leq (t-1)(|A|-1).\] Since  $|A| = n + 2p - t -1 $ and $|B| = n + t$, \[(n -t+3)(m+t) \leq (t-1)(n+2p-t-2).\] But now consider the difference between the left- and the right-hand sides of the inequality above: 

    \begin{equation*}
        \begin{split}
            (n -t+ 3)(m+t)- (t-1)(n+2p-t-2) &= 2 m^2 + m(p - t + 4) - p(2t - 3) + 4 t - 3 \\ & > 0,
        \end{split}
    \end{equation*}
    where above we use the fact that $m \geq p - 3$. This inequality is impossible, so we have a contradiction. This closes the case where $2 \leq t \leq p-1$ (equivalently $1 \leq k \leq p-2$).

    Now, suppose that $k= p-1$, so that $|A| = n + 2p-2$. Since $t=p -k = 1$, there exists a vertex $a' \in A$ such that $\text{red-}d_B(a')\geq m + t = m+1$. Choose any vertex $b \in S$, and recall that $S = B \cap R(a')$. Since $\text{red-}d_A(b) \geq p + k -1 $, we know that $b$ is connected to some vertex $a \in A\setminus \{ a' \}$ by a red edge. By Lemma \ref{lem:max_blue_degree}, we have that $\text{blue-}d_A(a) \leq n$ and thus $\text{red-}d_{A \setminus \{a^\prime\}}(a) \geq |A \setminus \{a^\prime\}| - n - 1 = (n + 2p - 3) - n - 1 = 2p - 4$ for all $a \in A$. By Lemma \ref{lem:min_path_degree}, every vertex in $A \setminus \{ a'\}$ is the endpoint of a red path on $(2p - 4) + 1 = 2p -3$ vertices within $A \setminus \{ a'\}$. Let $P$ be such a red path in $A\setminus\{a'\}$ with endpoint $a$. We then form a red $S_{2p}(n, m)$ as follows. We first connect $a^\prime \rightarrow b \rightarrow a$ and then follow the path $P$ to its other endpoint; connecting this endpoint to $u$ completes a red path on $2p$ vertices. This path uses exactly $2p - 2$ vertices in $A$, leaving behind $n$ unused vertices for our $n$-star centered at $u$. Similarly, the path uses exactly one vertex in $S$, leaving behind at least $m$ unused vertices in $S$. We then use $a^\prime$ as the center of our $m$-star, forming a red $S_{2p}(n,m)$, as desired.

    Therefore, if $1 \leq k \leq p-1$ and our 2-coloring of $K_r$ has no polarized pair, it must admit a monochromatic $S_{2p}(n,m)$.
\end{proof}

\begin{lemma}
    \label{lem:wrapper2}
    If $n \geq 2m + p - 1$, $m \geq p -3$, and $0 \leq k \leq p - 1$, then our 2-coloring of $K_r$ admits a polarized pair or a monochromatic $S_{2p}(n, m)$. 
\end{lemma}
\begin{proof}
    This is an immediate consequence of Lemmas \ref{lem:k_1_pminus1}, \ref{lem:k_equals_0}, and \ref{lem:k_1_pminus1_edge_case}.
\end{proof}

\begin{lemma}
    \label{lem:large_n_wrapper}
    If $n \geq 2m + p - 1$ and $m \geq p - 3$, then our 2-coloring of $K_r$ admits a polarized pair or a monochromatic $S_{2p}(n, m)$. 
\end{lemma}
\begin{proof}
    This follows immediately from Lemmas \ref{lem:p_leq_k_leq_m_plus_p_minus_1}, \ref{lem:k_geq_m_plus_p}, and \ref{lem:wrapper2}.
\end{proof}
\begin{lemma}\label{lem:little_helper}
    If our 2-coloring of $K_r$ contains a monochromatic $S_{2p-2}(n+2, m)$, then it either admits a polarized pair or a monochromatic $S_{2p}(n,m)$.
\end{lemma}
\begin{proof}
    Without loss of generality, suppose that our $2$-coloring of $K_r$ contains a red $S_{2p-2}(n+2, m)$. Let $u$ be the center of the $(n+2)$-star, and let $v$ be the immediate predecessor of $u$ on the path. In addition, let $X$ be the set of the $n+2$ spoke vertices red-adjacent to $u$.

    If $v$ is connected to every vertex $x \in X$ by a blue edge, then $u$ and $v$ define a polarized pair. Thus, we may assume that $v$ is connected to some vertex $x' \in X$ by a red edge. If $x'$ is red-adjacent to any other vertex $x \in X$, then the sequence $v \rightarrow x' \rightarrow x \rightarrow u$ extends the path of the $S_{2p-2}(n+2,m)$ by $2$, producing an $S_{2p}(n,m)$.

    Therefore, we may assume that $x'$ is connected to every other $x \in X$ by a blue edge. Since $|X \setminus \{ x^\prime \}| = n+1$, it follows that $u$ and $x'$ define a polarized pair.
\end{proof}

\begin{lemma}
    \label{lem:no_polarized_pair}
    If $m \geq p - 3$ and our 2-coloring of $K_r$ does not contain a polarized pair, then it admits a monochromatic $S_{2p}(n, m)$.
\end{lemma}
\begin{proof}
    We prove this result by performing a short induction over $p$. The base case $p = 3$ follows from the argument in Appendix \ref{sec:p=3}. For the inductive hypothesis, we suppose that the claim holds for some $p = i$; that is, if $p = i$, then every 2-coloring of $K_r$ (with $r = \max(2n + 2i - 1, n + 2m + 3i - 1)$) either admits a polarized pair or a monochromatic $S_{2i}(n, m)$. We will show that the claim holds for $p = i + 1$.

    Assume that $p = i + 1$, so that
    \begin{align*}
        r = \max(2n + 2p - 1, n + 2m + 3p - 1) = \max(2n + 2i + 1 , n + 2m + 3i + 2).
    \end{align*}
    We consider three cases involving the size of $n$. Note that every $2$-coloring of $K_{r^\prime}$ with $r^\prime =\max(2n + 2i + 3, n + 2m + 3i + 1)$ admits an $(n + 2)$-polarized (and therefore $n$-polarized) pair or a monochromatic $S_{2i}(n + 2, m)$ by our inductive hypothesis.
    
    If $n \leq 2m + p - 3 = 2m + i - 2$, then $r^\prime = n + 2m + 3i + 1 \leq r$. Therefore, every $2$-coloring of $K_r$ admits a monochromatic $S_{2p - 2}(n + 2, m)$.  It follows from Lemma \ref{lem:little_helper} that our $2$-coloring of $K_r$ either admits a polarized pair or a monochromatic $S_{2p}(n, m)$.
    
    The case when $n = 2m + p - 2 = 2m + i - 1$ is similar. Here, we have that $r^\prime = 2n + 2i + 3 =  r$, so once again, every 2-coloring of $K_r$ admits a monochromatic $S_{2p - 2}(n + 2, m)$. By Lemma \ref{lem:little_helper}, our $2$-coloring of $K_r$ either admits a polarized pair or a monochromatic $S_{2p}(n, m)$. 

    Finally, if $n \geq 2m + p - 1$, then Lemma \ref{lem:large_n_wrapper} gives us the desired result.

    By induction, we conclude that every $2$-coloring of $K_r$ either admits a polarized pair or a monochromatic $S_{2p}(n, m)$. The lemma statement follows immediately.
\end{proof}
Lemma \ref{lem:no_polarized_pair} above completes the no-polarized-pair-exists case. By combining with Lemma \ref{lem:polarized_pair}, we obtain our desired upper bound.
\begin{lemma}
    \label{lem:upper_bound}
    Suppose $m \geq p - 2$. Then, our $2$-coloring of $K_r$ admits a monochromatic $S_{2p}(n, m)$.
\end{lemma}
\begin{proof}
    The proof here is immediate. If our $2$-coloring of $K_r$ contains a polarized pair, then Lemma \ref{lem:polarized_pair} yields a monochromatic $S_{2p}(n, m)$. If there is no such polarized pair, then Lemma \ref{lem:no_polarized_pair} does the same. Therefore, our $2$-coloring of $K_r$ admits a monochromatic $S_{2p}(n, m)$, as desired.
\end{proof}

\section{Conclusion}\label{sec:conclusion}
Sections \ref{sec:lower_bounds} and \ref{sec:upper_bounds} give us matching upper and lower bounds for $r(S_{2p}(n,m))$ when $m \geq p -2 \geq 2$, allowing us to conclude the following result.
\begin{theorem}
    Suppose $n \geq m \geq p - 2 \geq 0$. Then,
    \[r(S_{2p}(n, m)) = \max(2n + 2p - 1, n + 2m + 3p - 1).\]
\end{theorem}
\begin{proof}
     The $p = 2$ case was solved by Burr and Erd\H{o}s \cite{BurrErdos}, and the $p = 3$ case is handled by Appendix \ref{sec:p=3}. By Lemma \ref{lem:lower_bounds}, we have that $r(S_{2p}(n, m)) \geq \max(2n + 2p - 1, n + 2m + 3p - 1)$. Likewise, for $p \geq 4$, Lemma \ref{lem:upper_bound} gives us that $r(S_{2p}(n, m)) \leq \max(2n + 2p - 1, n + 2m + 3p - 1)$. The theorem statement follows immediately.
\end{proof}
To close, we would like to offer the following open problems:
\begin{enumerate}
    \item $\text{[Full resolution.]}$ In this paper, we determine the Ramsey numbers of all ``short\rq\rq\ even-linked double stars $S_{2p}(n, m)$, that is, where $m \geq p - 2$. The ``long\rq\rq\ case, where $p \geq m + 3$, is still open. However, we simply conjecture that our formula $r(S_{2p}(n, m)) = \max(2n + 2p - 1, n + 2m + 3p - 1)$ continues to hold true in general.
    \begin{conjecture}
        \label{conj:conj1}
        The Ramsey numbers of even-linked double stars $S_{2p}(n, m)$ are given by \[r(S_{2p}(n, m)) =  \max(2n + 2p - 1, n + 2m + 3p - 1)\]
        for all $n, m, p \geq 2$.
    \end{conjecture}
    We know that this conjecture is true when $p = 2$ (see Burr and Erd\H{o}s, \cite{BurrErdos}), $p = 3$ (see Appendix \ref{sec:p=3}), and for all short even-linked double stars (i.e.\ $m \geq p - 2$, which absorbs the $p = 2$ and $p = 3$ cases). In 2026, Liu and Chen proved that the conjecture persists when $m = 2$ and either $n = 2$ or $p \geq n + 1$ \cite{LiuChen}. We also have a saturating lower bound that applies for all $n, m, p \geq 2$ (see Section \ref{sec:lower_bounds}). For these reasons, we are inclined to suspect that Conjecture \ref{conj:conj1} holds.
    \item $\text{[Multicolor Ramsey numbers.]}$
    The primary object that we use to enforce the existence of an even-linked double star is the polarized pair. Polarized pairs consist of two vertices $\pred$ and $\pblue$ whose respective red and blue neighborhoods sufficiently overlap. It seems natural to extend this definition to more colors---for instance, creating polarized \emph{triples} in red, blue, and now also green. With this idea in mind, does the argument presented in this paper generalize to more colors? Why or why not?
    \item $\text{[Odd-linked double stars.]}$ The Ramsey numbers of odd-linked double stars also have yet to be fully resolved. As noted in the introduction, Karamchedu and Klawe in 2022 proved that $r(S_{2p + 1}(n, m)) = 2n + 2m + 2p - 1$ when $n \geq 2p + 1$, and $r(S_{2p + 1}(n, 2)) = n + 3p + 3$ when $m = 2$ and $n \leq p - 2$ \cite{KaramcheduKlawe}. While these exact results do not handle all odd-linked double stars, their corresponding lower bound does. In particular, Karamchedu and Klawe prove that $r(S_{2p + 1}(n, m)) \geq \max(2n + 2m + 2p - 1, n + m + 3p + 1)$ for all $n, m, p \geq 2$. As such, we conjecture that this formula exactly determines the Ramsey numbers for all odd-linked double stars.
    \begin{conjecture}
        The Ramsey numbers of odd-linked double stars $S_{2p + 1}(n, m)$ are given by \[r(S_{2p + 1}(n, m)) =  \max(2n + 2m + 2p - 1, n + m + 3p + 1)\]
        for all $n, m, p \geq 2$.
    \end{conjecture}
    \item $\text{[Classical double stars.]}$ Finally, the Ramsey numbers of classical double stars $S_2(n, m)$ are still unresolved as well. In 1979, Grossman et al.\  proved that $r(S_2(n, m)) = \max(2n + 2, n + 2m + 2)$ when $m \geq 3$ and either $n \leq \sqrt{2}m$ or $3m \leq n$, leaving behind the unsolved gap $\sqrt{2}m < n < 3m$ \cite{GrossmanHararyKlawe}.\footnote{Grossman et al.\ also handle the $m \leq 2$ case, showing that $r(S_{2}(n, m)) = \max(2n + 1, n + 2m + 2)$ when $n$ is odd, and the general formula $\max(2n + 2, n + 2m + 2)$ when $n$ is even.} They conjectured that the same formula continues to hold for the Ramsey numbers in the gap, a claim disproved by Norin et al.\ in 2016 via a counterexample \cite{NorinSunZhao}. However, Norin et al.\ also used flag algebras to show that the conjecture still remains true for a portion of the gap, tightening the unsolved region to $1.699m < n < 3m$. The best known explicit upper bound inside this gap comes from Flores Dub\'o and Stein in 2025 \cite{FloresStein}, who show that
    \[r(S_2(n, m)) \leq \left\lceil \sqrt{3n^2 + mn + \frac{m^2}{4}}  + \frac{m}{2}\right\rceil + 1\]
    when $1.699m < n < 3m$. As such, we ask the following question:
    \begin{question}
        What are the Ramsey numbers of classical double stars $S_2(n, m)$? That is, what are their Ramsey numbers in the gap $1.699m < n < 3m$?
    \end{question}
\end{enumerate}

\section*{Acknowledgments}
We would like to extend our thanks to Jerrold W.\ Grossman of Oakland University for providing us with a copy of \cite{Grossman}. Mithra Karamchedu was supported by an Honorable Mention for the Giovanni Borrelli Mathematics Fellowship at Harvey Mudd College and by the NSF Graduate Research Fellowship under Award DGE-2437839; he would also like to thank Yuval Wigderson and the organizers of the 2025 Undergraduate Summer School at the Park City Mathematics Institute. Andrés R.\ Vindas Meléndez was partially supported by the National Science Foundation under Award DMS-2532321.
\newpage

\appendix

\section{The $p=3$ Case}\label{sec:p=3}
In this section, we will prove that 
\[r(S_{6}(n,m)) = \max(2n + 5, n + 2m + 8),\]
which represents the $p = 3$ case. The argument that follows makes use of the Ramsey numbers of two disjoint stars, which are given by Lemma \ref{lem:disjoint_stars} below.
\begin{lemma}(Due to Grossman, \cite{Grossman})
    \label{lem:disjoint_stars}
    The Ramsey numbers of disjoint stars $K_{1,s} \sqcup K_{1,t}$ are given by
    \[r(K_{1,s} \sqcup K_{1,t}) = \max(2s + 1, s + 2t, s + t + 3),\]
    where $s \geq t \geq 1$.
\end{lemma}
Now, let $r = \max(2n + 5, n + 2m + 8)$, and consider some fixed 2-coloring of $K_r$. By Lemma \ref{lem:disjoint_stars} (with $s = \max(n + 2, m + 3)$ and $t = \min(n + 2, m + 3)$), our $2$-coloring admits a monochromatic---without loss of generality, say, red---copy of $K_{1, n + 2} \sqcup K_{1, m + 3}$.
\subsubsection{Notation}
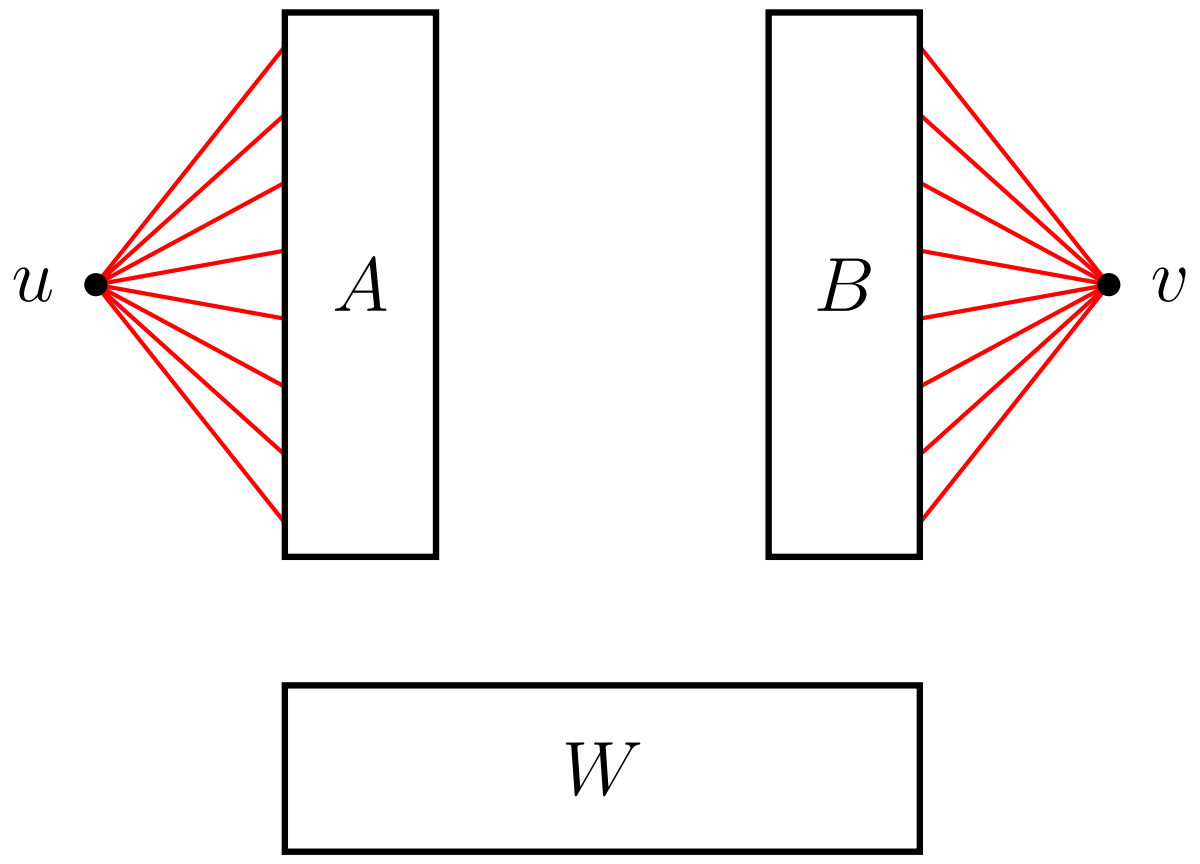
\begin{figure}[h]
\centering
\begin{tikzpicture}[
    box/.style={draw=black, very thick, fill=white},
    vertex/.style={circle, draw=black, fill=white, inner sep=1.4pt},
    rededge/.style={red, line width=0.8pt},
    every node/.style={font=\Large}
]

\def\Aleft{1.5}
\def\Aright{2.5}
\def\Bleft{4.7}
\def\Bright{5.7}

\def\BoxBottom{1.2}
\def\BoxTop{4.8}

\def\Wbottom{-0.75}
\def\Wtop{0.35}

\coordinate (u) at (0.25,3.0);
\coordinate (v) at (6.95,3.0);

\foreach \y in {1.425,1.875,2.325,2.775,3.225,3.675,4.125,4.575}{
    \draw[rededge] (u) -- (\Aleft,\y);
}

\foreach \y in {1.425,1.875,2.325,2.775,3.225,3.675,4.125,4.575}{
    \draw[rededge] (\Bright,\y) -- (v);
}

\draw[box]
    (\Aleft,\BoxBottom)
    rectangle
    (\Aright,\BoxTop);

\draw[box]
    (\Bleft,\BoxBottom)
    rectangle
    (\Bright,\BoxTop);

\node at (2.0,3.0) {$A$};
\node at (5.2,3.0) {$B$};

\draw[box]
    (\Aleft,\Wbottom)
    rectangle
    (\Bright,\Wtop);

\node at (3.6,-0.20) {$W$};

\node[vertex, fill=black] at (u) {};
\node[vertex, fill=black] at (v) {};

\node[left=4pt] at (u) {$u$};
\node[right=4pt] at (v) {$v$};

\end{tikzpicture}

\caption{The sets $A$, $B$, $W$, and the vertices $u,v$.}
\label{fig:appendix_fig}
\end{figure}

We define the following variables with respect to our chosen $2$-coloring; we will reuse them throughout the section; see Figure \ref{fig:appendix_fig}.
\begin{itemize}
    \item Let $u$ and $v$ be the centers of the $(n + 2)$- and $(m + 3)$-stars, respectively.
    \item Let $A$ contain the spokes of our $(n + 2)$-star and $B$ contain the spokes of our $(m + 3)$-star. Since \(n \geq 2\) and \(m \geq 2\), we have that \(|A| = (n + 2) \geq 4\) and \(|B| = (m + 3) \geq 5\). 
    \item Let $W = V(K_r) \setminus ( A \cup B \cup \{u,v\})$, so that $|W| = \max (m+1, n-m-2)$.
    \item Finally, let $A_0 = \{ a \in A : \text{red-}d_A(a) = 0\}$, and  let $A_1 = A \setminus A_0 = \{ a \in A : \text{red-}d_A(a) \geq 1\}$. Define the sets $B_0$ and $B_1$ analogously with respect to $B$.
\end{itemize}

We will now prove a series of lemmas that will together prove the desired result.

\begin{lemma}\label{lem:appendix_helper_1}
If the bipartite graph between $A$ and $B$ contains a red path on four vertices, then our 2-coloring of $K_r$ admits a red $S_6(n,m)$.
\end{lemma}

\begin{proof}
    Suppose there exists such a red path $a_1 - b_1 - a_2 - b_2$ between $A$ and $B$. Since $u$ is red-adjacent to $a_1 \in A$ and $v$ is red-adjacent to $b_2 \in B$, we can append $u$ and $v$ to either end of our path, extending it to six vertices: $u - a_1 - b_1 - a_2 - b_2 - v$. This path uses two vertices in $A$ and $B$ each, leaving behind $n$ vertices in $A$ and $m + 1$ vertices in $B$. Thus, we can use $u$ as the center of our $n$-star (with spokes in $A \setminus \{a_1, a_2\}$) and $v$ as the center of our $m$-star (with spokes in $B \setminus \{b_1, b_2\}$), completing a red $S_6(n, m)$.
\end{proof}
\newpage
\begin{lemma}\label{lem:appendix_helper_2}
    If there is a red edge between $A_1$ and $B_1$, then our 2-coloring of $K_r$ admits a red $S_6(n,m)$.
\end{lemma}

\begin{proof}
    Suppose there exists a red edge between any vertex $a \in A_1$ and $b \in B_1$. By the definitions of $A_1$ and $B_1$, there exist vertices $a' \in A$ and $b' \in B$ (with $a \neq a^\prime$ and $b \neq b^\prime$) such that $\{a,a'\}$ and $\{b,b'\}$ are red. Since $u$ is red-adjacent to $a^\prime$ and $v$ is red-adjacent to $b^\prime$, we can thus form the path $u - a' - a - b  - b' - v$ on six vertices. This path occupies two vertices in $A$ and $B$ each, allowing us to use $u$ as the center of our $n$-star (with spokes in $A \setminus \{a, a^\prime\}$) and $v$ as the center of our $m$-star (with spokes in $B \setminus \{b, b^\prime\}$); this completes a red $S_6(n, m)$, as desired.
\end{proof}

\begin{lemma}[Due to Faudree and Schelp, \cite{FaudreeSchelp}]
\label{lem:faudree_schelp}
    Let $G = X \sqcup Y$ be a complete bipartite graph with $\min(|X| , |Y|) \geq 3$. Then, any $2$-coloring of $G$ admits a monochromatic path on four vertices.
\end{lemma}

\begin{lemma}\label{lem:appendix_helper_3}
    If there is a blue complete bipartite graph between $A$ and $B$, then our 2-coloring of $K_r$ admits a monochromatic $S_6(n,m)$.
\end{lemma}

\begin{proof}
    Suppose there exists a blue edge between any vertex $b \in B$ and $w \in W \cup \{ u \}$. Then, we can complete a blue $S_{6}(n,m)$ as follows. First, use the blue complete bipartite graph between $A$ and $B$ to form the path
    \[b - a_1 - b_2 - a_2 - b_3 - a_3,\] which alternates between $A$ and $B$. This path occupies exactly three vertices in $A$ and three vertices in $B$. By virtue of the blue complete bipartite graph, $b$ is blue-adjacent to $n-1$ other vertices in $A$ that do not belong to the path. By similar reasoning, $a_3$ is blue-adjacent to $m$ vertices in $B$ also not on the path. Since $b$ is also blue-adjacent to $w$, we can use $b$ as the center of our $n$-star (with $n - 1$ spokes in $A \setminus \{a_1, a_2, a_3\}$ and $w$ as the last spoke) and $a_3$ as the center of our $m$-star (with spokes in $B \setminus \{b, b_2, b_3\}$), completing a blue $S_6(n, m)$.

    Therefore, we may assume that there is a red complete bipartite graph between $B$ and $W \cup \{ u \}$. But we can now form a red $S_{6}(n, m)$ as follows: choose any vertices $b_1^\prime, b_2^\prime, b_3^\prime \in B$ and $w_1, w_2 \in W$, and use the red complete bipartite graph to form the red path $u - b_1^\prime - w_1 - b_2^\prime -w_2 -b_3^\prime$; use $u$ as the center of our $n$-star (with spokes in $A$) and $b_3^\prime$ as the center of our $m$-star (with $m - 1$ spokes in $W \setminus \{w_1, w_2\}$ and $v$ as the last spoke). This completes a red $S_6(n,m)$, as desired.
\end{proof}

\begin{lemma}
    If $A_0 = \emptyset$ and $B_0 = \emptyset$, then our $2$-coloring of $K_r$ admits a monochromatic $S_6(n,m)$.
\end{lemma}

\begin{proof}
    If $A_0 = \emptyset$, it follows that $A_1 = A$. Similarly, if $B_0 = \emptyset$, then $B_1 = B$. By Lemma \ref{lem:appendix_helper_2}, we may assume that there are no red edges between $A_1 = A$ and $B_1 = B$; therefore, there is a blue complete bipartite graph between $A$ and $B$. By Lemma \ref{lem:appendix_helper_3}, our $2$-coloring admits a monochromatic $S_{6}(n,m)$, as desired.
\end{proof}

\begin{lemma}
    If $A_0 \neq \emptyset$ and $B_0 \neq \emptyset$, then our $2$-coloring of $K_r$ admits a monochromatic $S_6(n,m)$.
\end{lemma}

\begin{proof}
 If $A_0 \neq \emptyset$ and $B_0 \neq \emptyset$, then by definition, there exists a vertex $x \in A$ with $\text{red-}d_A(x) = 0$ and likewise a vertex $y \in B$ with $\text{red-}d_B(y) = 0$. It follows that $\text{blue-}d_A(x) = n + 1$ and $\text{blue-}d_B(y) = m+2$.

 First, suppose that $|A_1| \geq 2$ and $|B_1| \geq 2$. By Lemma \ref{lem:appendix_helper_2}, we may assume that there is a blue complete bipartite graph between $A_1$ and $B_1$ (otherwise, we would have a red $S_6(n,m)$). Using this complete bipartite graph, we can form the blue path \[x - a_1 - b_1 - a_2 - b_2 - y,\]
 where $a_1, a_2 \in A_1$ and $b_1, b_2 \in B_1$. Our path uses two vertices from $B \setminus \{y\}$, leaving behind $\text{blue-}d_B(y) - 2 = (m+2) - 2=m$ unused vertices in $B$; thus, we can use $y$ as the center of our $m$-star (with spokes in $B \setminus \{y, b_1, b_2\}$). Similarly, there are $\text{blue-}d_A(x) - 2 = (n + 1) - 2=n - 1$ unused vertices in $A \setminus \{x\}$; if $x$ is blue-adjacent to any vertex $w \in W \cup \{v\}$, then we can use $x$ as the center of our $n$-star (with $n - 1$ spokes in $A \setminus \{x, a_1, a_2\}$ and $w$ as the last spoke), completing a blue $S_6(n, m)$. Hence, we may assume that $x$ is red-adjacent to every vertex in $W \cup \{ v \}$. If there exists a red edge between any vertex $w \in W$ and $b \in B_1$, then the red path \[u - x - w - b - b'-v\] completes a red $S_6(n,m)$, with $n$-star centered at $u$ (with spokes in $A \setminus \{x\}$) and $m$-star centered at $v$ (with spokes in $B \setminus \{b, b^\prime\}$). Note that such a vertex $b^\prime \in B$ with $\{b, b^\prime\}$ being red must exist, since $b \in B_1$. As a result, we may assume that there is a blue complete bipartite graph between $W$ and $B_1$. However, we can now form the blue path $$x - a_1 - b_1 -w - b_2 - y$$
 for any vertex $w \in W$. This path uses one vertex from $A \setminus \{x\}$, leaving behind $\text{blue-}d_A(x) - 1 = (n + 1) - 1 = n$ unused vertices in $A \setminus \{x\}$; thus, we can use $x$ as the center of our $n$-star (with spokes in $A \setminus \{a_1\}$). 
 The path uses two vertices in $B \setminus \{y\}$, leaving behind $\text{blue-}d_B(y) - 2 = (m + 2) - 2 = m$ unused vertices, so we can use $y$ as the center of our $m$-star (with spokes in $B \setminus \{b_1, b_2\}$). This completes a blue $S_{6}(n, m)$.

 Therefore, we may assume that either $|A_1| \leq 1$ or $|B_1| \leq 1$, or both. Observe that $|A_1| \leq 1$ implies that $A_1 = \emptyset$ (this follows directly from the definition of $A_1$).  Similarly, $|B_1| \leq 1$ implies that $B_1 = \emptyset$. 

 Suppose that $B_1$ is empty, meaning that $B$ induces a blue complete graph.  By Lemma \ref{lem:faudree_schelp}, there is a monochromatic path on four vertices alternating between $A \setminus \{ x \}$ and $B\setminus\{ y \}$. If this path is red, then by Lemma \ref{lem:appendix_helper_1}, it completes a red $S_6(n,m)$. Hence, we may assume that the path---say, $x_1 - y_1 - x_2 -y_2$---is blue, with $x_1, x_2 \in A \backslash\{x \}$ and $y_1, y_2 \in B \backslash \{ y \}$. If there now exists a blue edge from $x$ to some vertex $w \in W \cup \{ v \}$, then the path $$x - x_1 -y_1 -x_2 -y_2 -y$$ completes a blue $S_6(n,m)$, with $x$ as the center of our $n$-star (with $n - 1$ spokes in $A \setminus \{x_1, x_2\}$ and $w$ as the last spoke), and $y$ as the center of our $m$-star (with spokes in $B \setminus \{y_1, y_2\}$). Therefore, we may assume that $x$ is red-complete to $W \cup \{ v\}$.
 
 Next, choose any three vertices $b_1, b_2, b_3 \in B \setminus \{ y_1 \}$, and consider the blue path $$x - x_1 - y_1 - b_1 - b_2 - b_3.$$ (Note that $b_1$, $b_2$, and $b_3$ are blue-adjacent in the path since $B$ induces a blue complete graph.) If $b_3$ is blue-adjacent to any vertex $w \in W \cup \{ u \}$, then the path above completes a blue $S_6(n,m)$. In particular, our path uses  one vertex in $A \setminus \{x\}$, leaving behind $\text{blue-}d_A(x) - 1 = (n + 1) - 1 = n$ unused vertices in $A \setminus \{x\}$; thus, we can use $x$ as the center of our $n$-star (with spokes in $A \setminus \{x, x_1\}$). Similarly, the path uses three vertices in $B \setminus \{b_3\}$, leaving behind $\text{blue-}d_B(b_3) - 3 = (m + 2) - 3 = m - 1$ unused vertices in $B \setminus \{y_1, b_1, b_2\}$, so we can use $b_3$ as the center of our $m$-star (with $m - 1$ spokes in $B \setminus \{y_1, b_1, b_2\}$ and $w$ as the last spoke). Hence, we may assume that $b_3$ is red-adjacent to all of $W \cup \{ u \}$. But our choice of $b_3 \in B \backslash \{ y_1 \}$ was arbitrary, so we may in fact assume that there is a red complete bipartite graph between $B \backslash \{ y_1 \}$ and $W \cup \{ u \}$.
 
But now, consider the red path \[u - x - w_1 - b_1 -w_2 -b_2,\] for some $w_1, w_2 \in W$. This path uses one vertex in $A$, leaving behind $\text{red-}d_A(u) - 1 = (n + 2) - 1 = n + 1$ vertices in $A$; thus, we  can use $u$ as the center of our $n$-star (with spokes in $A \setminus \{x\}$). Similarly, the path uses two vertices in $W$, leaving behind $\text{red-}d_W(b_2) - 2 = |W| - 2 \geq  m - 1$ vertices, so we can use $b_2$ as the center of our $m$-star (with $m - 1$ spokes in $W \setminus \{w_1, w_2\}$ and $v$ as the last spoke). Therefore, if $B_1 = \emptyset$, then our $2$-coloring of $K_r$ admits a monochromatic $S_6(n, m)$.

 On the other hand, suppose that  $A_1 = \emptyset$ and $B_1 \neq \emptyset$. By Lemma \ref{lem:faudree_schelp}, there is a monochromatic path on four vertices alternating between $A \backslash \{x\}$ and $B \backslash\{ y \}$; by Lemma \ref{lem:appendix_helper_1}, we may assume this path---say, $a_1 -b_1-a_2-b_2$---is blue, with $a_1, a_2 \in A \setminus \{x\}$ and $b_1, b_2 \in B \setminus \{y\}$.
 If $x$ is blue-adjacent to some vertex $w \in W \cup \{ v \}$, we can complete a blue $S_6(n,m)$ by connecting $x$ to $a_1$ and $y$ to $b_2$, using $y$ as the center of our $m$-star (with spokes in $B \setminus \{y, b_1, b_2\}$), and using $x$ as the center of our $n$-star (with $n - 1$ spokes in $A \setminus \{a_1, a_2\}$ and $w$ as the last spoke). Therefore, we may assume that $x$ is red-adjacent to all of $W \cup \{ v \}$.
 
 Finally, let $a \in A$. Since $A_1 = \emptyset$, we have that $\text{blue-}d_A(a) = n + 1$. By similar reasoning to before, we may assume there is a blue path $a_1 - b_1 - a_2 - b_2$ with $a_1, a_2 \in A \setminus \{a\}$ and $b_1, b_2 \in B \setminus \{y\}$. If there is a blue edge between $a$ and any vertex $w \in W \cup \{v\}$, then we can complete a blue $S_{6}(n, m)$ with the path $a - a_1 - b_1 - a_2 - b_2 - y$, using $a$ as the center of our $n$-star (with $n - 1$ spokes in $A \setminus \{a, a_1, a_2\}$ and $w$ as the last spoke) and $y$ as the center of our $m$-star (with spokes in $B \setminus \{y, b_1, b_2\}$). Thus, we can assume that there are no blue edges between $a$ and $W \cup \{v\}$, for all $a \in A$; that is, there is a red complete bipartite graph between $A$ and $W \cup \{v\}$. By Lemma \ref{lem:appendix_helper_3}, there must also exist a red edge $\{a^\prime, b\}$ between $A$ and $B$, so we can form the path
 \[u- a^{\prime\prime} - w^\prime - a^\prime - b - v,\]
 where $a^{\prime\prime} \in A$ (with $a^{\prime\prime} \neq a^\prime$) and $w^\prime \in W$. However, this completes a red $S_{6}(n, m)$, with $u$ as the center of our $n$-star (with spokes in $A \setminus \{a^\prime, a^{\prime\prime}\}$) and $v$ as the center of our $m$-star (with spokes in $B \setminus \{b\}$).

In all three cases of the sizes of $A_1$ and $B_1$, we find that our $2$-coloring of $K_r$ admits a monochromatic $S_6(n, m)$, as desired.
\end{proof}

\begin{lemma}
    If $A_0 \neq \emptyset$ and $B_0 = \emptyset$, then our 2-coloring of $K_r$ admits a monochromatic $S_6(n,m)$.
\end{lemma}

\begin{proof}
    Suppose there exists some vertex $z \in B$ with $\text{blue-}d_B(z) = m+1$ and thus $\text{red-}d_B(z) = 1$; let $z^\prime$ be the unique red neighbor of $z$ in $B$. Fix any vertex $x \in A_0$.
    
    By Lemma \ref{lem:faudree_schelp}, there is a monochromatic path on four vertices alternating between $A \setminus \{ x \}$ and $B\setminus\{ z, z^\prime \}$. If this path is red, then by Lemma \ref{lem:appendix_helper_1}, it completes a red $S_6(n,m)$. Thus, we may assume that the path---say, $a_1 - b_1 - a_2 -b_2$---is blue, with $a_1, a_2 \in A \setminus \{x \}$ and $b_1, b_2 \in B \setminus \{ z, z^\prime \}$.  Appending $x$ and $z$ to the ends of this path yields \[x-a_1-b_1-a_2-b_2-z.\] Observe that $x$ is blue-adjacent to $\text{blue-}d_A(x) - 2 = (n + 1) - 2 = n-1$ unused vertices in $A$, and $z$ is blue-adjacent to $\text{blue-}d_B(z) - 2 = (m + 1) - 2 = m-1$ unused vertices in $B$. Therefore, if $x$ and $z$ are blue-adjacent to two distinct respective vertices in $W$, then we can complete a blue $S_6(n, m)$.
    
    If $x$ and $z$ have a common red neighbor $w \in W$, then the following path completes a red $S_6(n, m)$:
    \[u - x - w - z - z^\prime - v,\]
    where $u$ is the center of our $n$-star (with spokes in $A \setminus \{x\}$) and $v$ is the center of our $m$-star (with spokes in $B \setminus \{z, z^\prime\}$). Hence, we can assume that $x$ and $z$ do not share any red neighbors in $W$, and by our prior reasoning, that they are not blue-adjacent to two distinct respective vertices in $W$. It follows that either $x$ is red-complete and $z$ is blue-complete to $W$, or $x$ is blue-complete and $z$ is red-complete to $W$.

    Suppose that $x$ is red-complete and $z$ is blue-complete to $W$.   If there is a red edge between some vertex $w \in W$ and $b \in B$, then the path 
    \[u - x - w - b - b^\prime - v\]
    completes a red $S_6(n, m)$, with $u$ as the center of our $n$-star (with spokes in $A \setminus \{x\}$) and $v$ as the center of our $m$-star (with spokes in $B \setminus \{b, b^\prime\}$). Note that such a vertex $b^\prime \in B$ with $\{b, b^\prime\}$ being red must exist, since $B_0 = \emptyset$. Consequently, we may assume that there is a blue complete bipartite graph between $B$ and $W$.

    If there is a red complete bipartite graph between $A \setminus \{x\}$ and $B \setminus \{z\}$, then it contains a red path on four vertices; by Lemma \ref{lem:appendix_helper_1}, our $2$-coloring of $K_r$ thus admits a red $S_6(n, m)$. Therefore, we may assume that there is a blue edge between some vertex $a \in A \setminus \{x\}$ and $b \in B \setminus \{z\}$. But we can now form the blue path
    \[x - a - b - w - b^\prime - z,\]
    for any vertex $w \in W$ and blue neighbor $b^\prime \in B \setminus \{b\}$ of $z$. This path completes a blue $S_{6}(n, m)$, with $x$ as the center of our $n$-star (with spokes in $A\setminus \{a\}$) and $z$ as the center of our $m$-star (with $m - 1$ spokes in $B \setminus \{b, b^\prime\}$ and any vertex in $W \setminus \{w\}$ as the last spoke).

    On the other hand, suppose that $x$ is blue-complete and $z$ is red-complete to $W$; the reasoning in this case is similar. If there is a red edge between some vertex $a \in A$ and $w \in W$, then the path
    \[u - a - w - z - z^\prime - v\]
    completes a red $S_6(n, m)$, with $u$ as the center of our $n$-star (with spokes in $A \setminus \{a\}$) and $v$ as the center of our $m$-star (with spokes in $B \setminus \{z, z^\prime\}$). Therefore, we may assume that there is a blue complete bipartite graph between $A$ and $W$.

    If there is a red complete bipartite graph between $A$ and $B_B(z)$ (the blue neighborhood  of $z$ in the set $B$), then it contains a red path on four vertices; by Lemma \ref{lem:appendix_helper_1}, our $2$-coloring of $K_r$ admits a monochromatic $S_6(n, m)$. Hence, we can assume that there is a blue edge between some vertex $a \in A$ and blue neighbor $b \in B$ of $z$. But now, consider the blue path
    \[w_1 - a^\prime - w_2 - a - b - z,\]
    where $w_1, w_2 \in W$ (with $w_1 \neq w_2$) and $a^\prime \in A \setminus \{a\}$. This path completes a blue $S_6(n, m)$, with $w_1$ as the center of our $n$-star (with spokes in $A \setminus \{a, a^\prime\}$) and $z$ as the center of our $m$-star (with spokes in $B \setminus \{b\}$).

    In either case, our $2$-coloring admits a monochromatic $S_6(n, m)$. As a result, we may assume that no such vertex $z \in B$ with $\text{blue-}d_B(z) = m + 1$ exists; hence, $\text{red-}d_B(b) \geq 2$ for all $b \in B$. It follows from Lemma \ref{lem:min_path_degree} that every vertex $b \in B$ is the endpoint of a red path on three vertices. By Lemma \ref{lem:appendix_helper_3}, we can also assume that there exists a red edge between some vertex $a \in A$ and $b \in B$, so we can form the red path $b-b^\prime - b^{\prime\prime}$ in $B$ (with $b$, $b^\prime$, and $b^{\prime\prime}$ being distinct). But now, the red path \[u - a - b - b^\prime - b^{\prime\prime}- v,\] completes a red $S_6(n,m)$, with $u$ as the center of our $n$-star (with spokes in $A \setminus \{a\}$) and $v$ as the center of our $m$-star (with spokes in $B \setminus \{b, b^{\prime}, b^{\prime\prime}\}$), as desired.
\end{proof}

\begin{lemma}
    If $A_0 = \emptyset$ and $B_0 \neq \emptyset$, then the 2-coloring admits a monochromatic $S_6(n,m)$.
\end{lemma}

\begin{proof}
    Fix some vertex $y \in B_0$. By the definition of $B_0$, we have that $\text{blue-}d_B(y) = m+ 2$. Now, suppose there exists a vertex $w \in W$ such that the edges $\{w, a\}$ and $\{w, b\}$ are red, for some $a \in A$ and $b \in B$. Let $a^\prime$ be any red neighbor of $a$ in $A$ (note that such a vertex $a' \in R_A(a)$  must exist, since $A_0 = \emptyset$). Then, the red path \[u - a' - a - w - b - v\] completes a red $S_6(n,m)$, with $u$ as the center of our $n$-star (with spokes in $A \setminus \{a, a^\prime\}$) and $v$ as the center of our $m$-star (with spokes in $B \setminus \{b\}$). Thus, we may assume that no such vertex $w \in W$ exists that is both red-adjacent to a vertex in $A$ and a vertex in $B$. Hence, every vertex $w \in W$ is either blue-complete to $A$, blue-complete to $B$, or both. As such, define
    \begin{equation*}
            W_A = \{ w \in W : A \subseteq B(w)\}\quad\text{and}\quad W_B = \{ w \in W : B \subseteq B(w)\}.\\
    \end{equation*}
    That is, $W_A$ and $W_B$ consist of the vertices in $W$ that are blue-complete to $A$ and $B$, respectively. Note that $W = W_A \cup W_B$.

    Suppose that $W_A$ and $W_B$ are both nonempty. Then, there exists a vertex $w \in W_A$ and $w^\prime \in W_B$, with $w \neq w^\prime$.
    By Lemma \ref{lem:faudree_schelp}, there is a monochromatic path on four vertices alternating between $A$ and $B$. If this path is red, then by Lemma \ref{lem:appendix_helper_1}, it completes a red $S_6(n,m)$. Hence, we may assume that the path---say, $a_1 - b_1 - a_2 -b_2$---is blue, with $a_1, a_2 \in A$ and $b_1, b_2 \in B$. Appending $w$ and $w^\prime$ to the ends of this path gives us \[w-a_1-b_1-a_2-b_2-w',\] which completes a blue $S_6(n,m)$, with $w$ as the center of our $n$-star (with spokes in $A \setminus \{a_1, a_2\}$) and $w^\prime$ as the center of our $m$-star (with spokes in $B \setminus \{b_1, b_2\}$).
    
    Hence, we may assume that either $W_A$ or $W_B$ is empty; let us first suppose that $W_B = \emptyset$. Since $W = W_A \cup W_B$ and thus $W = W_A$, there is a blue complete bipartite graph between $A$ and $W$. If there is a red complete bipartite graph between $A$ and $B \setminus \{y\}$, then it contains a red path on four vertices; by Lemma \ref{lem:appendix_helper_1}, this gives us a red $S_6(n, m)$.
    Therefore, we can assume that there is a blue edge between some vertex $a \in A$ and $b \in B \setminus \{y\}$. But now, choose any vertex $a' \in A \setminus \{a\}$, and any distinct $w, w' \in W$. The path \[w' - a' - w - a - b -y\] completes a blue $S_6(n,m)$, with $w^\prime$ as the center of our $n$-star (with spokes in $A \setminus \{ a, a^\prime \}$) and $y$ as the center of our $m$-star (with spokes in $B \setminus \{ b \}$).

    Finally, suppose that $W_A$ is empty. Since $W$ now equals $W_B$, there is a blue complete bipartite graph between $B$ and $W$. First, let us assume there is a red edge $\{b_1, b_2\}$ in $B\setminus\{ y \}$. If $b_1$ is red-adjacent to some other vertex $a \in A$, then the path \[u-a'-a-b_1-b_2-v\] completes a red $S_6(n,m)$. Note that such a vertex $a^\prime \in A$ with $\{a, a^\prime\}$ being red must exist, since $A_0 = \emptyset$. Therefore, we may assume that $b_1$ is blue-complete to $A$, and by similar reasoning, so is $b_2$. But now, consider the blue path \[
    b_1 - w_1 - b - w_2 - b^\prime - y,\]
    where $b,b' \in B \setminus\{ y,b_1 \}$ (with $b \neq b^\prime$) and $w_1, w_2 \in W$ (with $w_1 \neq w_2$). This path completes a blue $S_6(n,m)$, with $b_1$ as the center of our $n$-star (with spokes in $A$) and $y$ as the center of our $m$-star (with $m - 1$ spokes in $B \setminus \{b_1, b, b^\prime\}$ and any vertex in $W \setminus \{w_1, w_2\}$ as the last spoke). Consequently, we may assume that no such red edge $\{b_1, b_2\}$ exists inside $B \setminus \{y\}$, and therefore $B$ (including $y$) induces a blue complete graph.

    Suppose there exists a vertex $b_0 \in B$ with $\text{blue-}d_A(b_0) \geq n$. Then, choose three distinct vertices $b_1, b_2, b_3 \in B \setminus \{b_0\}$ and distinct $w_1, w_2 \in W$, and consider the path
    \[b_0 - w_1 - b_1 - w_2 - b_2 - b_3.\]
    This path completes a blue $S_6(n, m)$, with $b_0$ as the center of our $n$-star (with spokes in $A$) and $b_3$ as the center of our $m$-star (with $m - 1$ spokes in $B \setminus \{b_0, b_1, b_2\}$ and any vertex in $W \setminus \{w_1, w_2\}$ as the last spoke).
    
    As a result, we may assume that $\text{blue-}d_A(b) \leq n - 1$ and thus $\text{red-}d_A(b) \geq (n + 2) - (n - 1) = 3$ for all $b \in B$. If any two vertices $b, b^\prime \in B$ shared a red neighbor $a \in A$, then we can form the red path on four vertices
    \[a^\prime - b - a - b^\prime,\]
    where $a^\prime$ is any red neighbor of $b$ in $A \setminus \{a\}$. By Lemma \ref{lem:appendix_helper_1}, our $2$-coloring of $K_r$ admits a red $S_6(n, m)$. Therefore, we can also assume that the red $A$-neighborhoods of vertices in $B$ are pairwise disjoint. It follows that $3|B| \leq |A|$, so \[3(m + 3) \leq n + 2.\]
    Rearranging, we find that $n \geq 3m + 7$ and so $|W|= n - m- 2$. Since the red $A$-neighborhoods of vertices in $B$ are pairwise disjoint, we have that \begin{align*}
        \text{red-}d_A(b) \leq |A| - \sum_{b' \in B \setminus\{ b \} } \text{red-}d_A(b')  \leq |A| - 3(|B| - 1) = n - 3m -4
    \end{align*}
    for all $b \in B$; hence, $\text{blue-}d_A(b) \geq (n + 2) - (n - 3m -4) = 3m + 6$. Once again, consider the path \[b_0 - w_1 - b_1 - w_2 - b_2 - b_3.\] Choose some vertex $w_3 \in W \setminus \{w_1, w_2\}$. We can now use $b_3$ as the center of our $m$-star (with $m - 1$ spokes in $B$ and $w_3$ as the last spoke). This leaves behind $|W| - 3 = n - m -5$ blue neighbors of $b_0$ in $W$ (recall that $b_0$ is blue-complete to $W$). Furthermore, we have that $\text{blue-}d_A(b_0) \geq 3m + 6$, giving us
    \[(3m + 6) + (n - m - 5) = n + 2m + 1 \geq n\]
    unused vertices to complete the $n$-star at $b_0$ (with spokes in $(A \cup W) \setminus \{w_1, w_2, w_3\}$). This completes a blue $S_6(n, m)$, as desired.
\end{proof}

\newpage

\bibliographystyle{amsplain}
\bibliography{references}

\end{document}